\documentclass[11pt,twoside,leqno]{aomamlt2e} 
  \pageno{649}
\received{March 12, 2004}
 \renewcommand{\titleheadline}[1]{\def\one{#1}\ifx\one\empty\else
\gdef\thetitle{{\frenchspacing%
\let\\ \relax\scriptsize{#1}}}}

\newif\ifshort

\let\shorttitle\titleheadline
\newcounter{notes}
\newenvironment{my}
{\begin{list}{{\rm \arabic{notes}.}}{\usecounter{notes}%
\setlength\parsep{.25ex plus.1ex minus.1ex}
\setlength\itemsep{.25ex  plus.1ex minus.1ex}
\setlength\topsep{3pt  minus1pt}}}
{\end{list}}

\newcommand\lder{H{\hskip1.5pt\rm \"{\hskip-7pt\it o}}lder}
\newcommand{\ho}{{\hskip1.5pt\rm \"{\hskip-7pt \it o}}}
\newcommand\hfq{\hfill\qed}
 
\def\theequation{\theSubsec.\arabic{equation}}
 \renewcommand{\thesection}{\arabic{section}}
\makeatletter
\renewcommand{\@seccntformat}[1]{\csname
the#1\endcsname.\hspace{0.5em}\setcounter{Subsec}{0}\setcounter{Subsubsec}{0}}\makeatother

\newtheorem{theorem}{Theorem}[Subsec]
\newtheorem{proposition}[theorem]{Proposition}
\newtheorem{lemma}[theorem]{Lemma}
\newtheorem{corollary}[theorem]{Corollary}

\theorembodyfont{\upshape}

\newtheorem{definition}[theorem]{{\it Definition}}

\def\dbar{\bar\partial}
\def\bx{\square}

\def\C{{\mathbb C}}
\def\CalC{\mathcal C}
\def\CC{\mathcal C}

\def\H{\mathcal H}
\def\I{\vartheta}
\def\K{\mathcal K}

\def\N{\mathcal M}
\def\N{\mathcal N}

\def\R{{\mathbb R}}
\def\S{{\mathcal S}}
\def\T{{\mathcal T}}

\def \IM {\Im\text{\rm m}}
\def\norm{\big\vert\big\vert}

\def\da{d_1(p_1,q_1)}
\def\db{d_2(p_2,q_2)}
\def\Va{V_1(p_1,q_1)}
\def\Vb{V_2(p_2,q_2)}

\def\pso{(s,p_{1},q_{1})}
\def\pst{(s,p_{2},q_{2})}

\def\sh{{\#}}
\def\m{{(m)}}
\def\bW{{\overline W}}

\begin{document}
\currannalsline{164}{2006} 

 \title{The $\dbar_b$-complex on\\  decoupled
boundaries in $\C^n$}

 \acknowledgements{Research supported in part by grants from the National Science Foundation.}
\twoauthors{Alexander Nagel}{Elias M. Stein}

 \institution{University of Wisconsin,
Madison, WI\\
\email{nagel@math.wisc.edu}\\
\vglue-9pt
 Princeton University,
 Princeton, NJ\\
\email{stein@math.princeton.edu}
}

 \shorttitle{$\overline\partial_b$ \scshape\uppercase{on decoupled boundaries}}

\newcommand{\namelistlabel}[1] {\mbox{#1}\hfil}
\newenvironment{namelist}[1]{%
\begin{list}{}
{\let\makelabel\namelistlabel \settowidth{\labelwidth}{#1}
\setlength{\leftmargin}{1.1\labelwidth}}
}{%
\end{list}}

\vskip6pt
\centerline{\bf Contents}
\def\sni#1{\smallbreak\noindent{#1}.\hskip5pt}
\sni{1} Introduction
\sni{2} Definitions and statement of results
\sni{3} Geometry and analysis on $M_j$ and on $M_1\times \cdots\times M_n$
\sni{4} Relative fundamental solutions for $\Box_b$ on $M_1\times \cdots\times M_n$
\sni{5} Transference from $M_1\times \dots\times M_n$ to $M$ and $L^p$ regularity of $\mathcal{K}$
\sni{6} Pseudo-metrics on $M$
\sni{7} Differential inequalities for the relative fundamental solution $K$
\sni{8} H\"older regularity for $\mathcal{K}$
\sni{9} Examples
\smallbreak\noindent References

\section{Introduction}

The purpose of this paper is to prove optimal estimates for solutions of the Kohn-Laplacian for certain classes of model domains in several complex variables. This will be achieved by applying a type of singular integral operator whose novel features (related to product theory and flag kernels) differ essentially from the more standard Calder\'{o}n-Zygmund operators that have been used in these problems hitherto.

\Subsec{Background}\quad\label{Background}
We consider the Kohn-Laplacian on $q$-forms, $\Box^{(q)}_{b} = \Box_b = \bar{\partial}_b
\bar{\partial}^\ast_b + \bar{\partial}_b^\ast \bar{\partial}_b$, defined on the boundary
$M = \partial \Omega$ of a smooth bounded pseudo-convex domain $\Omega \subset \mathbb{C}^n$. Our objective is the study of the (relative) inverse operator $\mathcal{K}$ and the corresponding Szeg\"{o} projection $\S$ (when it exists), which satisfy $\Box_b
\, \mathcal{K} = \ \mathcal{K} \, \Box_b = I-\S$. By definition $\S$ is the orthogonal projection on the $L^2$ null-space of $\Box_b$.

In formulating the questions of regularity pertaining to the above, it is useful to recall Fefferman's hierarchy \cite{Fefferman95} of the levels of understanding of the problem, which we rephrase as follows:

\begin{enumerate}
\item Proof of $C^\infty$ regularity.

\item Derivation of optimal $L^p$, H\"older, and Sobolev-space estimates of solutions.
 
\item Analysis of singularities of the distribution kernels of the operators $\mathcal{K}$ and $S$ and derivation of the estimates in (2) from a corresponding theory of singular integrals.
\end{enumerate}

Now as far as the $C^\infty$ regularity is concerned, this has been resolved in the general situation where an appropriate ``finite-type'' condition holds (at least for the closely connected $\bar{\partial}$-Neumann problem) by the work of Kohn \cite{Ko72}, \cite{Kohn1979}, Catlin \cite{Catlin1983}, \cite{Catlin1987}, and D'Angelo \cite{D'Angelo1982}. However, the more refined results of (2) and (3) have been obtained only in a more restrictive setting. This was carried out in a series of developments beginning with the work of Folland and Stein (\cite{FoSt74}) in the strongly pseudo-convex case, and in later works of, among others, Christ (\cite{Christ91a}, \cite{Christ88}), Fefferman and Kohn (\cite{FefKoh88}), Kohn (\cite{Kohn85}), McNeal (\cite{McNeal89}), Nagel, Rosay, Stein, and Wainger (\cite{NaRoStWa89}), and Rothschild and Stein (\cite{RoSt}). This culminated in the work of Koenig (\cite{Koenig01}) on finite type domains whose Levi-form has comparable eigenvalues.

At the base of these results is a version of the Calder\'{o}n-Zygmund theory for the following class of singular integrals: One considers operators $\mathcal{T}$ of the form $\mathcal{T} ( f ) ( x ) = {\int} T ( x , y ) \, f (y ) dy$ whose kernels $T ( x , y )$ are distributions that are smooth away
from the diagonal, that satisfy the characteristic size estimates $\left| T ( x , y ) \right| \, \lesssim \, d(x,y)^a\,V(x,y)^{-1}$, and that satisfy corresponding differential inequalities and cancellation properties. Here $d ( x , y )$ is the control metric determined by the vector fields which are the real and imaginary parts of the tangential Cauchy-Riemann operators, and $V (x , y)$ denotes the volume of the ball centered at $x$ of radius $d(x,y)$. It can be shown that the relative fundamental solution $\K$ and the Szeg\"{o} projection $\S$ are of this type, with $a = 2$ for $\mathcal{K}$, and $a = 0$ for $\S$. As a result, one obtains for these operators maximal sub-elliptic estimates in $L^p$, \textit{etc}.

Unfortunately, while highly satisfactory, the above framework with a natural metric controlling all estimates
cannot carry over in general. In fact, in more general circumstances there seem to arise a number of
inequivalent metrics that control different aspects of the problem.  This appears to be connected with earlier
observations of Derridj \cite{Derr78} and Rothschild \cite{Ro80} that maximal  sub-ellipticity can hold only
if the eigenvalues of the  Levi-form are comparable. It is the purpose of this paper to make progress in the
resolution of problems (2) and (3) in an illustrative model case - that of decoupled domains.

\Subsec{A special case}\label{SSSpecialCase} To get a better grasp of these problems and 
the results we obtain, we take a closer look at the special case 
of a decoupled domain where 
$\Omega = \{ z \in \mathbb{C}^3$: $\IM [z_3] > | z_1 |^n + |z_2 |^m \}$, with $n,\, m$ even integers. Then $M = \partial \Omega$ can be identified with $ \{ ( z , t ) \in \mathbb{C}^2 \times \mathbb{R}, \, z = ( z_1 , z_2) \}$, and 
$$
\bar{Z}_1 = \frac{\partial}{\partial \bar{z}_1} - i\,\frac{n}{2} |z_1 |^{n-2} \, z_1 \frac{\partial}{\partial t}, \qquad \bar{Z}_2 = \frac{\partial}{\partial \bar{z}_2} - i\,\frac{m}{2} | z_2|^{m-2} z_2 \frac{\partial}{\partial t}
$$
form a basis for the tangential Cauchy-Riemann vector fields. The eigenvalues $\lambda_1$, $\lambda_2$ of the Levi-form at a point $(z_1,z_2,t)$ are essentially $| z_1|^{n-2}$ and $|z_2|^{m-2}$, and are not comparable. With $\bar{Z}_j = \frac{1}{2} ( X_j + i Y_j)$, we can consider $d_\Sigma$, the control metric defined by $X_1 , Y_1 , X_2 , Y_2$. 

However, the above domain is also convex, so that  there is another natural metric, which reflects the
``flatness'' of the boundary in different complex directions, the ``Szeg\"{o} metric'' $d_S$; (see McNeal
\cite{McN94},
\cite{McN94b}, and Bruna, Nagel and Wainger \cite{BrNaSt88} for a real analogue). In our special case, if
$n \leq m$, when we measure the distance of the point $p= (z_1, z_2, t)$ from the origin $0$ we have:
\begin{eqnarray*}
\noalign{\vskip-6pt}
d_\Sigma(0,p) &\approx& |z_1|\phantom{^m}+|z_2|\phantom{^n}+|t|^{1/m} ;\\
d_S(0,p) &\approx& |z_1|^m + |z_2|^n +|t|.
\end{eqnarray*}
Note that $d_{S}(0,p)^{1/m}\approx |z_{1}|+|z_{2}|^{n/m}+|t|^{1/m}$, and this is not equivalent to $d_{\Sigma}(0.p)$ if $n \neq m$. Thus these metrics, or powers of these metrics are in general not equivalent.

Now $d_\Sigma$ controls the inverse of the sub-Laplacian $\mathcal{L} = - \frac{1}{2} \,{\sum^{2}_{i = 1}} \, ( Z_i \bar{Z}_i \, + \, \bar{Z}_i Z_i)$, while $d_S$ controls the Szeg\"{o} kernel (the orthogonal projection on the null-space of the operator $- {\sum^{2}_{i = 1}} Z_i \bar{Z}_i$), and some mixture of $d_\Sigma$ and $d_S$ arises in the fundamental solution of the operator $\Box_b = -(Z_1 \bar{Z}_1 + \bar{Z}_2 Z_2) = \Box^1_b \, + \, \Box^2_b$, which is essentially the Kohn-Laplacian acting on $1$-forms.

With this we can state a part of our main result obtained below, formulated in this special case, as follows:

\vskip6pt {\scshape Theorem.} 
 {\it There is an operator $\mathcal{K}$ so that{\rm ,} when it is
 applied to smooth functions with compact support{\rm ,} there is 
the identity $\mathcal{K} \, \Box_b \, = \, \Box_b \, \mathcal{K} \, = I$. Moreover}
\begin{my}

\item[{\rm (a)}] {\it The four operators $Z_1\,\bar Z_1\,\K =
 \bx_b^1\, \K${\rm ,} \,\,$\bar Z_2\, Z_2\,\K = \bx_b^2\,\K${\rm ,} 
 $\bar Z_1\,\bar Z_1\,\K${\rm ,} and $Z_2\,Z_2\,\K$ are
bounded on $L^p(M)$ for $1 < p < \infty$.}
 
\item[{\rm (b)}] {\it Let $B_1,\,B_2$ be bounded functions on $M${\rm ,}
 and suppose there are constants $C_1,\,C_2$ so that}
\begin{eqnarray*}
\noalign{\vskip-6pt}
\lambda_1(z_1)\,B_1(z_1,z_2,t) &\leq& C_1\,\lambda_2(z_2);\\
\lambda_2(z_2)\,B_2(z_1,z_2,t) &\leq &C_2\,\lambda_2(z_1).
\end{eqnarray*}
{\it Then the two operators $B_1\,\bar Z_1\,Z_1\,\K = B_1\,\overline\bx_b^1\,\K$
 and $B_2\, Z_2\,\bar Z_2\,\K = B_2\,\overline\bx_b^2\,\K$ are bounded on $L^p(M)$ for $1 < p < \infty$.
Here $\lambda_1(z) = |z|^{m-2}$ and $\lambda_2(z) = |z|^{n-2}$ are the eigenvalues of the Levi form.}

\item[{\rm (c)}] {\it Let $B_1,\,B_2$ be bounded functions on $M${\rm ,}
 and suppose there are constants $C_1,\,C_2$ so that}
\begin{align*}
B_1(z_1,z_2,t) &\leq C_1\,\lambda_2(z_2);\\
B_2(z_1,z_2,t) &\leq C_2\,\lambda_2(z_1).
\end{align*}
{\it Then the two operators $B_1\, Z_1\,Z_1\,\K $ and $B_2\, \bar Z_2\,\bar Z_2\,\K$ 
are bounded on $L^p(M)$ for $1 < p < \infty$.}

\item[{\rm (d)}] {\it $\mathcal{K}$ maps $L^\infty(M)$ to the isotropic
\lder\ space $\Lambda^\alpha(M)${\rm ,} where} $$\alpha \, =
\min\left\{\frac{2}{n},\frac{2}{m}\right\}.$$ 
\end{my}
 
\vskip8pt

 The conclusion (b) is part of the optimal substitute for maximal sub-ellipticity that holds in this case.

\Subsec{Methods used} \label{MethodsUsed}
To describe the methods used we continue with the case considered above. We begin by considering
separately the component domains
\begin{align*}
M_1 &= \left\{(z_1,w_1)\in \C^2\,\big\vert\,\Im[w_1] = |z_1|^n\right\} \simeq \{(z_1,t_1) \in \C\times\R\}, \quad \text{and}\\
M_2 &= \left\{(z_2,w_2)\in \C^2\,\big\vert\,\Im[w_2] = |z_2|^m\right\} \simeq \{(z_2,t_2) \in \C\times\R\}.
\end{align*}
We denote by $\widetilde{M}$ the Cartesian product $ M_1 \times M_2$ and we let $\pi$ be the projection of $\widetilde{M}$ to $M$ given by $\pi: ( z_1 ,t_1) \times (z_2 , t_2) \, \to \, (z_1 , z_2 , t_1 + t_2 )$.

The idea is to deduce the results about regularity of $\Box_b$ on $M$ from corresponding results on $\widetilde{M}$. Moreover, passing to the product allows one to consider various combinations of the separate metrics on each factor of $\widetilde{M}$, which in effect account for the different metrics on $M$. Our analysis proceeds as follows.
\begin{my}

\item[(1)] \emph{Analysis on each $M_j$}:  
  Here the key point is the use of the nonhypoelliptic ``heat'' semi-group $e^{-s \Box_j}$ on $M_j$ where
$\Box_1 = Z_1 \bar{Z}_1$, $\Box_2 = \bar{Z}_2 Z_2$. (The needed estimates for this semi-group were
obtained in \cite{NaSt00}.) For later purposes one observes that if\ \begin{equation*}
\mathcal{K}_j = \displaystyle{\int^{\infty}_{0}} ( e^{-s \Box_j} - S_j ) ds,
\end{equation*}
where $S_j$ is the orthogonal projection on the null-space of $\Box^j$, then $\mathcal{K}_j \Box^j = \Box^j \mathcal{K}_j = I - S_j$.

\item[(2)] \emph{Results on the product $\widetilde{M} = M_1 \times M_2$}:   In finding a (relative)
inverse for $\Box_1 + \Box_2$ on $\widetilde{M}$ one considers 
$$
\widetilde{\mathcal{K}} \, = \, \displaystyle{\int^{\infty}_{0}} (e^{-s ( \Box_1 + \Box_2)} - S_1 \otimes S_2 ) d s 
$$
 and also a substitute version
$$
\widetilde{\mathcal{N}} = \displaystyle{\int^{\infty}_{0}} ( e^{-s \Box_1} - S_1 ) \, \otimes \, (e^{- s \Box_2} - S_2 ) ds .
$$
 Now $\widetilde{\N}$ is more tractable than $\widetilde{\K}$ since any second order derivative in $Z_j$
and $\bar Z_j$ of $\widetilde \N$ turns out to be a product-type singular integral on $M_1 \times M_2$. For
such singular integrals an $L^p$ theory has been worked out in \cite{NaSt00.3}. However,
$\widetilde{\mathcal{K}}$ is the desired relative inverse, since $( \Box_1 + \Box_2)
\widetilde{\mathcal{K}} = \widetilde{\mathcal{K}} ( \Box_1 + \Box_2 ) = I - S_1 \otimes S_2$; its
properties can ultimately be deduced from those of $\widetilde{N}$ because of the identity
$$
\widetilde{\mathcal{K}} = \widetilde{\mathcal{N}} \, + \, \mathcal{K}_1 \, \otimes \, S_2 \, + \, S_1 \, \otimes \, \mathcal{K}_2 \, .
$$

\item[(3)] \emph{Descent to $M$}: 
  The operators above on $\widetilde{M} = M_1 \times M_2$ are translation-invariant in the $t_1$ and $t_2$
variables. Each appropriate operator $T$ of this kind can be transferred by the projection $\pi$:
$\widetilde{M} \longrightarrow M$ to an operator $T^{\#}$ on $M$, via the identity
$$
T^{\#} ( f ) \, = \, J ( T ( f \, \circ \, \pi ))
$$
where $J ( F ) ( z_1,w_1,z_2,w_2,t ) \, = \, \displaystyle{\int^{\infty}_{- \infty}} \, F (z_1,w_1, t - s ,z_2,w_2, s ) \, ds$. This is then applied to $\widetilde{\mathcal{K}}$ to obtain $\mathcal{K} = ( \widetilde{\mathcal{K}} )^{\#}$, the inverse of $Z_1 \bar{Z}_1 + \bar{Z}_2
Z_2$ on~$M$.
\end{my}
 \vskip5pt

There is however a fundamental issue that arises at this point. Operators like $\widetilde{\mathcal{K}}$ and $\widetilde{\mathcal{N}}$ are not pseudo-local, because as product-like operators their kernels have singularities on the products of the diagonals of the $M_i$, and not just on the diagonal of $\widetilde{M}$. As a result the projections of such operators on $M$ are thus in general again not pseudo-local. Why then is the operator $\mathcal{K}$ pseudo-local? Connected with this is the question of obtaining the appropriate differential inequalities satisfied by the kernel of $\mathcal{K}$ away from diagonal.

The resolution of these problems is connected with the key idea of ``borrowing'', which allows one to pass from smoothness inherent in the $t_1$ (and $z_1$) variable to the $t_2$ (and $z_2$) variable, and vice-versa. This technique is used in several places below where it takes a number of different forms. A particularly transparent example is the identity
\begin{equation*}
\left( \frac{\partial}{\partial t_1} \, S_1 \, \otimes K_2\right)^{\#} \, = \, \left( S_1 \, \otimes \, \frac{\partial}{\partial t_2} \, K_2 \right)^{\#}
\end{equation*}
which is used in obtaining conclusion (b) of the theorem above.

\Subsec{Previous work}  Besides the results mentioned earlier which deal with the situation of comparable
eigenvalues of the Levi-form, several other situations have been previously studied. The case of a decoupled
domain in $\mathbb{C}^3$ with exactly one degenerate eigenvalue was dealt with in the paper of Machedon
\cite{Ma88} , where he also finds certain estimates for the fundamental solution which involve several
metrics. In addition, Fefferman, Kohn, and Machedon \cite{FefKohMa90} have obtained results on
H\"{o}lder regularity for $\Box_b$ on boundaries of diagonalizable domains (which is a larger class of
domains than we consider). In contrast, here we obtain sharp $L^p$ and H\"{o}lder estimates, and relevant
differential inequalities for the solving operators and Szeg\"{o} projections.

The general idea of ``lifting'' to a product (or ``simpler'' situation) is old, having already appeared in different forms in the study of the sub-Laplacian \cite{RoSt}, and in \cite{Ma88}. More recently it was used in \cite{MuRiSt95} to study certain operators on the Heisenberg group, and for $\Box_b$ on quadratic CR manifolds of higher-codimension in \cite{NaRiSt00}. The operators arising in \cite{MuRiSt95}, related to the boundary operator of the $\bar{\partial}$-Neumann problem for the ball, which occurred in \cite{PhongStein86}, already implicitly display the feature of the conflicting metrics which we have discussed above. There the kernels of the relevant operators arise as products of components that are homogeneous in different senses: the isotropic homogeneity reflecting the Euclidean metric, and the automorphic homogeneity of the Heisenberg group, reflecting the control metric.

\Subsec{Organization of the paper} Section 2 contains a review of background material and statements of
the main results of the paper. The needed aspects of the geometry and analysis of each of the factors $M_i$
and on their Cartesian product are set down in Section 3. Section 4 studies the various versions of the relative
fundamental solutions of $\Box_b$ on $\widetilde{M}$. This leads to $L^p$ results on $M$ via
transference, as is shown in Section 5. Section 6 deals with the various metrics on $M$ and the resulting
differential inequalities of the kernels are obtained in Section 7. In Section 8, we prove the H\"{o}lder
regularity of the solutions, and in Section 9 we give examples to show that our regularity results are optimal.

\section {Definitions and statement of results}
\vglue-12pt
\Subsec{Definitions}
A domain $\Omega \subset \C^{n+1}$ and its boundary $M$ are said to be \emph{decoupled} if there are
sub-harmonic, nonharmonic polynomials $P_j$
such that
\begin{equation}\label{E:1.0.2}
\begin{split}
\Omega &= \Big\{(z_1,\ldots,z_n,z_{n+1}) \in \C^{n+1}\,\Big\vert\,\IM[z_{n+1}] > \sum_{j=1}^n P_j(z_j)\Big\};\\
M &= \Big\{(z_1,\ldots,z_n, z_{n+1}) \in \C^{n+1}\,\Big\vert\,\IM[z_{n+1}] = \sum_{j=1}^n P_j(z_j)\Big\}.
\end{split}
\end{equation}
We call the integer  $m_{j}= 2 + \text{degree}(\triangle P_{j})$ the ``degree'' of $P_{j}$. (The actual degree of $P_{j}$ may be larger, but the addition of a harmonic polynomial to $P_{j}$ does not affect our analysis, and can be eliminated by a change of variables.)
We identify $M$ with $\C^n \times \R$ so that the point $\big( z_{1}, \ldots, z_n, t+i\big(\sum_j P_j(z_j)\big)\big) \in M$ corresponds to the point $(z_1, \ldots, z_n, t)\in \C^n \times \R$. $M$ has real dimension $2n+1$. When integrating on $M$, we take the measure to be Lebesgue measure on $\C^n\times\R$. 

In addition to the boundary of a decoupled domain as in (\ref{E:1.0.2}), we also consider Cartesian products of boundaries of domains in $\C^2$. For $1 \leq j \leq n$, let
\begin{equation}\label{E:1.2.1w}
\begin{split}
\Omega_j &= \left\{(z_j,w_j) \in \C^2\,\Big\vert\,\IM[w_j] > P_j(z_j)\right\};\\ 
M_j &= \left\{(z_j,w_j) \in \C^2\,\Big\vert\,\IM[w_j] = P_j(z_j)\right\}.
\end{split}
\end{equation}
As before, we identify $M_j$ with $\C\times \R$ so that the point $\big(z_j, t+iP_j(z_j)\big)$ corresponds to the point $(z_j,t)$. When integrating on $M_j$ we use Lebesgue measure on $\C\times \R$. The Cartesian product of these boundaries is
\begin{equation}\label{Equation2.1.3.19}
\widetilde M = M_1 \times \cdots \times M_n \subset \C^{2n}.
\end{equation}
Then $\widetilde M$ is the Shilov boundary of the product domain $\Omega_1 \times \cdots \times \Omega_n$. It has real dimension $3n$ and real codimension $n$. We can identify $\widetilde M$ with $\C^n \times \R^n$ so that the point $p= \big(z_1, t_1 + i P_1(z_1), \ldots, z_n, t_n + i P_n(z_n)\big) \in \widetilde M$ corresponds to the point $(z_1,\ldots,z_n,t_1,\ldots, t_n) = (z,t)
\in \C^n \times \R^n$. When integrating on $\widetilde M$, we take the measure to be Lebesgue measure on $\C^n\times\R^n$.

Let $\pi:\C^{2n} \to \C^{n+1}$ be the linear holomorphic mapping
\begin{equation*}\label{E:1.1.7o2}
\pi(z_1,\ldots,z_n,w_1,\ldots,w_n) = (z_1, \ldots, z_n,w_1 + \cdots + w_n).
\end{equation*}
This induces a mapping from $\widetilde M$ to $M$. In terms of the coordinates given by $\C^n \times \R^n$ and $\C^n \times \R$, we have
\begin{equation}\label{E:1.2.4d}
\pi(z_1,\ldots,z_n,t_1,\ldots,t_n) = (z_1,\ldots,z_n, t_1 + \cdots + t_n).
\end{equation}
The mapping $\pi$ allows us to transfer functions from $\widetilde M$ to $M$. If $\varphi \in\break
\mathcal C^\infty_0(\C^n\times \R^n)$, we define a function $\varphi^{\#} \in \mathcal
C^\infty_0(\C^n\times \R)$ by setting
\begin{equation}\label{E:1.1.7.1}
\begin{split}
\varphi^{\#}(z,t) &= 
\int_{\R^{n-1}} \varphi\Big(z, r_1, \ldots, r_{n-1}, t-\sum_{j=1}^{n-1}r_j\Big)\,dr_1\cdots dr_{n-1}\\
&\equiv \int_{r \in \Sigma(t)} \varphi(z,r)\,d\tilde r
\end{split}
\end{equation}
where $\Sigma(t) = \{(r_1, \ldots, r_n) \in \R^n\,\big\vert\,r_1 + \cdots + r_n = t\}$ and $d\tilde r$ is $(n-1)$-dimensional
\pagebreak
 Lebesgue measure on $\Sigma(t)$.

\Subsec{The $\dbar_b$-complex and the $\bx_b$ operator on $M$ and $\widetilde M$} Let $M$ be the
boundary of a decoupled domain as in (\ref{E:1.0.2}). Using coordinates $(z_1, \ldots, z_n,t)\break \in \C^n
\times
\R$, bases for the Cauchy-Riemann operators of type $(1,0)$ and $(0,1)$ are given by the operators
$\{Z_j,\,1\leq j \leq n\}$ and by $\{\bar Z_j,\,1\leq j \leq n\}$ where
\begin{equation}\label{E:CauchyRiemann}
\begin{split}
Z_j &= \frac{\partial}{\partial z_j} +i\,\frac{\partial P_j}{\partial z_j}(z_j)\,\frac{\partial}{\partial t} = X_j - i X_{n+j},
\\
\bar Z_j &= \frac{\partial}{\partial \bar z_j} -i\,\frac{\partial P_j}{\partial \bar z_j}(z_j)\,\frac{\partial}{\partial t} = X_j + i X_{n+j},
\end{split}
\end{equation}
where $\{X_1,\ldots, X_{2n}\}$ are real vector fields.

\demo{Remark} For future reference, note that the operators $Z_j$, $\bar Z_j$, and their sums and products
commute with translations in the variable $t$. The same will also be true of the inverses or relative inverses
we construct for such operators. Hence the corresponding distribution kernels $K\big((z,t),(w,s)\big)$ will be
of the form $K(z,w,t-s)$.
\Enddemo

We recall the formalism of the $\dbar_b$-complex on $M$. If $f$ is a function, then
\begin{equation*}
\dbar_b[f] = \sum_{j=1}^n \bar Z_j[f]\,d\bar z_j.
\end{equation*}
Let $\I_q$ denote the set of strictly increasing $q$-tuples of integers between $1$ and $n$. Let $J= \{j_1,\ldots,j_q\} \in \I_q$, and let $d\bar z_J$ denote the $(0,q)$-form $d\bar z_{j_1}\wedge\cdots \wedge d\bar z_{j_q}$. Then $\{d\bar z_J\}_{J\in \I_q}$ is a basis for the space of $(0,q)$ forms and
\begin{equation*}
 \dbar_b\Big[\sum_{J\in \I_q} f\,d\bar z_J \Big] =
 \sum_{J\in \I_q} \dbar_b[f]\,\wedge\,d\bar z_J.
\end{equation*}
One checks that $\dbar_b^2 = 0$. 

Let $\dbar_b^*$ denote the formal adjoint of $\dbar_b$ so that $\dbar_b^*$ maps $(0,q+1)$-forms to $(0,q)$-forms. Thus for compactly supported $(0,q)$ and $(0,q+1)$ forms $\varphi$
 and~$\psi$,
\begin{equation*}
\big\langle \dbar_b[\varphi],\psi\big\rangle_{q+1} =
\big\langle\varphi,\dbar_b^*[\psi]\big\rangle_q,
\end{equation*} 
where $\big\langle \,\cdot\,,\,\cdot\,\big\rangle_q$ is the $L^2$-inner product on $(0,q)$-forms defined so that the forms $\{d\bar z_J\}_{J\in \I_q}$ are orthonormal.

The Kohn-Laplacian
\begin{equation}\label{EKohn-Laplacian}
\bx_b =\dbar_b\,\dbar_b^* +\dbar_b^*\,\dbar_b
\end{equation}
is a second order system of partial differential operators which maps $(0,q)$-forms to $(0,q)$-forms. For the decoupled boundary $M$, $\bx_b$ acts as follows. For $1 \leq j \leq n$, let
\begin{equation}\label{E:1.3.2z}
\begin{split}
\bx_j^{(+)} &= - \bar Z_j\,Z_j;\\\bx_j^{(-)} &= - Z_j\,\bar Z_j.
\end{split}
\end{equation}
For $J \in \I_q$ and $1 \leq k \leq n$ set
\begin{equation*}
J(k) =
\begin{cases}
(+) & \text{if $k \in J$},\\ (-) & \text{if $k \notin J$}.
\end{cases}
\end{equation*}
The operator $\bx_b$ acts diagonally, and is given by
\begin{equation}\label{E:1.3.4z}
\bx_b\left(\sum_{J\in \I_q} \varphi_J\,d\bar z_J\right) = \sum_{J\in \I_q}
\bx_J(\varphi_J)\,d\bar z_J
\end{equation}
where
\begin{equation}\label{E:1.3.5z}
\bx_J = \sum_{k=1}^n \bx_k^{J(k)}.
\end{equation}
Thus, the study of the $\dbar_b$ complex on $M$ on $(0,q)$-forms is reduced to the study of the $\binom{n}{q}$ operators $\bx_J$ for $J \in \I_q$.

We can also consider the $\dbar_b$-complex on the product submanifold $\widetilde M$. Instead of the vector fields (\ref{E:CauchyRiemann}), we set
\begin{equation*}\label{E:CauchyRiemannB}
\begin{split}
Z_j &= \frac{\partial}{\partial z_j} +i\,\frac{\partial
P_j}{\partial z_j}(z_j)\,\frac{\partial}{\partial t_j};\\ \bar Z_j
&= \frac{\partial}{\partial \bar z_j} -i\,\frac{\partial
P_j}{\partial \bar z_j}(z_j)\,\frac{\partial}{\partial t_j}.
\end{split}
\end{equation*}
The $\dbar_b$ complex on $\widetilde M$ is defined in the exactly the same way as on $M$. If we then define operators $\bx_j^\pm$ as before, the operator $\bx_b = \dbar_b\,\dbar_b^* + \dbar_b^*\,\dbar_b$ has exactly the same form as in equations (\ref{E:1.3.4z}) and (\ref{E:1.3.5z}).

The mapping $\pi:\widetilde M \to M$ given in equation (\ref{E:1.2.4d}) induces a mapping from functions on $M$ to functions on $\widetilde M$, and hence induces a mapping $d\pi$ from tangent vectors on $\widetilde M$ to tangent vectors on $M$. In particular, if $T_j = \frac{\partial} {\partial t_j}$ on $\widetilde M$ and $T = \frac{\partial}{\partial t}$ on $M$, then $d\pi(T_j) = T$ for $1 \leq j \leq n$. This justifies our use of the same notation, 
i.e.\ $Z_j$ and $\bar Z_j$ for vectors fields and $\bx_b$ for the Kohn Laplacian, on both $\widetilde M$
and $M$. The adjoint mapping $d\pi^*$ which carries differential forms on $M$ to differential forms on
$\widetilde M$ commutes with the mappings $\dbar_b$.

We have also considered a mapping $\varphi \to \varphi^{\#}$ in (\ref{E:1.1.7.1}) which carries functions on $\widetilde M$ to functions on $M$. The following is then clear.

\begin{proposition} Let $\varphi \in \mathcal C^\infty_0(\widetilde M)$ so that $\varphi^{\#} \in \mathcal C^\infty_0(M)$.
 Then $T[\varphi^{\#}]\break = (T_j[\varphi])^{\#}$ and so in particular $(T_j[\varphi])^{\#} =
(T_k[\varphi])^{\#}$ for any $1 \leq j,k \leq n$.
\end{proposition}

\vglue-12pt
\Subsec{Outline of the argument}\label{SS:Outline}  We now expand   the discussion in Section
\ref{MethodsUsed} and describe the main ideas involved in the construction of relative fundamental
solutions for the operators $\{\bx_J\}$. Let $W_j$ denote either $Z_j$ or $\bar Z_j$ where $W_j$ is then a
first order differential operator on $M$ which depends only on the variables $z_j$ and $t$, and which
commutes with translation in $t$. Let $\bx_j = W_j^*\,W_j = - \overline W_j\,W_j$ be the corresponding
nonnegative, self-adjoint second order differential operator on $L^2(M)$. Let $\S_j$ be the orthogonal
projection of
$L^2(M)$ onto the null space of $\bx_j$. Note that this space is the same as the null space of $W_j$. Let
$\{e^{-s\bx_j}\}$ be the semi-group of contractions on $L^2(M)$ with infinitesimal generator $\bx_j$. Then $\bx = \sum_{j=1}^n\bx_j$ is one of the operators $\bx_J$. 

Let $M_j$ be the boundary in $\C^2$ given in equation (\ref{E:1.2.1w}). Then $\bx_j$ and $e^{-s\bx_j}$ also act on $L^2(M_j)$. From the theory of domains of finite type in $\C^2$, it is known that the projection $\S_j$ of $L^2(M_j)$ onto the null space of $\bx_j$ and the heat kernel $\H_j$ for the semi-group $\{e^{-s\bx_j}\}$ are given by operators
\begin{equation*}\label{E:1.2.1y}
\begin{split}
\S_j[f](z_j,t) &= \int_{\C\times\R} f(w_j,r) S_j(z_j,w_j,t-r) \,dw_j\,dr,\\ 
\H_{j}[f](s,z_j,t) &= \int_{\C\times\R} f(w_j,r) H_{j}(s, z_j,w_j,t-r) \,dw_j\,dr,
\end{split}
\end{equation*}
where $S(z,w,t)$ and $H_j(s,z,w,t)$ are distributions on $\C\times\C\times\R$ and on $(0,\infty)\times \C\times \C \times \R$.
We can think of the projection onto the null space of $\bx_j$ and the heat kernel $e^{-s\bx_j}$ as operators either on $M_j$ or on
the decoupled boundary $M$. In other words, when thinking of these operators acting on $M$ we can also write
\begin{eqnarray*}\label{Eq001}
\S_j[f](z_1,\ldots,z_n,t) &=& \int_{\C\times\R} f(z_1, \ldots,w_j,\ldots,z_n,r)\, S_j(z_j,w_j,t-r)\,dw_j\,dr;\\
\H_j[f] (z_1,\ldots,z_n,t) &=& \int_{\C\times\R} f(z_1, \ldots,w_j,\ldots,z_n,r)\, H_j(s,z_j,w_j,t-r)\,dw_j\,dr.
\end{eqnarray*}

We remark that the distribution kernels $S_j$ and $H_{j}$ are the limits of certain distributions $S_j^\varepsilon$ and $H_{j}^\varepsilon$ which
are given by integration against infinitely differentiable functions with compact support. All the estimates we shall make on $S_j$ and $H_j$ hold for the approximating kernels, and the estimates are uniform in $\varepsilon$. 

Now define an operator $\S$ on $L^2(M)$ by setting 
\begin{equation*}\label{Eq002}
\S = \prod_{j=1}^n \S_j.
\end{equation*}
This is the orthogonal projection onto the intersection of the null spaces of the operators $\{\bx_1, \ldots, \bx_n\}$, which is the same as the projection onto the null space of the operator $\bx = \bx_1 + \cdots +\bx_n$. Also note that since the operators $\{\bx_j\}$ commute, the operator $\H = \exp\big[-s\sum_{j=1}^n \bx_j\big]$ is just a product:
\begin{equation*}\label{Eq003}
\H = e^{-s\,\big[\sum_{j=1}^n \bx_j\big]} = \prod_{j=1}^n e^{-s\,\bx_j} = \prod_{j=1}^n \H_j.
\end{equation*}
It thus follows that the operators $\S$ and $\H$ are then given on $M$ by integral operators
\begin{equation*}
\begin{split}
\S[f](z,t) &= \int_{\C^n\times \R} f(w,r) \,S(z,w,t-r)\,dw\,dr\\
 e^{-s\bx}[f](z,t) &= \int_{\C^n\times \R} f(w,r) \,H(s,z,w,t-r)\,dw\,dr
\end{split}
\end{equation*}
where the kernels $S(z,w,t)$ and $H(s,z,w,t)$ are given by the convolutions
\begin{equation*}
\begin{split}
S(z,w,t) &= \int_{\Sigma(t)} \prod_{j=1}^n S_j(z_j,w_j,r_j)\, d\tilde r;
\\ H(s,z,w,t) &= \int_{\Sigma(t)} \prod_{j=1}^n H_j(s,z_j,w_j,r_j)\, d\tilde r.
\end{split}
\end{equation*}
Here $\Sigma(t)$ and $d\tilde r$ are defined in equation (\ref{E:1.1.7.1}).

Now by the spectral theorem, the operator $e^{-s\big[\sum_{j=1}^n \bx_j\big]}$ converges strongly to $\S$ as $s\to \infty$, and
hence it will be easy to check that
\begin{equation*}
\K = \int_0^\infty \Big[e^{-s\,\big[\sum_{j=1}^n
\bx_j\big]}-\S\Big]\,ds
\end{equation*}
is a relative fundamental solution for $\bx = \sum_j\bx_j$ in the sense that
\begin{equation*}
\K\,\bx = \bx\,\K = I - \S.
\end{equation*}
Also, if we set
\begin{equation*}
\N = \int_0^\infty \prod_{j=1}^n \Big(e^{-s\bx_j} -\S_j\Big)\,ds,
\end{equation*}
then $\N$ is also a relative fundamental solution for $\bx$ in the sense that
\begin{equation*}
\N\,\bx = \bx\,\N = \prod_{j=1}^n (I-\S_j).
\end{equation*}

Since the distribution kernels $\{S_j(z_j,w_j,t)\}$ and $\{H_j(s,z_j,w_j,t)\}$ are known, we have the following explicit formulas for the distribution kernels $K$ and $N$ for relative fundamental solutions for $\bx$ on $M$:
\begin{align*}
K(z,w,t) &= \int_0^\infty \int_{\Sigma(t)}\Big[\prod_{j=1}^n H_j(s,z_j,w_j,r_j) - \prod_{j=1}^n S_j(z_j,w_j,r_j)\Big]\,d\tilde r\,ds,
\\ N(z,w,t) &= \int_0^\infty\int_{\Sigma(t)} \prod_{j=1}^n \Big(H_j(s,z_j,w_j,r_j) - S_j(z_j,w_j,r_j)\Big)\,d\tilde r\,ds.
\end{align*}

There is one further idea in analyzing these last integrals. We write
\begin{align*}
\widetilde K(z,w,r) &= \int_0^\infty \Big[\prod_{j=1}^n H_j(s,z_j,w_j,r_j) - \prod_{j=1}^n S_j(z_j,w_j,r_j)\Big]\,ds,\\
\widetilde N(z,w,r) &= \int_0^\infty \prod_{j=1}^n \Big(H_j(s,z_j,w_j,r_j) - S_j(z_j,w_j,r_j)\Big)\,ds.
\end{align*}
These are the kernels of operators which are relative fundamental solutions for the operator $\sum_{j=1}^n \bx_j$ acting not on $M$, but on the Cartesian product $\widetilde M = M_1\times \cdots \times M_n$. Then, at least formally, we have
\begin{align*}
K(z,w,t) = \int_{\Sigma(t)}\widetilde K(z,w,r)\,d\tilde r = (\widetilde K)^{\#}(z,w,t)\\
N(z,w,t) = \int_{\Sigma(t)}\widetilde N(z,w,r)\,d\tilde r = (\widetilde N)^{\#}(z,w,t).
\end{align*}

We shall first analyze these relative fundamental solutions $\widetilde \K$ and $\widetilde \N$ on the product, and then
integrate over $\Sigma(t)$ to obtain the relative fundamental solution on $M$. By doing this, we can take
advantage of the product structure of the operators in the integrand, and in fact use product theory to establish
$L^p$ regularity. Transference methods show that integration over $\Sigma(t)$ then gives the same $L^p$
regularity for $\K$.

\Subsec{Statement of results} We summarize our main results about relative inverses to the operators
$\bx_J$ on the decoupled manifold $M$. The Kohn-Laplacian operator $\bx_b$ has an infinite dimensional
null space in $L^2(M)$ when acting on $(0,0)$-forms or on $(0,n)$-forms, but has no null space in
$L^2(M)_{(0,r)}$ when acting on $(0,r)$-forms for $1 \leq r \leq n-1$. Let $\S_0$ denote the orthogonal
projection of $L^2(M)$ onto the null space of $\bx_b$ acting on functions, and let $\S_n$ denote the
orthogonal projection of $L^2(M)_{(0,n)}$ onto the null space of $\bx_b$ acting on $(0,n)$-forms. We
show that each $\Box_J$ has an inverse modulo the relevant projection.

\begin{theorem} \label{TH1}
For each of the $2^n$ possible operators $\{\bx_J\}${\rm ,}
 we construct a distribution $K_J$ on $M\times M$ so
that if $\mathcal{K}_J$ denotes the linear operator\/\footnote{By abuse of notation, this means that
$\K_J[\varphi]$ is the distribution on $M$ such that $\langle \K_J[\varphi],\psi\rangle = \langle
K_j,\varphi\otimes \psi\rangle$ for $\varphi, \psi \in \CC^\infty_0(M)$.}
\begin{equation*}\label{E:1.0.3}
\K_J[\varphi](p)= \int_M \varphi(q)\,K_J(p,q)\,dq,
\end{equation*}
then
$$
\K_J\,\bx_J = \bx_J\,\K_J =
\begin{cases}
I-S_0 &\text{if \, $\bx_J$ acts on functions};
\\ I & \text{if \, $\bx_J$ acts on a $(0,r)$-form with $1\leq r \leq n-1$};
\\ I-S_{n} &\text{if \, $\bx_J$ acts on $(0,n)$-forms}.
\end{cases}
$$
\end{theorem}

We study the regularity of the operators $\K_J$ in both $L^p$ and H\"older spaces on $M$. Since $\K_J$ inverts a second order operator (modulo bounded projections), one expects that $\K_J$ should behave like an operator smoothing of order two. However as
we have already pointed out in Section \ref{Background}, Derridj \cite{Derr78} showed that maximal hypoelliptic estimates are possible only if the eigenvalues of the Levi form degenerate at the same rate. In particular, for decoupled domains in $\C^{n+1}$ with $n>1$, the operator $\bx_b$ on $M$ fails to be maximally subelliptic near $p \in M$ whenever $\triangle P_j(p) = 0$ for some $j$. Thus $Q(Z, \bar Z)\,\K_J$ cannot be a bounded operator on $L^2(M)$ for an arbitrary quadratic combination of ``good'' derivatives $\{Z_1, \ldots, Z_n,\bar Z_1, \ldots, \bar Z_n\}$.

\begin{theorem}\label{TH2}
Write $\Box_J = W_1\,\bar W_1+ \cdots + W_n\,\bar W_n$ 
where each $W_j$ is one of the two operators $\{Z_j, \bar Z_j\}${\rm ,}
 and $\bar W_j$ is the other. Let $\K_{J}$
be the distribution constructed in Theorem {\rm \ref{TH1}}. Then\/{\rm :}\/
\begin{enumerate}

\item When $1 \leq k,l\leq n$ with $l\neq k${\rm ,} 
  the operators $W_k\,\bar W_k\,\K_J = -\bx_k\, \K_J${\rm ,} 
$\bar W_k\,\bar W_k\,\K_J${\rm ,} and $\bar W_l\,\bar W_k\,\K_J$ 
extend to bounded linear operators on $L^p(M)$ for 
$1 < p < \infty$.

\item If $B_k$ is a bounded smooth function on $M$ and if there are constants $C_{k,l}$
 so that for all $p = (z_1, \ldots, z_n,t) \in M$ 
\begin{equation*}
  |B_k(p)|\,\triangle P_k(z_k) \leq C_{k,l} \,\triangle P_l(z_l)
\end{equation*}
for $1 \leq k \leq n$ and all $l${\rm ,} then the operator 
$B_k \,\bar W_k\,W_k\,\K_J = - B\,\overline \bx_k\, \K_J$ 
extends to a bounded linear operator on $L^p(M)$ for $1 < p < \infty$.

\item If $B_k$ is a bounded smooth function on $M$ and if there are constants $C_{k}$
 so that for all $p = (z_1, \ldots, z_n,t) \in M$
 \begin{equation*}
  |B_k(p)| \leq C_{k} \inf_{l\neq k} \triangle P_l(z_l)
 \end{equation*}
then the operator $B_k \,W_k\,W_k\,\K_J$ 
extends to a bounded linear operator on $L^p(M)$ for $1 < p < \infty$.
\end{enumerate}
\end{theorem}

 Note that this theorem does not make any assertion about operators of the form $W_{l}\,\bar
W_{k}\K_{J}$. Theorems \ref{TH1} and \ref{TH2} will be proved in Section 
\ref{SSDistributions} and Section \ref{SSLPregularity} below.  That
the estimates are optimal is shown in Section 9.

Precise Lipschitz regularity results require the introduction of various metrics on the space $M$. 
This is done below in Section \ref{S:Further estimates}. At this stage, however, we can state the following
global result involving the standard isotropic Lipschitz spaces.

\begin{theorem}\label{TH3} Let $m = \max\{m_j\}$ be the largest of the 
degrees of the polynomials $P_j$. \/{\rm (}\/This is the \/{\rm ``}\/type\/{\rm ''}\/ of the boundary $M$.\/{\rm )}\/ 
Assume $m > 2${\rm ,} and 
suppose that $f$ is a function bounded and supported on a ball of radius one in~$M$. 
Then for all $J$ there is a constant $C$ so
that if $h \in \C^n\times \R \cong M${\rm ,} then
\begin{equation*}
\left\vert \K_J[f](p+h) - \K_J[f](p) \right\vert \leq C\,\vert h \vert^{\frac{2}{m}}.
\end{equation*}
\end{theorem}

There is also a corresponding result when $m = 2$ (the strongly
pseudo-convex case).  These results are presented in Section 8. 
 
We shall show that the distributions $K_J$ are singular only on the diagonal of $M \times M$, and we obtain estimates on the size of these distributions and their derivatives away from these singularities. These estimates involve two different pseudo-metrics on $M$. The first, which we call $d_\Sigma$, is the control metric associated to the collection of vector fields which are the real and imaginary parts of the vector fields $\{\bar Z_1, \ldots, \bar Z_n\}$. The second pseudo-metric, which we call $d_S$, describes the singularities of the singular integral which gives the Szeg\"o projections $S_0$ and $S_n$. The corresponding balls are denoted by $B_\Sigma(p,\delta)$ and $B_S(p,\delta)$. These metrics are described in Section \ref{S:Geometries}. The estimates of $K_J$ also involve functions $\mu_j(p,\delta)$ which are defined in Section \ref{SS:2.2a}, equation (\ref{E:2.2.6.3}) below.

\begin{theorem} \label{TH4} 
The distribution $K_J(p,q)$ is given by integration against a $\CC^\infty$ function away from the diagonal
$\{p=q \in M\}${\rm ,}
 where there are  the following estimates. Let $\partial_j^{\alpha_j}$ be a derivative of order
$|\alpha_j|$ made up of the vector fields $Z_j$ and $\bar Z_j$ in which each acts in either the variables 
$p_j$
or $q_j$. Then for all $\alpha = (\alpha_1,\ldots, \alpha_n)$ there is a constant $C_\alpha= C$ so that
\begin{eqnarray*}
&&
\Big\vert \Big[\prod_{j=1}^n \partial_j^{\alpha_j}\Big] K_J(p,q)\Big\vert \\
&&\qquad \leq C
\frac{\left[\sum_{j=1}^n \mu_j(p,d_S(p,q))\right]}{\Big\vert B_S(p,d_S(p,q))\Big\vert}^2\,\log\left[2 +
\frac{\left[\sum_{j=1}^n \mu_j(p,d_S(p,q))\right]}{d_\Sigma(p,q)}\right]
\\&&
\qquad\qquad\prod_{j=1}^n \left[\mu_j(p,d_S(p,q))^{-1} + d_\Sigma(p,q)^{-1}\right]^{|\alpha_j|} .
\end{eqnarray*}
\end{theorem}

\section{Geometry and analysis on $M_j$ and on $M_1\times \cdots \times M_n$}\label{S:hypersurfaces}
\vglue-6pt

In this section we summarize some geometric and analytic results which we require later in the paper. For $1 \leq j \leq n$, let $M_j$ be the hypersurface given in equation (\ref{E:1.2.1w}). Let $\widetilde M$ be the Cartesian product as in (\ref{Equation2.1.3.19}). Subsections \ref{SS:2.2a} through \ref{SS:2.4a} deal with the study of the model hypersurfaces $M_j$ in $\C^2$. Subsections \ref{SS:2.5a} and \ref{SS:NIS} deal with geometry and analysis on $\widetilde M$.

\Subsec{The control metric on $M_j$}
\label{SS:HeatKernel1:2.1}\label{SS:2.2a}  Recall that we write the complex vector field $\bar Z_j = X_j +i
X_{n+j}$ where $\{X_j, X_{n+j}\}$ are real vector fields on $M_j$. Define a metric $d_j$ on $M_j$ as
follows. If $p,q \in M_j$ and $\delta >0$, let $AC(p,q,\delta)$ denote the set of absolutely continuous
mappings $\gamma:[0,1]\to M_j$ such that $\gamma(0) = p$ and $\gamma(1) = q$, and such that for almost
all $t \in [0,1]$ we have $\gamma'(t) = \alpha_j(t)\,X_j(\gamma(t)) + \alpha_{n+j}(t)\,X_{n+j}(\gamma(t))$
with $|\alpha_j(t)|^2 + |\alpha_{n+j}(t)|^2 < \delta^2$. Then we define
\begin{equation*}\label{E:2.2.3.0.1}
d_j(p,q) = \inf\left\{\delta>0\,\Big\vert \, 
{AC}(p,q,\delta) \neq \emptyset\right\}.
\end{equation*}
The corresponding nonisotropic ball is
\begin{equation*}
B_j(p,\delta) =\left\{ q \in M_j\,\Big\vert\,d_j(p,q) <
\delta\right\},
\end{equation*}
and $\big\vert B_j(p,\delta)\big\vert$ denotes its volume. Set
\begin{equation*}\label{E:2.2.3.0.11}
V_j(p,q) = \big\vert B_j\big(p,d_j(p,q)\big)\big\vert.
\end{equation*}

The volume of the ball $B(p,\delta)$ is essentially a polynomial in $\delta$ with coefficients that depend on $p$. Let $T = \frac{\partial}{\partial t}$ so that at each point of $M_j$ the tangent space is spanned by the vectors $\{X_j, X_{n+j}, T\}$. Write the commutator
\begin{equation}\label{E:2.2.3.0.1.1}
[X_j,X_{n+j}] = \lambda_{j}\,T + a_j\,X_j + a_{n+j}\,X_{n+j}
\end{equation}
where $\lambda_j, a_j, a_{n+j} \in \CalC^\infty(M_j)$. If $\alpha = (\alpha_{1},\ldots,\alpha_{k})$ is a $k$-tuple with each $\alpha_{j}$ equal to $j$ or $n+j$, let $|\alpha| = k$ and let $X^{\alpha} = X_{\alpha_{1}}\,\cdots\,X_{\alpha_{k}}$ denote the corresponding $k^{\text{th}}$ order differential operator. For $k \geq 2$ set
\begin{equation*}
\Lambda_j^k(p) = \sum_{|\alpha|\leq k-2}|X^{\alpha}\lambda_j(p)|,
\end{equation*}
where $\lambda_j$ is as  defined in (\ref{E:2.2.3.0.1.1}), and set
\begin{equation*}
\Lambda_j(p,\delta)= \sum_{k=2}^{m_j}\Lambda_j^k(p)\,|\delta|^k.
\end{equation*}
We now have

\begin{proposition}\label{P:2.2.1.1} There are constants $C_1,\,C_2$ depending only on $m_j$ so that for $p\in M_j$ and $\delta > 0$
$$
C_1\,\delta^2\,\Lambda_j(p,\delta) \leq \big |B_j(p,\delta)\big| \leq C_2\,\delta^2\,\Lambda_j(p,\delta). 
$$
Also{\rm ,} $V_j(p,q)\approx V_j(q,p) \approx d_j(p,q)^2 \, \Lambda_j\left(p,d_j(p,q)\right)$ where $A \approx B$ means that the ratio $A/B$ is bounded and bounded away from zero.
\end{proposition}

There is an alternate description of the balls $\{B_j(p,\delta)\}$ and metric $d_j$ given in terms of explicit inequalities. For $z, w \in \C$ let
$$
T_j(w,z) = 2\,\IM\Big[\sum_{k=1}^{m_j} \frac{\partial^k P_j}{\partial z^k}(w)\,\frac{(z-w)^k}{k!}\Big].
$$
Then with $p = (w,s)\in M_j$, set
$$
\widetilde B_j (p,\delta ) = \left\{(z,t) \in M_j\,\Big\vert\,
\text{$|z-w|<\delta$ and $|t-s+T_j(w,z)|<\Lambda_j(w,\delta)$}\right\}.
$$

Note that for $\delta > 0$, $\delta \to \Lambda_j(p,\delta)$ is a monotone increasing function. Hence there is a unique inverse function $\mu_j(p,\delta)$ such that for $\delta \geq 0$ we have $\Lambda_j\big(p,\mu_j(p,\delta)\big) = \mu_j\big(p,\Lambda_j(p,\delta)\big) = \delta$. We have
\begin{equation}\label{E:2.2.6.3}
\mu_j(p,\delta)^{-1} \approx \sum_{k=2}^{m_j}
\Lambda_j^k(p)^{\frac{1}{k}}\,|\delta|^{-\frac{1}{k}}.
\end{equation}

\begin{proposition}\label{P:2.2.1.2} There are constants $C_1$ and $C_2$ depending only on $m_j$ so that for $p \in M_j$ and $\delta > 0$
\begin{equation*}
\widetilde B_j(p,C_1\,\delta) \subset B_j(p,\delta) \subset
\widetilde B_j(p,C_2,\delta).
\end{equation*}
Moreover{\rm ,} if $(z,t), (w,s) \in M_j${\rm ,}
$$
d_j\big((z,t),(w,s)\big) \approx |z-w| + \mu_j\big(w,\big\vert
t-s-T_j(w,z)\big\vert\big).
$$
\end{proposition}

\vskip8pt

  Proofs of Propositions \ref{P:2.2.1.1} and \ref{P:2.2.1.2} can be found in \cite{NaStWa85}.

\vskip3pt
We say that the manifold $M_j$ is normalized at a point $(z,t)$ if $\frac{\partial^k P_j}{\partial z^k}(z) = 0$ for $1 \leq k \leq m_j$. If $p = (a,s+iP_j(a)) \in M_j \subset \C^2$, one can always make a biholomorphic change of variables of $\C^2$ which moves $p$ to the origin, and which carries the manifold $M_j$ to a new manifold $M_j^p$ of the same type which is normalized at the origin. (See Subsection \ref{SS:Sum of squares}
 below for further details.) It follows from Proposition \ref{P:2.2.1.2} that if the domain $M_j$ is normalized
at the origin, then the balls and distances have the particularly simple form:
\begin{equation*}
\begin{split}
\widetilde B_j((0,0),\delta) &= \Big\{(z,t) \in \C\times \R\,\Big\vert\,\text{ $|z| < \delta$ and $|t| < \Lambda_j(0,\delta )$}\Big\};\\ d_j\big((z,t),(0,0)\big) &\approx |z| + \mu_j(0,t).
\end{split}
\end{equation*}

\Subsec{{\rm NIS} operators on $M_j$}\label{SS:HeatKernel1:2.2} We briefly review the definition of the class
of nonisotropic smoothing (NIS) operators on $M_j$. These were introduced in \cite{NaRoStWa89}, and the
definition below is taken from Koenig \cite{Koenig01}. Let 
\begin{equation*}
\T[f](p) = \int_M T(p,q)\,f(q)\,dq
\end{equation*}
where $T$ is a distribution on $M_j\times M_j$.

\begin{definition}\label{D:NIS} $\T$ is an NIS operator smoothing of order $k$ if the distribution $T$ is given away from the diagonal of $M_j\times M_j$ by integration against a $\CC^\infty$ function and:
\begin{enumerate}
\item There exists $\beta < \infty$, and for $s \geq 0$ there exists $\alpha(s) < \infty$ such that if $\zeta, \zeta'
\in \CC^\infty_0(M_j)$ with $\zeta' \equiv 1$ on the support of $\zeta$, then there exists a constant
$C_{s,\zeta,\zeta'}$ so that for all $f \in \CC^\infty_0(M_j)$
\begin{equation*}
\norm \zeta\,\T[f]\norm_s \leq C_{s,\zeta,\zeta'}\left[\norm \zeta'\,f\norm_{\alpha(s)} 
+ \norm f\norm_\beta\right].
\end{equation*}

\item Let $X^\alpha_p$ and $X^\beta_q$ be derivatives of order $|\alpha|$ and $|\beta|$
 in the vector fields $X_j$ and $X_{n+j}$ acting on the variables $p$ and $q$. There exist  constants
$C_{\alpha,\beta}$ so that for $p \neq q$
\begin{equation*}
\big\vert X^\alpha_p\,X^\beta_q T(p,q)\big\vert \leq
 C_{\alpha,\beta}\, d_j(p,q)^{k-|\alpha|-|\beta|}\,V_j(p,q)^{-1}.
\end{equation*}

\item For each integer $\ell$ there are  an integer $N_{\ell}$ and a constant $C_{\ell}$ so that if $\varphi$ is
a
$\CalC^\infty$ function supported on $B_j(p,\delta)$, then for all $\varepsilon > 0$ and all $\alpha$ with
$|\alpha| = \ell$
\begin{equation*}
\big\vert X^\alpha \T[\varphi](p)\big\vert \leq C_{\ell}\, \delta^{k-\ell} 
\sup_{q\in M} \sum_{|J|\leq N_{\ell}} \delta^{|J|} \big\vert X^J\varphi(q)\big\vert.
\end{equation*}

\item The above conditions also hold for the adjoint operator $\T^*$ with distribution kernel $T(y,x)$.
\end{enumerate}
The constants $C_{s,\zeta,\zeta'}$, $C_{\alpha,\beta}$ and $C_\ell$ are called the NIS constants of the
operator~$T$. \end{definition}
 
\begin{definition}  With regard to condition (3), if $\varphi$ is a smooth function
supported  on $B_j(p,\delta)$, we say that $\varphi$ is {\it a normalized
bump function}  if  
\begin{equation*}
\sup_{q\in M}\,\sum_{|J|\leq N_{\ell}} \delta^{|J|} \big\vert X^J\varphi(q)\big\vert \leq 1.
\end{equation*}
\end{definition}

\vglue-12pt
\Subsec{The $\dbar_b$-complex on $M_j$} \label{SS:2.4a} We shall need the following basic results
concerning the $\dbar_b$ complex on $M_j$. Let $Z_j$ and $\bar Z_j$ be the tangential Cauchy-Riemann
operators of type $(1,0)$ and $(0,1)$ on $M_j$ as given in equation (\ref{E:CauchyRiemann}). Then
$\dbar_b[f] = \bar Z_j[f]\,d\bar z$, and the formal adjoint is $\dbar_b^*[g\,d\bar z] = - Z_j[g]$. We have
operators $\bx_j^{(\pm)}$ as in equation (\ref{E:1.3.2z}). The Kohn-Laplacian is then $\bx_j^{(-)}$ when
acting on functions, and is $\bx_j^{(+)}$ when acting on $(0,1)$-forms. Each operator $\bx_j^\pm$ extends
to a closed, densely defined, nonnegative self-adjoint operator on $L^2(M_j)$. Let $\S_j^{(\pm)}$ denote
the orthogonal projection of $L^2(M_j)$ onto the null space of $\bx_j^{(\pm)}$. The operator
$\S_j^{(\pm)}$ is induced by a distribution $S_j^{(\pm)}(p,q)$ on $M_j\times M_j$ which is given away
from the diagonal by integration against a $\CC^\infty$ function. Let $e^{-s\bx_j^{(\pm)}}$ denote the
semi-group generated by $\bx_j^{(\pm)}$. For each $s > 0$, there is a distribution heat kernel
$H_j^{(\pm)}(s,p,q)$ on $ M_j \times M_j$ such that
\begin{equation*}\label{E:2.2.3.1}
e^{-s\bx_j^{(\pm)}}[f](p) = \int_{M_j} H_j^{(\pm)}(s,p,q)\,f(q)\,dq.
\end{equation*}

The analysis of $\bx_j^{(-)}$ and $\bx_j^{(+)}$ is similar, and from now on we shall omit the superscript. Thus $\bx_j$ will stand for either $\bx_j^{(-)}$ or $\bx_j^{(+)}$, $\S_j$ will denote the projection onto its null space, and $H_j(s,p,q)$ will denote the corresponding heat kernel.

\begin{theorem}\label{T:2.2.3} The operator $\S_j$ is an {\rm NIS}
 operator on $M_j$ smoothing of order zero. Moreover{\rm ,}
 there is an {\rm NIS} operators $\K_j$ smoothing of order
two such that
$$
\K_j\,\bx_j =\bx_j\,\K_j= I - \S_j.
$$
The distribution kernel of the operator $\K_j$ is related to the corresponding
 heat kernel $H_j$ and the projection $S_j$ by the formula
$$
K_j(p,q) = \int_0^\infty H_j(s,p,q) - S_j(p,q)\,ds.
$$
\end{theorem}

  Proofs can be found in \cite{Christ88}, \cite{Christ91a}, \cite{Christ91}, \cite{ChNaSt92},
\cite{NaRoStWa89}, and \cite{NaSt00}.

\begin{theorem} Let $B_j\big(z_j,w_j,\zeta\big)$
 be the Bergman kernel for the domain $\Omega_j$. Then the Szeg\ho\ kernel is given by
$$
S_j(z_j,w_j,t) = \int_0^\infty B_j(z_j,w_j,t+ir)\,dr.
$$
If  
$$
S^\varepsilon_j(z_j,w_j,t) = \int_\varepsilon^\infty
B_j(z_j,w_j,t+ir)\,dr
$$
then $S^\varepsilon_j \to S_j$ as $\varepsilon \to 0${\rm ,}
 and for $\varepsilon > 0$, the kernels $S_j^\varepsilon$ are smooth bounded functions on $M_j \times M_j$.
\end{theorem}

  Further discussion, and the relevant estimates for the Bergman kernels, can be found in \cite{NaRoStWa89}.

\begin{theorem}\label{T:1.1.1} There is a function $G_j\in \CC^{\infty} \big((0,\infty) \times M_{j}\times M_{j}\big)$ so that
$$
H_j(s,p,q) = G_j(s,p,q) + S_j(p,q)
$$
where $S_j(p,q)$ is the distribution kernel for the orthogonal projection operator $\S_j$. In particular{\rm ,}
 the distribution $H_j(s,p,q)$ is given by integration against a $\CC^\infty$ function away from the
diagonal.  There is a constant $C_{\alpha,l}$ with the following property.
 Let $\partial_{p,q}^\alpha$ denote
a differentiation of total order $|\alpha|$ in $Z_j$ and $\overline Z_j${\rm ,} acting either on the $p$ or $q$
variables. Then for $p \neq q$
\begin{equation}\notag
\big\vert\partial^{l}_{s}\,\partial_{p,q}^{\alpha} G_j(s,p,q)\big\vert
\leq
\begin{cases}
C_{\alpha,l}\,d_j(p,q)^{-2l-|\alpha|}\,V_j(p,q)^{-1} & \text{if $s\leq d_j(p,q)^{2}$}
\\\\
C_{\alpha,l}\,s^{-l-\frac{1}{2}|\alpha|}\big\vert B_j(p,\sqrt s)\big\vert^{-1} &
\text{if $s \geq d_j(p,q)^{2}$}.
\end{cases}
\end{equation}
Moreover{\rm ,} for every nonnegative integer $N$ there is a constant $C_{N,\alpha,l}$ so that
\begin{equation}\notag
\big\vert\partial^{l}_{s}\,\partial_{p,q}^\alpha H_j(s,p,q)\big\vert \leq C_{N,\alpha,l}\,
 \frac{d_{M_{j}}(p,q)^{-2l-|\alpha|}}{ V_j(p,q)} \left[\frac{s^{N}}{s^{N} + d_j(p,q)^{2N}}\right].
\end{equation}
\end{theorem}
\vskip8pt

{\it Remark}.  It follows from the estimates on $G_j$ and $|B_j(p,\sqrt s)|$ that for $s \geq
d_j(p_j,q_j)^2$,
\begin{equation}\label{E:2.3.1e}
\Big\vert G_j(s,p_j,q_j)\big\vert \lesssim \left[\big\vert\lambda_j(p_j) \big\vert s^2 +
s^{\frac{1}{2}m_j+1}\right]^{-1}.
\end{equation}
In particular, each function $G_j(s,p_j,q_j)$ is integrable in $s$ at infinity.

\begin{theorem}\label{T:1.2} The operators
\begin{align*}
H_{j,s}[f](p) &= \int_{M_{j}}H_j(s,p,q)\,f(q)\,dq,\\ G_{j,s}[f](p) &=
\int_{M_{j}}G_j(s,p,q)\,f(q)\,dq
\end{align*}
are {\rm NIS} operators smoothing of order zero{\rm ,} and the associated {\rm NIS}
constants are uniformly bounded for $s > 0$.
\end{theorem}

  Proofs of Theorems \ref{T:1.1.1} and \ref{T:1.2} can be
found in \cite{NaSt00}.

\Subsec{Geometry on Cartesian products}\label{SS:2.5a} For $1 \leq j \leq n$, let $M_j$ be a hypersurface
in $\C^2$ as defined in equation (\ref{E:1.2.1w}), and let $\widetilde M = M_1 \times \cdots \times M_n$.
Each of the nonisotropic distances $d_j$ on $M_j$ can be regarded as a function on $\widetilde M$ which
depends only on the variables $(z_j,t_j)$.

In addition, there is a nonisotropic metric $d_{\Sigma}$ on $\widetilde M$ induced by all the real vector fields $\{X_1,\ldots, X_{2n}\}$. If $p,q \in M_j$ and $\delta >0$, let $AC(p,q,\delta)$ denote the set of absolutely continuous mappings $\gamma:[0,1]\to \widetilde M$ such that $\gamma(0) = p$ and $\gamma(1) = q$, and such that for almost all $t \in [0,1]$ we have $\gamma'(t) = \sum_{j=1}^{2n} \alpha_j(t)\,X_j(\gamma(t))$ with $\sum_{j=1}^{2n}|\alpha_j(t)|^2 < \delta^2$. Then 
$$d_{\Sigma}(p, q) = \inf \left\{\delta>0 \, \big\vert\,\text{AC}(p,q,\delta) \neq \emptyset \right\}.$$

This metric is appropriate for describing the fundamental solution of the operator ${\mathcal L}= \sum_{j=1}^{2n} X_j^2$, and we refer to it as the {\it sum of squares} metric. (See \cite{NaStWa85} for a discussion of the relationship between the metric and the operator ${\mathcal L}$.) The metric $d_{\Sigma}$ can be explicitly described as follows. Let $p = (z_1,t_1, \ldots, z_n,t_n) \in \widetilde M$. We can assume without loss of generality that each manifold $M_j$ is normalized at the origin. We denote the origin of $\widetilde M$
by~$\bar 0$. Then
$$
d_\Sigma(\bar0,p) \approx \sum_{j=1}^n \left[|z_j| + \mu_j(0, |t_j|)\right].
$$
The ball centered at $\bar 0$ of radius $\delta$ is, up to constants, given by
$$
B_\Sigma(\bar 0, \delta) = \Big\{(z,t) \in \widetilde M\,\Big\vert\,\text{$|z_j|<\delta$ and $|t_j| < \Lambda_j(0,\delta)$ for $1\leq j \leq n$}\Big\}.
$$
We have
$$
\big\vert B_\Sigma\left(\bar 0,\delta\right) \big\vert \approx
\delta^{2n}\,\prod_{j=1}^n \Lambda_j(0,\delta),
$$
and
$$
\big\vert B_\Sigma\big(\bar 0,d_\Sigma(z,t)\big)\big\vert \approx \big[\sum_{j=1}^n |z_j| + \mu_j(0,|t_j|)\big]^{2n}\,\prod_{j=1}^n \Lambda_j\Big(0, \big[\sum_{j=1}^n |z_j| + \mu_j(0,|t_j|)\big]\Big).
$$

\Subsec{Product singular integral operators on $M_1\times \cdots \times M_n$} \label{SS:NIS}  We
introduce a class of operators on $\widetilde M$  (product NIS operators)
which is the analogue of the standard product singular integrals in the 
Euclidean setting. (For a discussion of the Euclidean case, see for example \cite{NaRiSt00}.) The definition
will involve differential inequalities  on the distribution kernel, and certain cancellation conditions expressed
in terms of the action of the operator on normalized bump functions.  These are generalizations of the
single-factor NIS operators arising in Section 3.2.
Because of the complicated cancellation conditions, it seems easiest to give the definition of product NIS operators on $\widetilde M$ by induction on the number $n$ of factors. When $n = 1$, we are in the situation discussed in Section \ref{SS:HeatKernel1:2.2}, and a product singular integral operator will just mean a standard NIS operator smoothing of order zero.

In general, product operators on $\widetilde M$ will be induced by distributions which are smooth functions away from the {\it product diagonal} given by
$$
\widetilde D = \left\{\big((p_1,\ldots,p_n); (q_1,\ldots q_n)\big)
\in \widetilde M\times \widetilde M\,\Big\vert\, \text{$p_j = q_j$ for some $1 \leq
j\leq n$}\right\}.
$$

\begin{definition} \label{D:5.5.1} Assume that product singular integral operators have been defined on
products where the number of factors is less than $n$. Let $\widetilde M$ be a product of $n$ hypersurfaces
$M_j$. Then $\widetilde \T$ is  {\it a product singular integral operator on} $\widetilde M$ if $\widetilde
\T$ is a continuous linear mapping from $\mathcal C^\infty_0(\widetilde M)$ to $\mathcal D'(\widetilde M)$
and if:
 \smallbreak 
  (1) \quad The Schwartz kernel $\widetilde T(p,q)$ associated to $\widetilde \T$ is a distribution which is a
$\CC^\infty$ function on the set $(\widetilde M\times \widetilde M\, \backslash\, \widetilde D)$. In
particular, if $\varphi_j, \psi_j \in \CC^\infty_0(M_j)$ have disjoint supports for $1 \leq j \leq n$,
$$
\Big\langle \widetilde \T(\varphi_1\otimes\cdots\otimes\varphi_n), \, \psi_1 \otimes\cdots\otimes
\psi_n\Big\rangle = \iint
\widetilde T(p,q)\Big[\prod_{j=1}^n\varphi_1(q_j)\,\psi_1(p_j)\Big]\, dp\,dq.
$$

 \smallbreak  (2) \quad The function $\widetilde T$ satisfies the following differential inequalities on
the\break
\vskip-12pt\noindent set
$(\widetilde M \times \widetilde M \,\backslash\, \widetilde D)$. For any $(\alpha_1,\ldots, \alpha_n)$ there
is a constant $C= C_{\alpha_1,\ldots, \alpha_n}$ with the following property. Suppose
$\partial^{\alpha_j}_{p_j,q_j}$ is a differential operator of total order $|\alpha_j|$ made out of the operators
$Z_j$ and $\bar Z_j$ acting on either $p_j$ or $q_j$. Then on $\widetilde M \times \widetilde M - \widetilde
D$
$$
\Big\vert\prod_{j=1}^n\partial^{\alpha_j}_{p_j,q_j}\, \,\widetilde T(p,q)\Big\vert \leq C\, \prod_{j=1}^n d_j(p_j,q_j)^{-|\alpha_j|} \Big[\prod_{j=1}^n V_j(p_j,q_j)\Big]^{-1}.
$$
 \smallbreak  (3) \quad For each normalized bump function $\varphi_n$ supported on a ball of radius
$\delta_n$ in $M_n$ and for each point $p_n \in M_n$ there is {\it a product singular integral operator}
$\widetilde {\T}^{\varphi_n, p_n}$ (of the $(n-1)$ factor type) defined on $M_1\times \cdots \times
M_{n-1}$ so that $p_n \to \widetilde {\T}^{\varphi_n, p_n}$ is smooth and
\begin{multline*}
\big\langle \widetilde {\T}(\varphi_1\otimes \cdots \otimes \varphi_n),  \psi_1\otimes \cdots \otimes
\psi_n\big\rangle \\  = \int_{M_n} \big\langle \widetilde{\T}^{\varphi_n,p_n}(\varphi_1\otimes \cdots
\otimes \varphi_{n-1}),(\psi_1\otimes \cdots \otimes \psi_{n-1})\big\rangle \,\psi_n(p_n)\,dp_n.
\end{multline*}
Moreover, the operator $\widetilde {\T}^{\varphi_n, p_n}$ satisfies all the conditions for product singular integrals with $(n-1)$ factors, uniformly in $\varphi_n$ and $p_n$, as do all operators $\delta_n^k\,\partial^k_{p_n}\,\widetilde {\T}^{\varphi_n, p_n}$. Here $\partial^k_{p_n}$ is a differential operator of order $k$ made out of $Z_n$ 
and~$\bar Z_n$.
 \smallbreak  (4) \quad Condition (3) holds for any permutation of the indices $\{1,
\,2,\,\ldots,\,n\}$.

 \smallbreak  (5) \quad If for $1 \leq j \leq n$, $\varphi_j$ 
is {\it a normalized bump function} (in the sense of Definition 3.2.2) 
supported on a ball 
$B_j(q_j,\delta_j)\subset M_j$, then
$$
\Big\vert \prod_{j=1}^n \partial^{\alpha_j}_{p_j} \widetilde
\T(\varphi_1\otimes \cdots \otimes \varphi_n)\Big\vert \lesssim
\prod_{j=1}^n \delta_j^{-|\alpha_j|}
$$
where the constants implied by the symbol $\lesssim$ are uniform.
 \smallbreak (6) \quad Properties (1) through (5) are also satisfied for all possible transposes of $\widetilde
T$; i.e.\ those operators which arise by interchanging some collection of $p_j$ and $q_j$.
\end{definition}
 
{\it Remark.} If $T_j$ is an NIS operator 
smoothing of order zero on $M_j$ for $1\leq j \leq n$, 
then the operator $T_1\otimes \cdots\otimes T_n$ is a product singular integral operator on $\widetilde M$.
\Enddemo

The main result on product singular integral operators, which is proved in \cite{NaSt00.3}, is the following.

\begin{theorem}\label{ProductSingularIntegralTheorem}
 Suppose that $\widetilde\T$ is a product singular integral operator on $\widetilde M$ in the sense of
Definition {\rm \ref{D:5.5.1}.}
 Then $\widetilde\T$ is bounded on $L^p(\widetilde M)$ for $1 < p < +\infty$.
\end{theorem}

{\it Remark.} In addition to the notion of product singular integral operators defined in Section \ref{SS:NIS},
there are also NIS operators smoothing of order zero with respect to the metric $d_\Sigma$. Both these
families of operators might be called {\it singular integrals} on $\widetilde M$. In fact, NIS operators of
order zero relative to the metric $d_\Sigma$ are also product NIS operators. The corresponding result for
Euclidean singular integrals is noted in \cite[Remark 2.1.6]{NaRiSt00}.

\section{Relative fundamental solutions for $\bx_b$ on $ M_1\times \cdots \times M_n$}
\label{S:Fundamental Solutions} 

We now construct two relative fundamental solutions $\widetilde\K_J$ and $\widetilde\N_J$ for each of the $2^n$ Kohn-Laplacian operators $\bx_J$ on the product $\widetilde M$. Recall that for each increasing tuple $J$ of integers between $1$ and $n$, $\bx_J = \sum_{j=1}^n \bx_j^{J(j)}$. On $\widetilde M$ there is essentially no difference between these operators. Thus we shall fix $J$ and abbreviate our notation so that we write $ \bx_J = \bx_b = \sum_{j=1}^n \bx_j$ where each $\bx_j$ is either $Z_j\,\bar Z_j$ or $\bar Z_j\,Z_j$. As usual, $H_j$ and $S_j$ are the heat kernel and projection associated to the particular choice of $\bx_j$. We will write the corresponding relative fundamental solutions for $\bx_b$ as $\widetilde\K$ and $\widetilde \N$.

\Subsec{The relative fundamental solutions $\widetilde{\mathcal K}$ and $\widetilde{\mathcal N}$}
Following the discussion in Section \ref{SS:Outline}, we make the following definitions.

\begin{definition}\label{D:2.3.1} For $(p,q) \in \widetilde M\times \widetilde M$ set
\begin{align}
\widetilde K(p,q) &= \int_0^{\infty} \Big(\prod_{j=1}^n H_j(s,p_j,q_j) - \prod_{j=1}^n S_j(p_j,q_j)\Big)\,ds\label{E:2.1.1.a}
,\\
\widetilde N(p,q) &= \int_0^\infty \prod_{j=1}^n\Big(H_j(s,p_j,q_j) - S_j(p_j,q_j)\Big)\,ds.\label{E:2.2.1.b}
\end{align}
Also for any proper subset $A \subset \{1,\ldots,n\}$, set
\begin{align}
\widetilde N_A(p,q) &= \int_0^\infty \prod_{j\in A} \Big(H_j(s,p_j,q_j) -
S_j(p_j,q_j)\Big)\,ds,\label{E:2.2.1.c}\\
\widetilde S_A(p,q) &= \bigotimes_{j\notin A} S_j(p_j,q_j).
\end{align}
If $A = \emptyset$ then $ \widetilde N_A = 0$, and we set $S = \widetilde S_A = \prod_{j=1}^n S_j$. Note that $\widetilde N_A$ only involves variables $(p_j,q_j)$ with $j \in A$, and $\widetilde S_A$ only involves variables $(p_j,q_j)$ with $j \notin A$. 
\end{definition}

\begin{proposition} The integrals {\rm (\ref{E:2.1.1.a})}
 through {\rm (\ref{E:2.2.1.c})} and $\widetilde S_A$ define distributions
 on $\widetilde M\times \widetilde M$.
For $\varepsilon > 0$ and $\emptyset \neq A \subset \{1,\ldots,n\}${\rm ,} the integral
\begin{align*}
\widetilde N_{A,\varepsilon}(p,q) &= \int_\varepsilon^\infty \prod_{j\in
A} \Big(H_j(s,p_j,q_j) - S_j(p_j,q_j)\Big)\,ds
\end{align*}
converges absolutely to an infinitely 
differentiable function with bounded derivatives. As distributions{\rm ,} $\widetilde N_{A,\varepsilon} \to \widetilde
N_{A}$ as $\varepsilon \to 0$.
\end{proposition}

\Proof  We have
$$
\widetilde N_{A,\varepsilon}(p,q) = \int_\varepsilon^\infty \prod_{j\in
A} G_j(s,p_j,q_j) \,ds.
$$
The estimates of Theorem \ref{T:1.1.1} and equation (\ref{E:2.3.1e}) show that $\widetilde N_{A,\varepsilon}$ is smooth and bounded with bounded derivatives as long as $\varepsilon > 0$. Also, expanding $\prod_{j\in A}(H_j-S_j)$, 
we see that
\begin{eqnarray*}
\widetilde N_{A,\varepsilon}(p,q) &=& \int_1^\infty \prod_{j\in A}
G_j(s,p_j,q_j)\,ds + \int_{\varepsilon}^1 \prod_{j\in A}
H_j(s,p_j,q_j)\,ds \\&&+(-1)^{|A\backslash B|} \sum_{B\subsetneqq A}\Big[
\int_\varepsilon^1 \prod_{j\in B} H_j(s,p_j,q_j)\,ds \Big]\bigotimes
\prod_{j\in A\backslash B} S_j(p_j,q_j).
\end{eqnarray*}
It follows from Theorem \ref{T:1.2} that $\int_\varepsilon^1 \prod_{j\in B} H_j(s,p_j,q_j)\,ds$ is a family of distributions which converges as $\varepsilon \to 0$. Hence the integral (\ref{E:2.2.1.c}) converges to a distribution. This establishes the statements about the distributions $\widetilde N_{A,\varepsilon}$ and $\widetilde N_A$. Next $\widetilde S_A$ is a distribution since it is a tensor product of distributions. The statement that the integral in
equation (\ref{E:2.1.1.a}) defines a distribution follows\footnote{This is a consequence of the Taylor expansion
 \begin{equation*}
 \prod_{j=1}^n x_j = \prod_{j=1}^n y_j + \underset{A
 \subsetneqq \{1,\ldots,n\}}{\sum} \big[\prod_{j\in A}
 y_j\,\prod_{k\notin A}(x_k-y_k)\big].
 \end{equation*}} 
from the identity 
\begin{multline*}
\int_\varepsilon^\infty \Big(\prod_{j=1}^n H_j(s,p_j,q_j) -
\prod_{j=1}^n S_j(p_j,q_j)\Big)\,ds \\
= \widetilde N_\varepsilon(p,q)
 + \sum_{A}
\widetilde N_{A,\varepsilon}(p,q)\otimes \widetilde \S_A(p,q).
\end{multline*}
Since the right-hand side defines a distribution which has a limit as $\varepsilon \to 0$, the same is true of the left-hand side. This completes the proof.
\hfq

\begin{lemma}\label{P:2.3.2}
The distributions $\widetilde K${\rm ,}
 $\widetilde N${\rm ,} and $\widetilde N_{A}$ induce operators $\widetilde\K${\rm ,}
 $\widetilde \N${\rm ,} and $\widetilde
\N_{A}$ on $\widetilde M$ which satisfy
\begin{equation}\label{E:First0}
\displaystyle \widetilde\K = \widetilde \N + \sum_{A} \widetilde \N_{A} \otimes \widetilde \S_A
\end{equation}
where the sum is take over all proper{\rm ,}
 nonempty subsets $A$ of $\{1,\ldots,n\}$. Moreover{\rm ,} $\widetilde\K$ and $\widetilde \N$ are relative
fundamental solutions of $\bx_b$ in the sense that
\begin{align}
 \widetilde\K\,\bx_b &= \bx_b\,\widetilde\K = I - \bigotimes_{j=1}^n \S_j;\label{E:First}\\ 
 \widetilde \N\,\bx_b &= \bx_b\,\widetilde \N = \bigotimes_{j=1}^n\left(I - \S_j\right).\label{E:First2}
\end{align}
\end{lemma}

\Proof  The identity (\ref{E:First0}) is again just a consequence of the formula for the Taylor expansion for the
function $\prod_{j=1}^n x_j$. To check identity (\ref{E:First}), note that $\bx_j\,\S_j \equiv 0$, and $\bx_j\,H_j(s,p_j,q_j) = - \frac{\partial H_j}{\partial s}(s,p_j,q_j)$. Thus, working with distributions, we have
\begin{equation*}
\begin{split}
(\sum_{j=1}^n \bx_j)\,K(p,q) &= -\int_0^\infty \frac{\partial (\prod_{j=1}^n H_j)}{\partial s}(s,p,q)\,ds= \prod_{j=1}^n\delta(p_j,q_j) - \prod_{j=1}^n S_j(p_j,q_j).
\end{split}
\end{equation*}
A similar argument gives identity (\ref{E:First2}).
\Endproof\vskip4pt  

{\it Remark.} Note that $\widetilde\K$ inverts $\bx_b$ modulo the projection onto the intersection of the
null spaces of the operators $\{\bx_j\}$ which is just the null space of $\bx_b$. On the other hand, $\widetilde
\N$ inverts $\bx_b$ modulo the (much larger) subspace which is the direct sum of the null space of the
operators $\{\bx_j\}$. Thus the operator $\widetilde\K$ is the {\it natural} relative fundamental solution
which inverts the operator except on the smallest possible subspace. On the other hand, as we shall see, the
operator $\widetilde \N$ has better regularity properties than the operator $\widetilde\K$. The identity
(\ref{E:First0}) in Lemma
\ref{P:2.3.2} provides the link between the two.

\Subsec{Analysis of the distributions $\widetilde N$}
\label{SS:widetilde K} In this section we show that we have global maximal hypoelliptic estimates for the
relative fundamental solution $\widetilde\N$ by showing that any two good derivatives of $\widetilde N$
give a product singular integral in the sense of Definition \ref{D:5.5.1}. It then follows from Theorem
\ref{ProductSingularIntegralTheorem} that such operators are bounded on $L^p(\widetilde M)$ for
$1<p<\infty$.

\begin{theorem}\label{T:5.1.1.1} Let $Q(Z, \bar Z)$ be any quadratic expression in the vector fields $\{Z_1,\,\bar Z_1, \ldots, Z_n,\,\bar Z_n\}$. Then $Q(Z, \bar Z)\,\widetilde \N$ is a product singular integral operator on $\widetilde M$ and consequently is bounded on $L^p(\widetilde M)$ for all $1 < p < +\infty$. 
\pagebreak
In 
 particular{\rm ,} there is a constant $C = C_{p,Q}$ so that if $\bx_b[u] = g$ and if\break
$\bigotimes_{j=1}^n(I- \S_j)[u] = u${\rm ,} then for $1 <p < \infty$
$$
\norm Q(Z,\bar Z)[u]\norm_{L^p(\widetilde M)} \leq C\,\norm
g\norm_{L^p(\widetilde M)}.
$$
\end{theorem}
\vskip8pt

Before presenting the proof of Theorem \ref{T:5.1.1.1}, we make some remarks, and state two useful inequalities.

\demo{Remark $1$} The condition that $\bigotimes_{j=1}^n(I- \S_j)[u] = u$ is equivalent to the statement
that
$u$ is orthogonal to each of the null spaces $\big\{u\,\big\vert\,\bx_j[u] = 0\big\}$, $1 \leq j \leq n$.

\demo{Remark $2$} The estimate in Theorem \ref{T:5.1.1.1} follows from the boundedness of the operator
$Q(Z, \bar Z)\N$ and the identity (\ref{E:First2}).

\demo{Remark $3$} Let $\chi \in \CC^\infty_0(\R)$ 
with $\chi(t) \equiv 1$ for $|t| \leq \frac{1}{2}$ and 
$\chi(t) \equiv 0$ for $|t| \geq 1$. 
Let $\rho_j$ be a  ``regularized distance'' on $M_j$, that is, a 
smooth function on 
$M_j\times M_j$ such that $\rho_j(p_j,q_j) \approx d_j(p_j,q_j)$ and 
$\big\vert X^\alpha \rho_j(p_j,q_j) 
\big\vert \lesssim d_j(p_j,q_j)^{1-|\alpha|}$, for any derivative
$X^\alpha$ of order $| \alpha |$ in the vector fields $X_j$ and
$X_{n+j}$ acting on either $p_j$ or $q_j$. 
(The existence of such distances is
established in \cite{NaSt00.2}.) For each $R>0$ define
$$
\chi_R(p,q) = \prod_{j=1}^n
\chi\left(\frac{\rho_j(p_j,q_j)}{R}\right).
$$
Our proof of Theorem \ref{T:5.1.1.1} will actually show that if we consider the kernel $\widetilde N_\varepsilon(p,q)\,\chi_R(p,q)$, then any two good derivatives composed with the corresponding operator
yield a product NIS operator on $\widetilde M$, with constants independent of $\varepsilon$ and $\R$. For
$\varepsilon > 0$ and $R < +\infty$, the kernel $\widetilde N_\varepsilon(p,q)\,\chi_R(p,q)$ is bounded and
has compact support. This observation will be important when we use transference results in Section
\ref{S:Descending} below to obtain information about operators on the decoupled boundary $M$.
 \Enddemo

The following two elementary estimates will be used frequently in what follows. 

\begin{proposition}\label{Prop3.2.4}  {\rm (a)}   If $F$ is a positive{\rm ,}
 monotone decreasing function on
$(0,\infty)${\rm ,}
 and if there exists $\varepsilon > 0$ such that $F(2t) \leq 2^{-1-\varepsilon}\,F(t)${\rm ,} then there is
a constant $C$ depending on $\varepsilon$ so that for all $a > 0${\rm ,}
\begin{equation*}
\int_a^\infty F(t)\,dt \leq C\, a\, F(a).
\end{equation*}

\smallskip

  {\rm (b)}   If $F$ is a positive{\rm ,} monotone decreasing function on $(0,\infty)${\rm ,}
 and if there exists $\varepsilon
> 0$ such that $F(t/2) \leq 2^{-1+\varepsilon}\,F(t)${\rm ,}
 then there is a constant $C$ depending on $\varepsilon$
so that for all $a > 0${\rm ,}
$$
\int_0^a F(t)\,dt \leq C\, a\, F(a).
$$
\end{proposition}

Now we turn to the proof of Theorem \ref{T:5.1.1.1}. For simplicity of exposition, we confine ourselves to the proof in the case that $n = 2$. The general situation only involves additional notational difficulty. Recall that in this two-dimensional case, the distributional kernel for $\widetilde \N$ is given by
$$
\widetilde N(p_1, p_2, q_1, q_2) = \int_0^\infty \prod_{j=1}^2
G_j(s,p_j,q_j)\,ds.
$$
We only need to show that the operator $Q(Z, \bar Z)\N$ is a product singular integral operator in the sense of Definition \ref{D:5.5.1}. In what follows, we establish the differential inequalities required for condition (2) in Lemma \ref{T:2.1.1}, the estimates required for condition (3) in Lemma \ref{T:3.1.1}, and the estimates for condition (5) in Lemma \ref{T:3.2.3}.

\begin{lemma}\label{T:2.1.1} On the set $(M_1\times M_2)\times (M_1\times M_2) - \widetilde D${\rm ,}
 the distribution $\widetilde N$ is given by integration against an infinitely differentiable function. For $j =
1, 2${\rm ,}
 let $\partial_j^{\alpha_j}$ be a derivative of order $|\alpha_j|$ in the vector fields $Z_j$ or $\bar Z_j$
acting either on the variables $p_j$ or $q_j$. If $p_1 \neq q_1$ and $p_2 \neq q_2${\rm ,} the integral
$$
\partial^{\alpha_1}_1\,\partial^{\alpha_2}_2\, \widetilde N(p_1,q_1,p_2,q_2) = \int_0^\infty \partial^{\alpha_1}_1 G_1(s,p_1,q_1)\,\partial^{\alpha_2}_2G_2(s,p_2,q_2)
\,ds 
$$
converges absolutely, and there is a constant $C= C_{\alpha_1,\alpha_2}$ such that
\begin{eqnarray}\label{E:4.2.1}
&&\Big\vert \partial_1^{\alpha_1} \partial_2^{\alpha_2} \widetilde N(p_1,q_1,p_2,q_2)\Big\vert
\\&\leq &C\,\left[\frac{d_1(p_1,q_1)^{-|\alpha_1|}\,d_2(p_2,q_2)^{-|\alpha_2|}
\left[\min\{d_1(p_1,q_1),d_2(p_2,q_2)\}\right]^2} {V_1\big(p_1,d_1(p_1,q_1)\big)\,
V_2\big(p_2,d_2(p_2,q_2)\big)}\right]. \nonumber
\end{eqnarray}
\end{lemma}

\Proof  Fix $p_1\neq q_1$ and $p_2 \neq q_2$, and assume without loss of generality that $0 < d_1(p_1,q_1) \leq d_2(p_2,q_2)$. Write
\begin{align*}
&\int_0^\infty \partial^{\alpha_1}_1 G_1(s,p_1,q_1)\, \partial^{\alpha_2}_2G_2(s,p_2,q_2)\,ds
\\ &= \left[ \int_0^{d_1(p_1,q_1)^2}+\int_{d_1(p_1,q_1)^2}^{d_2(p_2,q_2)^2} +
\int_{d_2(p_2,q_2)^2}^\infty\right]\left[\partial^{\alpha_1}_1
G_1(s,p_1,q_1)\, \partial^{\alpha_2}_2G_2(s,p_2,q_2)\right]\,ds
\\
&=\quad I \quad+\quad II \quad +\quad III.
\end{align*}
According to Theorem \ref{T:1.1.1},  we can estimate integral $I$ by
$$
\big\vert I \big\vert \leq C\,\da^2\,\da^{-|\alpha_1|}\Va^{-1}
\db^{-|\alpha_2|}\Vb^{-1}.
$$
For integral $II$ we use Theorem \ref{T:1.1.1} and Proposition \ref{Prop3.2.4}(a) to obtain
\begin{eqnarray*}
\big\vert II \big\vert &\leq& C\, \db^{-|\alpha_2|}\Vb^{-1}
\int_{d_1(p_1,q_1)^2}^\infty s^{-|\alpha_1|/2}\, \left\vert
B_1(p_1,\sqrt s)\right\vert^{-1}\,ds
\\
&\leq& C\,\da^2\,\da^{-|\alpha_1|}\Va^{-1}
\db^{-|\alpha_2|}\Vb^{-1}
\end{eqnarray*}
since $F(t) = t^{-1-|\alpha_1|/2} \Lambda_1\big(p_1,\sqrt
t)\big)^{-1}$ satisfies $F(2t) \leq 2^{-2}\,F(t)$.  Finally, for integral $III$, we again use
 Proposition \ref{Prop3.2.4}(a) and obtain
\begin{eqnarray*}
\big\vert III \big\vert &\leq& C \int_{d_2(p_2,q_2)^2}^\infty
s^{-(|\alpha_1|+|\alpha_2|)/2} \left\vert B_1(p_1,\sqrt
s)\right\vert^{-1}\, \left\vert B_2(p_2,\sqrt
s)\right\vert^{-1}\,ds\\ &\leq&
C\,\db^{2-(|\alpha_1|+|\alpha_2|)}\, \left\vert
B_1(p_1,d_2(p_2,q_2))\right\vert^{-1}\, \left\vert
B_2(p_2,d_2(p_2,q_2))\right\vert^{-1}\,\\ &\leq&
C\,\db^{-|\alpha_1|}\,\Lambda_1(p_1,\db)^{-1}
\,\db^{-|\alpha_2|}\,V_2(p_2,\db)^{-1}\\ &\leq&
 C\,\da^2\,\da^{-|\alpha_1|}\Va^{-1}
\db^{-|\alpha_2|}\Vb^{-1}.
\end{eqnarray*}
Combining the estimates for $I$, $II$, and $III$ completes the proof since $\da\break = \min\left\{\da,
\db\right\}$.
\Endproof\vskip4pt

We next study cancellation properties of the distribution $\widetilde N$ which are expressed
 by the action of the operator $\widetilde \N$ on normalized bump functions on one factor, say $M_1$.

\begin{lemma} \label{T:3.1.1} Let $\varphi$ be a normalized bump function supported on
$B_1(p_1,\delta)\break \subset M_1$. Let $\partial_j^{\alpha_j}$ be a derivative of order $|\alpha_j|$ in the
vector fields $Z_j$ or $\bar Z_j$ acting either on the variables $p_j$ or $q_j$. There is a constant $C =
C(\alpha_1,\alpha_2)$ so that if $p_2 \neq q_2${\rm ,}
\begin{equation*}
\begin{split}
\Big\vert \int_0^\infty\Big[\int_{M_1}\partial_1^{\alpha_1} G_1(s,p_1,q_1)\,&\varphi(q_1)\,dq_1\Big]\,\partial_2^{\alpha_2}
G_2(s,p_2,q_2)\,ds\,\Big\vert \\ &\leq
C\,\delta^{-|\alpha_1|}\,\frac{d_2(p_2,q_2)^{-|\alpha_2|}}
{V_2(p_2,q_2)}\, \min\{\delta^2,d_2(p_2,q_2)^2\}.
\end{split}
\end{equation*}
\end{lemma}

\Proof  We have the following estimates:
\begin{align*}\tag{a}
\left\vert \int_{M_1}\partial^{\alpha_1}_1G_1(s,p_1,q_1)\,\varphi(q_1)\,dq_1\right\vert
&\leq C
\begin{cases}
\delta^{-|\alpha_1|}&\text{if $s\leq \delta^2$},\\\\
s^{-|\alpha_1|/2}\frac{V_1(p_1,\delta)}{V_1(p_1,\sqrt s)} &
\text{if $s \geq \delta^2$}.
\end{cases}\\\\
\tag{b} \left\vert \partial^{\alpha_2}_2G_2(s,p_2,q_2)\right\vert &\leq C
\begin{cases}
\frac{d_2(p_2,q_2)^{-|\alpha_2|}}{V_2(p_2,q_2)} &\text{if $s\leq
d_2(p_2,q_2)^2$},\\\\ \frac{s^{-|\alpha_2|/2}}{V_2(p_2,\sqrt s)} &
\text{if $s \geq d_2(p_2,q_2)^2$}.
\end{cases}
\end{align*}
(The first part of assertion (a) follows from Theorem \ref{T:1.2}. The second part 
\pagebreak
of (a) as well as assertion (b) follow from Theorem \ref{T:1.1.1}.) Write

\begin{eqnarray*} \noalign{\vskip-16pt}
&&\hskip-12pt \int_0^\infty\left[\int_{M_1}\partial^{\alpha_1}_1
G_1(s,p_1,q_1)\,\varphi(q_1)\,dq_1\right]\,\partial^{\alpha_2}_2
G_2(s,p_2,q_2)\,ds
\\&&\quad=\left[ \int_0^{\delta^2}+ \int_{\delta^2}^\infty\right]
\left[\int_{M_1}\partial^{\alpha_1}_1
G_1(s,p_1,q_1)\,\varphi(q_1)\,dq_1\right]\,\partial^{\alpha_2}_2
G_2(s,p_2,q_2)\,ds
\\& &
\quad= A + B.
\end{eqnarray*}
In integral A, $s \leq \delta^2$ so according to (a) we obtain
$$
\big\vert A \big\vert \leq
C\,\delta^{-|\alpha_1|}\,\int_0^{\delta^2} \big\vert
\partial_2^{\alpha_2} G_2(s,p_2,q_2)\big\vert\,ds.
$$
The analysis now depends on the relative sizes of $d_2(p_2,q_2)$ and $\delta$.

\demo{Case {\rm 1:} $d_2(p_2,q_2) \leq \delta$}  In this
case, using estimate (b) we get
\begin{eqnarray*}
\big\vert A \big\vert &\leq &C\,
\delta^{-|\alpha_1|}\left[\int_0^{d_2(p_2,q_2)^2}
\frac{d_2(p_2,q_2)^{-|\alpha_2|}}{V_2(p_2,q_2)}\,ds +
\int_{d_2(p_2,q_2)^2}^\infty
\frac{|s|^{-|\alpha_2|/2}}{V_2(p_2,\sqrt s)}\,ds\right]\\ &\leq&
C\, \delta^{-|\alpha_1|}\frac{d_2(p_2,q_2)^{-|\alpha_2|}}
{V_2(p_2,q_2)}\,d_2(p_2,q_2)^2\\ &=& C\, \delta^{-|\alpha_1|}\,
\frac{d_2(p_2,q_2)^{-|\alpha_2|}}{V_2(p_2,q_2)}\,
\min\{\delta^2,d_2(p_2,q_2)^2\}.
\end{eqnarray*}

\demo{Case {\rm 2:} $d_2(p_2,q_2) \geq \delta$}  In this
case, the region of $s$ integration only involves $s \leq
d_2(p_2,q_2)^2$, and so using the first case of estimate (b), we
obtain
\begin{eqnarray*}
\big\vert A \big\vert &\leq& C\,
\delta^{-|\alpha_1|}\,\frac{d_2(p_2,q_2)^{-|\alpha_2|}}{V_2(p_2,q_2)}\,
\delta^2\\ &=&
\delta^{-|\alpha|}\,\frac{d_2(p_2,q_2)^{-|\alpha_2|}}{V_2(p_2,q_2)}\,
\min\{\delta^2,d_2(p_2,q_2)^2\}.
\end{eqnarray*}

Next consider the integral $B$. Here the range of integration is
$s\geq \delta^2$ and so by the first estimate in (a), we have
$$
\big\vert B \big\vert \leq C\,\int_{\delta^2}^\infty
s^{-|\alpha_1|/2}\frac{V_1(p_1,\delta)}{V_1(p_1,\sqrt
s)}\,\partial^{\alpha_2}_2G_2(s,p_2,q_2)\,ds.
$$
Again, the analysis depends on the relative size of $d_2(p_2,q_2)$
and $\delta$.

\demo{Case {\rm 1:} $d_2(p_2,q_2) \leq \delta$}  In this
case we integrate where $s \geq d_2(p_2,q_2)^2$, and according to
(b) we have
\begin{equation*}
\begin{split}
\big\vert B \big\vert &\leq C\,\int_{\delta^2}^\infty
s^{-(|\alpha_1|+|\alpha_2|)/2}\,\frac{V_1(p_1,\delta)}{V_1(p_1,\sqrt
s)\,\,V_2(p_2,\sqrt s)}\,ds \\ &\leq C\, \int_{\delta^2}^{\infty}
\frac{|s|^{-(|\alpha_1|+|\alpha_2|)/2}}{V_2(p_2,\sqrt s)}\,ds\\
&\leq C\,\delta^{2-|\alpha_1|-|\alpha_2|}\,V_2(p_2,\delta)^{-1}\\
&\leq C\,
\delta^{-|\alpha_1|-|\alpha_2|}\,\Lambda_2(p_2,\delta)^{-1}\\
&\leq C\,
\delta^{-|\alpha_1|}\,\frac{d_2(p_2,q_2)^{-|\alpha_2|}}{V_2(p_2,q_2)}\,
\min\{\delta^2,d_2(p_2,q_2)^2\}.
\end{split}
\end{equation*}

\demo{Case {\rm 2:} $d_2(p_2,q_2) \geq \delta$}  In this
case, we write
\begin{equation*}
\begin{split}
B &= \left[\int_{\delta^2}^{d_2(p_2,q_2)^2} +
\int_{d_2(p_2,q_2)^2}^\infty \right] \left(s^{-|\alpha_1|/2}
\frac{V_1(p_1,\delta)}{V_1(p_1,\sqrt
s)}\,\partial^{\alpha_2}_2G_2(s,p_2,q_2)\,\right)ds\\\\& = I + II.
\end{split}
\end{equation*}

\medskip

In integral $I$, we integrate where $s \leq d_2(p_2,q_2)^2$ and so
by the first estimate in (b) we have
\begin{eqnarray*}
\big\vert I \big\vert &\leq& C\,
\frac{d_2(p_2,q_2)^{-|\alpha_2|}}{V_2(p_2,q_2)}\,\int_{\delta^2}^\infty
s^{-|\alpha_1|/2}\,\frac{V_1(p_1,\delta)}{V_1(p_1,\sqrt s)}\,ds\\
&\leq &C\,
\delta^{-|\alpha_1|}\,\frac{d_2(p_2,q_2)^{-|\alpha_2|}}{V_2(p_2,q_2)}\,\delta^2\\
&=& C\,
\delta^{-|\alpha_1|}\,\frac{d_2(p_2,q_2)^{-|\alpha_2|}}{V_2(p_2,q_2)}\,
\min\{\delta^2,d_2(p_2,q_2)^2\}.
\end{eqnarray*}
Finally, in integral $II$, we integrate where $s \geq
d_2(p_2,q_2)^2$ and so by the second estimate in (b) we have
\begin{equation*}
\begin{split}
\big\vert II \big\vert &\leq C\, \int_{d_2(p_2,q_2)^2}^\infty
s^{-(|\alpha_1|+|\alpha_2|)/2}\frac{V_1(p_1,\delta)}{V_1(p_1,\sqrt
s)}\,V_2(p_2,\sqrt s)^{-1}\,ds \\
&\leq C\,
\frac{d_2(p_2,q_2)^{2-|\alpha_1|-|\alpha_2|}}{V_2(p_2,d_2(p_2,q_2))}
\, \frac{\delta^2\,\Lambda_1(p_1,\delta)}
{d_2(p_2,q_2)^2\,\Lambda_1\big(p_1,d_2(p_2,q_2)\big)}\\ &\leq C\,
\frac{d_2(p_2,q_2)^{-|\alpha_1|-|\alpha_2|}}
{V_2(p_2,d_2(p_2,q_2))}\,\delta^2\\ &= C\,
\delta^{-|\alpha_1|}\,\frac{d_2(p_2,q_2)^{-|\alpha_2|}}{V_2(p_2,q_2)}\,
\min\{\delta^2,d_2(p_2,q_2)^2\}.
\end{split}
\end{equation*}
This completes the proof of Lemma \ref{T:3.1.1}.\Endproof\vskip4pt  

\medskip

We finally need to study the action of the distribution $\widetilde
N$ on pairs of bump functions.

\begin{lemma}\label{T:3.2.3} Suppose $\varphi_j$ is a normalized
bump function supported on the ball $B_j(p_j,\delta_j)\subset
M_j$. There is a constant $C$ depending on $\alpha_j$ but
independent of $\delta_j$ so that
\begin{equation*}
\begin{split}
\Bigg\vert \int_0^\infty \int_{M_1} \int_{M_2}
\partial^{\alpha}_1G_1(s,p_1,q_1)\,&\partial^{\beta}_2G_2(s,p_2,q_2)\,\varphi_1(q_1)
\,\varphi_1(q_2)\,dq_1\,dq_2\,ds\Bigg\vert\\ &\leq C\,
\delta_1^{-|\alpha_1|}\,\delta_2^{-|\alpha_2|}\,
\min\{\delta_1^2,\delta_2^2\}.
\end{split}
\end{equation*}
\end{lemma}
\vglue8pt

\Proof  Suppose without loss of generality that
$\delta_1 \leq \delta_2$. Write
\begin{equation*}
\int_0^\infty = \int_0^{\delta_1^2} +
\int_{\delta_1^2}^{\delta_2^2} + \int_{\delta_2^2}^\infty = I + II
+ III.
\end{equation*}
We use the estimates
\begin{equation*}
\begin{split}
\Bigg\vert\int_{M_j} \partial^{\alpha_j}
G_j(s,p_j,q_j)\,\varphi_j(q_j)\,dq_j \Bigg\vert \leq C
\begin{cases}
\delta^{-|\alpha_j|}&\text{if $s\leq \delta_j^2$},\\\\
s^{-|\alpha_j|/2}\frac{V_j(p_j,\delta_j)}{V_j(p_j,\sqrt s)} &
\text{if $s \geq \delta_j^2$}.
\end{cases}
\end{split}
\end{equation*}
For integral $I$, we have the estimate
\begin{equation*}
\big\vert I \big\vert \lesssim
\delta_1^{-|\alpha_1|}\,\delta_2^{-|\alpha_2|}\,\delta_1^2 =
\delta_1^{-|\alpha_1|}\,\delta_2^{-|\alpha_2|}\,
\min\{\delta_1^2,\delta_2^2\}.
\end{equation*}
For integral $II$ we have the estimate
\begin{equation*}
\begin{split}
\big\vert II \big\vert &\lesssim
\delta_2^{-|\alpha_2|}\int_{\delta_1^2}^\infty s^{-|\alpha_1|/2}\,
\frac{V_1(p_1,\delta_1)}{V_1(p_1,\sqrt s)}\,ds\\ &\lesssim
\delta_1^{-|\alpha_1|}\,\delta_2^{-|\alpha_2|}\,\delta_1^2 =
\delta_1^{-|\alpha_1|}\,\delta_2^{-|\alpha_2|}\,
\min\{\delta_1^2,\delta_2^2\}.
\end{split}
\end{equation*}
Finally, for integral $III$ we have the estimate
\begin{equation*}
\begin{split}
\big\vert III \big\vert &\lesssim \int_{\delta_2^2}^\infty
s^{-(|\alpha_1|+|\alpha_2|)/2}
\frac{V_1(p_1,\delta_1)}{V_j(p_1,\sqrt s)}
\frac{V_2(p_2,\delta_2)}{V_j(p_2,\sqrt s)}\,ds\\ &\lesssim
\delta_2^{-|\alpha_1|+|\alpha_2|}\,
\frac{\delta_1^2\,\Lambda_1(p_1,\delta_1)}
{\delta_2^2\,\Lambda_j(p_1,\delta_2)}\delta_2^2\\ &\lesssim
\delta_1^{-|\alpha_1|}\,\delta_2^{-|\alpha_2|}\,\delta_1^2\\ &=
\delta_1^{-|\alpha_1|}\,\delta_2^{-|\alpha_2|}\,
\min\{\delta_1^2,\delta_2^2\}.
\end{split}
\end{equation*}
This completes the proof of Lemma \ref{T:3.2.3}, and consequently
the proof of Theorem \ref {T:5.1.1.1}.\hfq

\Subsec{Analysis of the distribution $\widetilde K$}
\label{SS:Analysis of K} Equation (\ref{E:First0}) gives a decomposition of the {\it
operator} $\widetilde\K$. It will also be important to have a
different kind of decomposition of the {\it kernel} $\widetilde K$
which is given by
\begin{equation*}
\begin{split}
\widetilde K(p,q) = \int_0^\infty \Big[\prod_{j=1}^n
H_j(s,p_j,q_j) - \prod_{j=1}^n S_j(p_j,q_j)\Big]\,ds.
\end{split}
\end{equation*}

\begin{theorem}\label{T:MainEstimate} There are distributions
$\widetilde K_0${\rm ,} $\widetilde K_\infty${\rm ,} and for each nonempty
proper subset $A \subset \{1,\ldots,n\}$ a distribution
$\widetilde K_A$ on $\widetilde M \times \widetilde M$ so that
away from the product diagonal $\widetilde D$ on $\widetilde
M\times \widetilde M${\rm ,}
$$
\widetilde K(p,q) = \widetilde K_0(p,q) +
\sum_A\Big(\prod_{j\notin A} S_j(p_j,q_j) \Big) \widetilde
K_A(p,q) + \Big(\prod_{j=1}^n S_j(p_j,q_j)\Big)\widetilde
K_\infty(p,q).
$$
Moreover $\widetilde K_0${\rm ,} $\widetilde K_A${\rm ,} and $\widetilde
K_\infty$ are locally integrable functions{\rm ,} and the size of their
derivatives can be estimated by the control metric $d_\Sigma$.
Thus for all derivatives $X^\alpha$ of total order $|\alpha|$ in
$X_1, \ldots,X_{2n}$ acting on the variables $p_j$ or $q_j${\rm ,} 
\begin{equation}\label{E:6.1.1a}
\begin{split}
\Big\vert X^{\alpha}\widetilde K_{0}(p,q)\Big\vert &\leq
C_{\alpha}\,\frac{d_\Sigma(p,q)^{2-|\alpha|}}{\prod_{j=1}^n
d_\Sigma(p,q)^2 \Lambda_j(p_j,d_\Sigma(p,q))};\\ \Big\vert
X^{\alpha} \widetilde K_{A}(p,q)\Big\vert &\leq
C_{\alpha}\,\frac{d_\Sigma(p,q)^{2-|\alpha|}}{ \prod_{j\in
A}d_\Sigma(p,q)^2 \Lambda_{j}\big(p_j;d_\Sigma(p,q)\big)};\\
\Big\vert X^{\alpha}\widetilde K_{\infty}(p,q) \Big\vert &\leq
C_{\alpha}\,d_\Sigma(p,q)^{2-|\alpha|}.
\end{split}
\end{equation}
\end{theorem}
\vskip8pt

To prove the estimates in Theorem \ref{T:MainEstimate}, we
decompose $H_{j}(s,x,y)$ into two parts. One part is supported
close to the diagonal and is singular there. The other part is
smooth everywhere.

\begin{definition} Let $\rho_j$ be {\it a regularized distance function
on} $M_j$ as defined in Remark 3 of Section 4.2. 
Let $\chi \in C^{\infty}(\R)$ satisfy $\chi(t)\equiv 1$ if $t\leq
\frac{1}{2}$ and $\chi(t)\equiv 0$ if $t\geq 1$. Let
\begin{equation*}
\chi_j(s,p_j,q_j) = \chi \left(\frac{\rho_j(p_j,q_j)^2}{s}\right)
\end{equation*}
and set
\begin{align*}
\widetilde S_{j}(s,p_{j},q_{j}) &=
\chi_j(s,p_j,q_j)\,S_{j}(p_{j},q_{j});\\ \Phi_{j}(s,p_j,q_{j}) &=
H_{j}(s,p_j,q_{j}) - \widetilde S_{j}(s,p_j,q_{j})
\end{align*}
\end{definition}

\begin{proposition}\label{P:6.2.2.1}
$H_{j}(s,p_j,q_{j}) = \Phi_{j}(s,p_j,q_{j}) + \widetilde
S_{j}(s,p_j,q_{j})$. Moreover
\begin{equation*}
\begin{split}
\Phi_{j}(s,p_j,q_{j}) &= \begin{cases} G_{j}(s,p_j,q_{j})&
\text{when \quad $s \geq 2\,\rho_j(p_j,q_{j})^{2}$},\\
H_{j}(s,p_j,q_{j})& \text{when \quad $s \leq
\rho_j(p_j,q_{j})^{2}$}.
\end{cases}\\
\widetilde S_j(s,p_j,q_j) &= \begin{cases} S_j(p_j,q_j)\phantom{s,}\,\, &
\text{when \quad $s \geq 2 \rho_j(p_j,q_j)^2$},\\ 0 & \text{when
\quad $s \leq \rho_j(p_j,q_j)^2$}.
\end{cases}
\end{split}
\end{equation*}
For every derivative $X^\alpha$ of order $|\alpha|$ in the vector
fields $X_j$ and $X_{n+j}$ and every integer $N$ there is a
constant $C_{\alpha,N}$ so that
\begin{eqnarray*}
&& |X^{\alpha}\Phi_{j}(s,p_j,q_j)| \\&&\qquad \leq C_{\alpha, N}
\begin{cases}
s^{-\frac{1}{2}|\alpha|}\,|B_j(p_j,\sqrt s)|^{-1} & \text{if $s
\geq d_{j}(p_j,q_j)^2$},
\\\\
s^N\,d_{j}(p_j,q_j)^{-2-|\alpha|-2N}\,
\Lambda_{j}(p_j,d_j(p_j,q_j))^{-1} & \text{if $s\leq
d_{j}(p_j,q_j)^{2}$}.
\end{cases}
\end{eqnarray*}
Also{\rm ,}
$$
\left\vert X^{\alpha}\chi_j(s,p_j,q_j)\right\vert\leq
C_{\alpha}\,d_{j}(p_j,q_{j})^{-|\alpha|}.
$$
Moreover{\rm ,} as soon as $\alpha$ is different from zero{\rm ,}
$X^{\alpha}\chi_j(s,p_j,q_j)$ is supported where
$d_{j}(p_j,q_{j})^{2}\approx s$.
\end{proposition}

\Proof  The estimates follow from the chain rule and the basic estimates in Theorem \ref{T:1.1.1}.
\Endproof\vskip4pt  

The following result is an immediate consequence of the definitions:

\begin{proposition}\label{P:B1} Suppose that $d_j(p_j,q_j) \leq
d_\ell(p_\ell,q_\ell)$. For $M$ large enough depending only on the
type{\rm ,} 
$$
\Lambda_{j}\big(p_{j};d_j(p_j,q_j)\big)\,
\left(\frac{d_j(p_j,q_j)}{d_\ell(p_\ell,q_\ell)}\right)^{M}
\lesssim \Lambda_{\ell}(p_\ell,d_\ell(p_\ell,q_\ell)).
$$
\end{proposition}

\demo{Proof of Theorem {\rm \ref{T:MainEstimate}}}
The integrand in the integral defining $K(p,q)$ in equation
(\ref{E:2.1.1.a}) can be written
\begin{align*}
\prod_{j=1}^n H_{j}(s,p_{j},q_j) &- \prod_{j=1}^n S_{j}(p_{j},q_j)
\\ &=
\prod_{j=1}^n \left(\Phi_j(s,p_j,q_j) + \widetilde
S_j(s,p_j,q_j)\right)- \prod_{j=1}^n S_{j}(p_{j},q_j)\\ &=
\prod_{j=1}^n \Phi_j(s,p_j,q_j)
+\sum_A \left[\prod_{j\in A}\Phi_j(s,p_j,q_j)\,\prod_{j\notin A}
\widetilde S_j(s, p_j,q_j)\right]\\ &\qquad\qquad\qquad\qquad
+\left[\prod_{j=1}^n \widetilde S_j(s,p_j,q_j) - \prod_{j=1}^n
S_j(p_j,q_j)\right]
\end{align*}
where the sum is taken over all nonempty, proper subsets $A
\subset \{1,\ldots,n\}$.\break Note that

\begin{eqnarray*}
\noalign{\vskip-16pt}
&&\int_0^\infty \prod_{j\in A}\Phi_j(s,p_j,q_j)\,\prod_{j\notin A}
\widetilde S_j(s, p_j,q_j)\,ds \\&&\qquad = \Big(\prod_{j\notin
A}S_j(p_j,q_j)\Big)\int_0^\infty \prod_{j\in
A}\Phi_j(s,p_j,q_j)\,\prod_{j\notin
A}\chi\left(\frac{\rho_{j}(p_j,q_{j})^{2}}{s}\right)\,ds
\end{eqnarray*}
and that
\begin{equation*}
\begin{split}
\int_0^\infty \Big[\prod_{j=1}^n \widetilde S_j(s,p_j,q_j) &-
\prod_{j=1}^n S_j(p_j,q_j)\Big]\,ds \\&= \Big(\prod_{j=1}^n
S_j(p_j,q_j)\Big) \int_0^\infty \prod_{j=1}^n
\left[1-\chi\left(\frac{\rho_{j}(p_j,q_{j})^{2}}{
s}\right)\right]\,ds.
\end{split}
\end{equation*}
Set
\begin{equation}\label{E:6.2.1a}
\begin{split}
\widetilde K_0(p,q) &= \int_0^\infty \prod_{j=1}^m
\Phi_j(s,p_j,q_j)\,ds,\\ \widetilde K_A(p,q) &= \int_0^\infty
\prod_{j\in A}\Phi_j(s,p_j,q_j)\,\prod_{j\notin
A}\chi\left(\frac{\rho_{j}(p_j,q_{j})^{2}}{s}\right)\,ds,\\
\widetilde K_\infty(p,q) &= \int_0^\infty
\prod_{j=1}^n\left[1-\chi\left(\frac{\rho_{j}(p_j,q_{j})^{2}}{
s}\right)\right]\,ds.
\end{split}
\end{equation}
To establish Theorem \ref{T:MainEstimate}, we need to show that
the integrals defined in equation (\ref{E:6.2.1a}) satisfy the
estimates stated in equation (\ref{E:6.1.1a}).

\Subsubsec{Estimates for $\widetilde K_0$} 

\begin{lemma}  
\begin{equation*}
\Big\vert\int_0^\infty \prod_{j=1}^n
X^{\alpha_j}\Phi_j(s,p_j,q_j)\,ds \Big\vert\leq C_\alpha
\frac{d_\Sigma(p,q)^{2-|\alpha|}}{\prod_{j=1}^n \Big[d_\Sigma(p,q)^2
\Lambda_j(p_j,d_\Sigma(p,q))\Big]}.
\end{equation*}
\end{lemma}
\vskip8pt

\Proof  The case $n=2$ is entirely typical. We may suppose,
without loss of generality, that $$d_{1}(p_{1},q_{1}) \leq
d_{2}(p_{2},q_{2}),$$ so that $d_\Sigma(p,q) \approx d_{2}$. We
split the integral into three parts; first where $0\leq s\leq
d_{1}(p_1,q_1)^{2}$, next where
 $d_{1}(p_1,q_1)^{2}\leq
s \leq d_{2}(p_2,q_2)^{2}$, and finally where
$d_{2}(p_2,q_2)^{2}\leq s <+\infty$. Using the estimates from
Proposition \ref{P:6.2.2.1} and also Proposition \ref{P:B1}, we
have
\begin{align*}
\Big\vert \int_{0}^{d_{1}(p_1,q_1)^{2}}
&X^{\alpha}_{1}\Phi_{1}\pso
\,X^{\beta}_{2}\Phi_{2}\pst\,ds\Big\vert
\\&\lesssim
\frac{d_{1}(p_1,q_1)^{-2-|\alpha|-2M}\,\,d_{2}(p_2,q_2)^{-2-|\beta|-2N}}{
\Lambda_{1}(p_1,d_{1}(p_1,q_1))\,
\Lambda_{2}(p_2,d_{2}(p_2,q_2))}\, \int_{0}^{d_{1}(p_1,q_1)^{2}}
s^{M+N}\,ds \\ &\lesssim
\frac{d_{1}(p_1,q_1)^{-|\alpha|}\,\,d_{2}(p_2,q_2)^{-|\beta|}}{
\Lambda_{1}(p_1,d_{1}(p_1,q_1))\,\,d_{2}(p_2,q_2)^{2}\,\,\Lambda_{2}(p_2,
d_{2}(p_2,q_2))}\left(\frac{d_{1}(p_1,q_1)}{
d_{2}(p_2,q_2)}\right)^{2M}\\ &\lesssim
d_\Sigma(p,q)^{-2-|\alpha|-|\beta|}\,\Lambda_{1}(p_1,d_\Sigma(p,q))^{-1}\,
\Lambda_{2}(p_2,d_\Sigma(p,q))^{-1}.
\end{align*}
Similarly
\begin{align*}
\Big\vert\int_{{d_1(p_1,q_1)}^{2}}^{{d_2(p_2,q_2)}^{2}}\,
&X^{\alpha}_{1}\Phi_{1}\pso X^{\beta}_{2}\Phi_{2}\pst
\,ds\Big\vert
\\&\leq
\frac{{d_2(p_2,q_2)}^{-2N-|\beta|}}{
{d_2(p_2,q_2)}^{2}\,\Lambda_{2}(p_{2};{d_2(p_2,q_2)})}\,
\int_{{d_1(p_1,q_1)}^{2}}^{{d_2(p_2,q_2)}^{2}}s^{N-2-|\beta|}\,
\Lambda_{1}\big(p_{1};\sqrt s\big)^{-1}ds\notag\\ \notag\\ &\leq
\frac{{d_2(p_2,q_2)}^{-2-|\alpha|-2N-|\beta|}}{
\Lambda_{2}(p_{2};{d_2(p_2,q_2)})}\,
\Lambda_{1}(p_{1};{d_2(p_2,q_2)})^{-1}\,{d_2(p_2,q_2)}^{2N}\notag\\
\notag\\ &\lesssim
 d_\Sigma(p,q)^{-2-|\alpha|-|\beta|}\,\Lambda_{1}(p_{1};
 d_\Sigma(p,q))^{-1}
 \,\Lambda_{2}(p_{2}; d_\Sigma(p,q))^{-1} 
 \notag
\end{align*}
since for $N$ large enough, the function
$s^{N-2-|\beta|}\,\Lambda_{1}\big(p_{1};\sqrt s\big)^{-1}$ is
monotone increasing, and hence we can estimate the integral by
taking this function at the upper end point.

Finally
\begin{align*}
\Big\vert \int_{{d_2(p_2,q_2)}^{2}}^{+\infty}\,
&X^{\alpha}_{1}\Phi_{1}(s,p_{1},y_{1})
X^{\beta}_{2}\Phi_{2}(s,p_{2},y_{2})\,ds \Big\vert
\\&\leq
\int_{{d_2(p_2,q_2)}^{2}}^{+\infty} s^{-\frac{1}{
2}(|\alpha|+|\beta|)-2}\Lambda_{1}(p_{1};\sqrt s)^{-1}\,
\Lambda_{2}(p_{2};\sqrt s)^{-1}\,ds \notag\\\notag\\ &\lesssim
d_\Sigma(p,q)^{-2-|\alpha|-|\beta|}\,
\Lambda_{1}(p_{1};d_\Sigma(p,q))^{-1}\,
\Lambda_{2}(p_{2};d_\Sigma(p,q))^{-1}.
\end{align*}
\vglue-32pt \phantom{hi}
\hfq\vskip6pt

\Subsubsec{Estimates for $\widetilde K_A$}\hfill

\begin{lemma}\label{Lemma4.3.6} Let $A\subset\{1,\ldots, n\}$ be a
proper{\rm ,} nonempty
subset. Then
\begin{equation}\label{E:A1}
\begin{split}
\Big\vert \int_0^\infty \Big(\prod_{j\in A}
&X^{\alpha_j}\Phi_j(s,p_j,q_j)\Big)\Big(\prod_{j\notin A}
X^{\beta_j} \chi_j(s,p_j,q_j)\Big)\,ds\Big\vert\\ &\lesssim
\frac{d_\Sigma(p,q)^{2-|\alpha|-|\beta|}}{ \prod_{j\in A}d_\Sigma(p,q)^2
\Lambda_{j}\big(p_j;d_\Sigma(p,q)\big)}.
\end{split}
\end{equation}
\end{lemma}
\vskip8pt

\Proof 
Let $\delta_j = \rho_j(p_j,q_j)$. Then $d_\Sigma(p,q) \approx
\max_{j}\{\delta_j\}$. Also set $\delta_A = \max_{j\notin
A}\{\delta_j\}$. Then the support of the integrand in equation
(\ref{E:A1}) is $s \geq \delta_A^2$.

\demo{Remark} If any of the integers $\{\beta_j\}$
is nonzero, the support of the integrand is between two multiples
of the corresponding $\delta_j^2$. Thus if some $\beta_j > 0$, all
the $j \notin A$ for which $\beta_j >0$ have the property that
$\delta_j \approx \delta_A$.

\demo{Case {\rm 1:} $\beta_j = 0$ for every $j\notin A$} 
There are now two subcases to consider, depending on the relative
size of $\delta_A$ and $d_\Sigma(p,q)$. We always have $\delta_A
\lesssim d_\Sigma(p,q)$.
\demo{Case {\rm  1a:} $\delta_A \approx d_\Sigma(p,q)$} 
This is the case in which an index $j$ for which $\delta_j$ is comparable to the
maximum does not belong to the set $A$.
\Enddemo

Now we must study
\begin{equation}\label{E:A2}
\Big\vert \int_{d_\Sigma(p,q)^2}^\infty \Big(\prod_{j\in A}
X^{\alpha_j}\Phi_j(s,p_j,q_j)\Big)\Big(\prod_{j\notin A}
 \chi_j(s,p_j,q_j)\Big)\,ds\Big\vert.
\end{equation}
The decay of the integrand at infinity allows us to estimate the
integral by $d_\Sigma(p,q)^2$ times the value of the integrand at
the lower endpoint. Since\break $|\chi_j(s,p_j,q_j)|\leq 1$, we obtain
from Proposition \ref{P:6.2.2.1} the estimate
$$
d_\Sigma(p,q)^2  \Big(\prod_{j\in A}
X^{\alpha_j}\Phi_j(d_\Sigma(p,q)^2,p_j,q_j)\Big)  \leq
 \frac{d_\Sigma(p,q)^{2-\sum_{j\in
A}|\alpha_j|}}{\prod_{j\in A} d_\Sigma(p,q)^2\,\Lambda_j\big( p_j,
d_\Sigma(p,q)\big)}
$$
since for every $j \in A$, $ d_j(p_j,q_j)\leq d_\Sigma(p,q)$. This
gives the desired estimate in this case.

\demo{Case {\rm 1b:} $\delta_A \ll  d_\Sigma(p,q)$}  This is
the case in which every index $j$ for which $\delta_j$ is comparable to the maximum
does belong to the set $A$.
\Enddemo

This time we must study
\begin{equation}\label{E:A3}
\Big\vert \int_{\delta_A^2}^\infty \Big(\prod_{j\in A}
X^{\alpha_j}\Phi_j(s,p_j,q_j)\Big)\Big(\prod_{j\notin A}
 \chi_j(s,p_j,q_j)\Big)\,ds\Big\vert.
\end{equation}
We write
\begin{equation*}
\int_{\delta_A^2}^\infty = \int_{\delta_A^2}^{d_\Sigma(p,q)^2} +
\int_{d_\Sigma(p,q)^2}^\infty.
\end{equation*}
The second integral is handled in the same way as Case 1a. In
dealing with the first integral, we would like to be able to take
out the integrand at the {\it upper} endpoint rather than at the
lower endpoint. However, note that there is an index $j \in A$ for
which $\delta_j \approx d_\Sigma(p,q)$. For this $j$ we use the
estimate from Proposition \ref{P:6.2.2.1}:
\begin{equation*}
\big\vert X^\alpha_j \Phi_j(s,p_j,q_j)\big\vert \leq C_{\alpha, N}
s^N\,d_\Sigma(p,q)^{-2-|\alpha_j|- 2N} \Lambda_j\big(p_j,
d_\Sigma(p,q)\big)^{-1}.
\end{equation*}
If we take $N$ large enough, we can make the entire integrand in
equation (\ref{E:A3}) monotone increasing, and we get the desired
estimate
$$
\frac{d_\Sigma(p,q)^{2-\sum_{j\in A}|\alpha_j|}}{\prod_{j\in A}
d_\Sigma(p,q)^2\,\Lambda_j\big( p_j, d_\Sigma(p,q)\big)}.
$$

We now must consider what happens if some $\beta_j > 0$.

\vskip6pt {\it Case {\rm 2:} Some $\beta_j > 0$.}  By the remark following Lemma \ref{Lemma4.3.6}, we
are led to study
\begin{equation}
\prod_{j\notin A}
\delta_A^{-|\beta_j|}\int_{\frac{1}{2}\delta_A^2}^{2\delta_A^2}
\big\vert \prod_{j\in A} X^{\alpha_j}\Phi_j(s,
p_j,q_j)\big\vert\,ds.
\end{equation}

\vskip2pt 
{\it Case {\rm 2a:} $\delta_A \approx d_\Sigma(p,q)$.}  As in Case 1a, we can estimate the integral by
$d_\Sigma(p,q)^2$ times the value of the integrand at the lower endpoint. This gives
us the estimate
$$
\frac{d_\Sigma(p,q)^{2-\sum_{j\in A}|\alpha_j|-\sum_{j\notin
A}|\beta_j|}}{\prod_{j\in A} d_\Sigma(p,q)^2\,\Lambda_j\big( p_j,
d_\Sigma(p,q)\big)}
$$
which is the desired result.

\vskip5pt {\it Case {\rm 2b:} $\delta_A \ll  d_\Sigma(p,q)$.} Again, as in Case 1b, we know there exists $j
\in A$ for which we have the estimate
$$
\big\vert X^\alpha_j \Phi_j(s,p_j,q_j)\big\vert \leq C_{\alpha, N}
s^N\,d_\Sigma(p,q)^{-2-|\alpha_j|- 2N} \Lambda_j\big(p_j,
d_\Sigma(p,q)\big)^{-1}.
$$
When we integrate $s^n\,d_\Sigma(p,q)^{-2N}$ between two multiples
of $\delta_A^2$, we get the factor
\begin{equation*}
\delta_A^2\,\left(\frac{\delta_A^2}{d_\Sigma(p,q)^2}\right)^N \leq
d_\Sigma(p,q)^2
\,\left(\frac{\delta_A^2}{d_\Sigma(p,q)^2}\right)^N.
\end{equation*}
If we take $N$ large enough, by using Proposition \ref{P:B1}, we
can replace each factor $\delta_A^{-|\beta_{j}|}$ with
$d_\Sigma(p,q)$, and again we get the required estimate. This
completes the proof. \hfq

\vglue-4pt
\Subsubsec{Estimates for $\widetilde K_\infty$} 
 
\vglue-20pt
\phantom{up}
\begin{lemma} 
$$
\Big\vert\int_0^\infty \prod_{j=1}^n X^{\alpha_j}\left(1 -
\chi_j(s,p_j,q_j)\right]\,ds \Big\vert\lesssim
d_\Sigma(p,q)^{2-\sum|\alpha_j|}.
$$
\end{lemma}

\Proof  The integrand is zero unless $s \leq
d_\Sigma(p,q)^2$, and each factor $\delta_j^{-|\alpha_j|}$ coming
from a derivative is dominated by $d_\Sigma(p,q)^{-|\alpha_j|}$.
This gives the desired estimate, and completes the proof.
\hfill\qed\vskip4pt

Thus the proof of Theorem \ref{T:MainEstimate} is
complete.\hfq 

\Subsec{Quadratic derivatives of $\widetilde K$} If $Q(Z, \bar Z)$ is a quadratic expression in the vector
fields $\{Z_1, \ldots, Z_n, \bar Z_1, \ldots, \bar Z_n\}$, it need not be the case that the operator $Q(Z, \bar
Z)\,\widetilde \K$ is bounded on $L^2(\widetilde M)$. In this section we obtain certain replacements for this
loss of maximal hypoellipticity. 

To describe the results, suppose that we are studying the relative fundamental solution operator $\widetilde K$ for the differential operator $\bx_b = \sum_{j=1}^n \bx_j$ where $\bx_j = W_j\,\bW_j$, and each $W_j$ is one of $\{Z_j, \overline Z_j\}$ so that $\bW_j$ is the other. Let $b$ be a bounded function on $\widetilde M$ (where we will write $b(p) = b(p_1, \ldots, p_n)$). We obtain conditions on $b$ that guarantee that the operators $b\,W_k\,W_l\,\widetilde \K$, $b\,W_k\,\bW_l\,\widetilde \K$, $b\,\bW_k\,W_l\,\widetilde \K$, and $b\,\bW_k\,\bW_l\,\widetilde \K$ are bounded on $L^p(\widetilde M)$ for $1 < p < \infty$. 
The size conditions that we need to impose on $b$ will depend on which of these forms we consider, and are given in terms of the quantities $\big\{\lambda_j(p_j) = \frac{\partial^2 P_j}{\partial z_j\partial \bar z_j}(p_j)\big\}$
 which\break\vskip-12pt\noindent are the eigenvalues of the Levi form on the decoupled boundary $M$.

\begin{theorem}\label{T:QuadraticDerivatives} Let $1 \leq k,l\leq n$. Then
 
\begin{enumerate}
\item The operators $W_k\,\bW_k\,\widetilde \K${\rm ,}
 $\bW_k\,\bW_k\,\widetilde \K${\rm ,} and $\bW_k\,\bW_l\,\widetilde \K$ are product singular integrals on
$\widetilde M$ and hence are bounded on $L^p(\widetilde M)$ for $1 < p < \infty$.

\item Let $b$ be a bounded function on $\widetilde M$ and suppose there exists a constant $C$ so that 
$$
\lambda_k(p_k)\,|b(p_1, \ldots, p_n)| \leq C \,\inf_{l\neq k}\lambda_l(p_l).
$$
Then  
$$
b\,\bW_k\,W_k\,\widetilde \K = \sum_\alpha b_\alpha\,\mathcal M_\alpha \,+ \,\sum_{\substack{l,A}} \big(T_k-T_{l}\big) b_{l,A}\,\lambda_k
\, \widetilde\N_A \, \otimes \widetilde \S_A.
$$
In the first sum on the right-hand side{\rm ,}
 each $b_\alpha$ is a bounded function on $\widetilde M${\rm ,} each $\mathcal M_\alpha$ is a product singular
integral on $\widetilde M${\rm ,} $T_j = \frac{\partial}{\partial t_j}${\rm ,} and the sum is finite. In the second
summation{\rm ,} each $b_{l,A}$ is a bounded function{\rm ,}
 and the sum is over all $1 \leq l\leq n${\rm ,} and all subsets $A
\subset \{1,\ldots, n\}$ with $l \in A$ and $k \notin A$.

\item Let $b$ be a bounded function on $\widetilde M$ which is independent of the variable $p_k${\rm ,}
 and suppose there exists a constant $C$ so that 
$$
|b(p_1, \ldots, p_n)| \leq C \,\inf_{l\neq k}\lambda_l(p_l).
$$
Then there are {\rm NIS}
 operators $\{P_\alpha, P_{l,A}\}$ of order zero acting only in the variable $p_k$ so that  
$$
b\,W_k\,W_k\,\widetilde \K = \sum_\alpha P_\alpha b_\alpha\,\mathcal M_\alpha \,+ \,\sum_{\substack{l,A}} P_{l,A}\,\big(T_k-T_{l}\big) b_{l,A}\,\lambda_k
\, \widetilde\N_A \, \otimes \widetilde \S_A.
$$
Again{\rm ,}
 each $b_\alpha$ and $b_{l,A}$ is a bounded function on $\widetilde M${\rm ,}
 each $\mathcal M_\alpha$ is a
product singular integral on $\widetilde M${\rm ,}
 $T_j = \frac{\partial}{\partial t_j}${\rm ,} and the second sum is over
all $1 \leq l\leq n${\rm ,} and all subsets $A \subset \{1,\ldots, n\}$ with $l \in A$ and $k \notin A$.
\end{enumerate}
\end{theorem}

  \emph{Proof of} (1). The key point is that $W_k\,\bW_k\,\S_k = 0$, $\bW_k\,\bW_k \,\S_k = 0$,
and $\bW_k\,\bW_l\,\S_k\otimes\S_l = 0$. The proof in all three cases is similar, so we only give the details
for the first. According to equation  (\ref{E:First0}) we have $\widetilde K = \widetilde \N +
\sum_A\widetilde \N_A\otimes\widetilde\S_A$, and hence
$$
W_k\,\bW_k\,\widetilde \K = W_k\,\bW_k\,\widetilde \N + \sum_A W_k\,\bW_k\,(\widetilde \N_A \otimes \widetilde \S_A).
$$
But $W_k\,\bW_k\,\widetilde \N$ is a product singular integral operator by Theorem \ref{T:5.1.1.1}. If $k \notin A$ then
$$
W_k\,\bW_k\,(\widetilde \N_A \otimes \widetilde \S_A) = \widetilde \N_A \otimes W_k\,\bW_k\,\widetilde \S_A = 0.
$$
On the other hand, if $k \in A$, then 
$$
W_k\,\bW_k\,(\widetilde \N_A \otimes \widetilde \S_A) = W_k\,\bW_k\,\widetilde \N_A \otimes \widetilde \S_A.
$$
But $W_k\,\bW_k\,\widetilde \N_A$ is a product singular integral operator in the variables coming from the set $A$, as follows from Theorem \ref{T:5.1.1.1} applied to fewer variables. Thus $W_k\,\bW_k\,(\widetilde \N_A \otimes \widetilde \S_A)$ 
(and hence $W_k\,\bW_k\,\widetilde \K$) is a product singular integral operator. This completes the proof of
(1).

\demo{Proof of {\rm (2)}} 
 We now consider the operator $b\,\bW_k\,W_k\,\widetilde \K =b\, \overline\bx_k\,\widetilde \K$. Arguing
as in the proof of (1), we have
\begin{eqnarray*}
b\,\overline \bx_k\,\widetilde \K &=& b\,\overline \bx_k\, 
\widetilde\N \,\,+\,\, \sum_{A\subset\{1,\ldots,n\}} b\,\overline\bx_k\big(\widetilde \N_A
\otimes \widetilde\S_A\big)\\ &=& b\,\overline \bx_k \, \widetilde\N \,\,
+\sum_{\substack {k\in A\subset\{1,\ldots, n\}}} \!\!\!\! b\,\big(\overline
\bx_k\,\widetilde \N_A\big)\otimes \widetilde \S_A +
\sum_{\substack {k\notin A\subset\{1,\ldots, n\}}}\!\!\!\! b\,\widetilde\N_A \otimes
\big(\overline \bx_k \widetilde \S_A\big).
\end{eqnarray*}

Now $\overline{\bx}_k\,\N$ is a product
NIS operator by Theorem \ref{T:5.1.1.1}. Also
$\overline\bx_k\widetilde\N_A$ is a product singular integral operator in the
variables $k \in A$, and so $(\overline\bx_k\widetilde\N_A)\otimes
\widetilde\S_A$ is a product singular integral operator on $\widetilde M$. Thus for any bounded function $b$ on $\widetilde M$, the first two terms on the right in the last equation are bounded operators on $L^p(\widetilde M)$ for $1 < p < \infty$. 

The difficult terms are those involving $\widetilde\N_A \otimes
\big(\overline \bx_k \widetilde \S_A\big)$ where $k \notin A$,
since in terms of the product structure, $\big(\overline \bx_k
\widetilde \S_A\big)$ is smoothing of order $-2$, while
$\widetilde\N_A$ is smoothing of order $+2$. In general, such a
product does not yield a product NIS operator.

In these bad terms, we transfer the two extra
derivatives from the right-hand side of the tensor product to the
left-hand side. Note that $\overline \bx_k = \bx_k \pm
\lambda_k\,T_k$ where $T_k = \frac{\partial}{\partial t_k}$. Since $k \notin A$, we have
$\bx_k\,\widetilde \S_A = 0$. Choose an index $l \in A$, so that in particular $l \neq k$. Then
\begin{eqnarray*}
b\,\widetilde\N_A \otimes \big(\overline \bx_k \widetilde \S_A\big)
&=& \pm b\, \lambda_k\,T_k\, \widetilde\N_A \,
\otimes \widetilde \S_A\\ &=&
\mp \,\Big(\frac{b\,\lambda_k}{\lambda_{l}}\Big)
\Big[\lambda_{l}\,T_{l}\, \widetilde\N_A \,
 \otimes \widetilde \S_A\Big] \pm b\,\lambda_k\big(T_k-T_{l}\big)
\, \widetilde\N_A \, \otimes \widetilde \S_A.
\end{eqnarray*}
Now $\lambda_{l}\,T_{l}\,\widetilde\N_A \,
\otimes \widetilde \S_A$
is a product singular integral operator on $\widetilde M$ since $\lambda_l\,T_l\,\widetilde \N_A$ is the commutator of two good derivatives applied to $\widetilde \N_A$.
 Thus $$
\Big(\frac{b\,\lambda_k}{\lambda_{l}}\Big)
\Big[\lambda_{l}\,T_{l}\, \widetilde\N_A \,
 \otimes \widetilde \S_A\Big]
$$
is a bounded operator on $L^p(\widetilde M)$ provided that there is a constant $C$ so that \begin{equation*}
\lambda_k(p_k)\,|b(p_1, \ldots, p_n)| \leq C \lambda_l(p_l).
\end{equation*}
Finally, note that the operator $b\,\lambda_{k}(T_{k}-T_{l}) - (T_{k}-T_{l})b\,\lambda_{k} $ is a bounded function.
This completes the proof of (2).

\demo{Proof of {\rm (3)}} 
 Recall that $\bar Z_j = X_j + i X_{n+j}$ where $\{X_j, X_{n+j}\}$ are real vector fields on the manifold
$M_j$. It follows that $\bx_k + \overline\bx_k = Z_j\,\bar Z_j + \bar Z_j\,Z_j = 2(X_k^2+ X_{n+k}^2)$, and
it is known that this operator can be inverted with an NIS operator on $M_k$, smoothing of order $2$. It then
follows that there are NIS operators $P_1$ and $P_2$ on $M_k$, smoothing of order zero, such that
$$
W_k\,W_k = P_1\bx_k + P_2\overline\bx_k.
$$
We can also regard $P_1$ and $P_2$ as operators on $\widetilde M$ which act only in the variable~$p_k$.
Thus if $B$ is a bounded function on $\widetilde M$ which is independent of the variable~$p_k$, the
operator which is multiplication by $B$ and the operator $P_j$ commute. It follows that we have
$$
B\,W_k\,W_k\,\widetilde \K = P_1\,B\,(\bx_k\,\widetilde \K) + P_2\,B\,(\overline \bx_k\,\widetilde \K).
$$
The proof of part (3) now follows immediately from parts (1) and (2).

\section{Transference from $M_1\times \cdots \times M_n$ to $M$ and $L^p$ regularity of $\K$}\label{S:Descending}

In order to pass from results about operators on the product $\widetilde M = M_1\times \cdots \times M_n$ to results about
operators on the decoupled boundary $M$, we use the mapping $\pi:\widetilde M \to M$ given by $\pi(z_1,\ldots,z_n,t_1, \ldots, t_n) = (z_1,\ldots,z_n,t_1+\cdots + t_n)$. We have already observed that $d\pi$ maps the $\dbar_b$-complex on $\widetilde M$ to the $\dbar_b$-complex on $M$. As discussed in Section \ref{SS:Outline}, we can also use $\pi$ to transfer the relative fundamental solutions for $\bx_b$ on $\widetilde M$ to relative fundamental solutions for $\bx_b$ on $M$. In Section \ref{SS:Transference} below, we use standard transference techniques to show that the resulting operators have the same $L^p$ norm on $M$ as the original operators have on $\widetilde M$. In Section \ref{SS:PFT} we study what this transference does to products of projections on $\widetilde M$. In Section \ref{SSDistributions} we show the existence of  relative fundamental solutions $\K$ and $\N$, and in Section \ref{SSLPregularity} we obtain $L^p$-regularity results for the relative
fundamental solutions $\N$ and $\K$ for $\bx_b$ on the decoupled boundary $M$.

\Subsec{A general transference result}
\label{SS:Transference} Let $\widetilde T$ be a measurable function on $\C^n\times \C^n\times \R^n$ with
compact support. Suppose that
\begin{equation*}\label{E:7.1.1}
\begin{split}
\sup_{z\in \C^n}\iint_{\C^n\times\R^n}|\widetilde T(z,w,t)|\,dw\,dt &= C_1 < +\infty;\\ 
\sup_{w\in \C^n}\iint_{\C^n\times\R^n}|\widetilde T(z,w,t)|\,dz\,dt &= C_2 < +\infty.
\end{split}
\end{equation*}
Define an operator $\widetilde \T$ acting on functions on $\C^n\times\R^n$ by setting
$$
\widetilde \T[F](z,t) = \iint_{\C^n\times \R^n} \widetilde T(z,w,t-s)\,F(w,s)\,dw\,ds.
$$
$\widetilde \T$ is then a bounded operator on $L^p(\C^n\times\R^n)$ with the bound at most $C_1^{\frac{1}{p}} \, C_2^{\frac{1}{p'}}$ where $\frac{1}{p}+\frac{1}{p'} = 1$. 
\vskip2pt

Next, if $(z,w,t) \in \C^n\times\C^n \times \R$, set
\begin{eqnarray*}\label{E:7.1.6}
T(z,w,t) &= &\int_{\Sigma(t)}\widetilde T(z,w,s) \,d\widetilde s\\
&=& \int \widetilde T(z,w,t-s_2-\cdots-s_n, s_2,\ldots,s_n)\,ds_2\,\cdots\,ds_n, 
\end{eqnarray*}
where as before $\Sigma(t)$ is the affine hyperplane $\{s\in\R^n\,\big\vert\,\sum_{j=1}^n s_j = t\}$, and $d\tilde s$ is $(n-1)$ dimensional measure on $\Sigma(t)$. Given a measurable function $f$ on $\C^n \times \R$, for $(z,t) \in \C^n\times \R$ 
define
\begin{eqnarray*} 
\T[f](z,t) &\equiv& \iint_{\C^n\times \R} T(z,w,t-s) \,f(w,s)\,dw\,ds\\
&=& \iint_{\C^n\times \R^n} \widetilde T(z,w,s)\,f(w,t-s_1-\cdots - s_n) \,dw\,ds_1\,\cdots\,ds_n\\
&=& \iint_{\C^n\times \R^n} \widetilde T(z,w,s)\,R_s[f](w,t)\,dw\,ds_1\,\cdots\,ds_{n}
\end{eqnarray*}
where for $s = (s_1, \ldots, s_n)$, $R_s[f](z,t) = f(z,t-s_1 - \cdots - s_n)$. Note that $\T[f]\circ\pi = \widetilde \T[f\circ\pi]$.
It follows that if $f \in L^\infty(\C^n\times \R)$, these integrals converge absolutely for all $(z,t) \in \C^n\times \R$.
Moreover, 
\begin{equation*}\label{E:7.1.8}
\begin{split}
\sup_{z\in \C^n}\iint_{\C^n\times\R}|T(z,w,t)|\,dw\,dt &= C_1 <
+\infty;\\ \sup_{w\in \C^n}\iint_{\C^n\times\R}|T(z,w,t)|\,dz\,dt
&= C_2 < +\infty.
\end{split}
\end{equation*}
Thus $\T$ is bounded on $L^p(\C^n\times \R)$ with norm at most $C_1^{\frac{1}{p}}\,C_2^{\frac{1}{p'}}$. 

We now have the following transference result which shows that a better {\it a priori} $L^p$ bound for $\widetilde \T$ gives the same bound for $\T$.

\begin{theorem}\label{T:7.1.1.5} Suppose there is a constant $A_p$ such that
$$
\norm \widetilde \T[F]\norm_{L^p(\C^n\times \R^n)} \leq A_p\,\norm
F\norm_{L^p(\C^n\times\R^n)}.
$$
Then the operator $\T$ satisfies
$$
\norm \T[f]\norm_{L^p(\C^n\times\R)} \leq A_p \, \norm
f\norm_{L^p(\C^n\times\R)}.
$$
\end{theorem}

\Proof  We follow the argument in \cite[Ch.~2]{CoifmanWeiss}. Let $E\subset \R^n$ denote the (compact)
projection onto $\R^n$ of the
 compact support of the function $\widetilde T$. Thus $\widetilde T(z,w,s) \neq 0$ implies $s \in E$. Choose
$\varepsilon > 0$, and choose a (large) bounded open set $V \subset \R^n$ so that if
$$
V + E = \left\{x\in \R^n\,\Big\vert\,\text{ $x = v +e$ with $v \in
V$ and $e \in E$}\right\},
$$
then
$$
\frac{\left\vert V+E\right\vert}{\left\vert V\right\vert} \leq 1 +
\varepsilon.
$$

We have
$$
R_y\big[\T[f]\big](z,t) = \underset{\C^n\times\R^n}{\iint} \widetilde T(z,w,s-y)\,R_s[f](w,t)\,dw\,ds.
$$
Let $\chi$ be the characteristic function of $V+E$. Since for any $y \in \R^n$,  
$$
\norm R_y[f]\norm_{L^p(\C^n\times \R)} = \norm f\norm_{L^p(\C^n\times\R)},
$$
we can average over $V$ and obtain
\begin{small}
\begin{eqnarray*}
 ||\T[f]||_{L^p(\C^n\times R)}^p  
&  =& \frac{1}{|V|} \int_V  ||R_y\big[\T[f]\big]||^p_{L^P(\C^n\times
\R)}\,dy
\\ &  =& \frac{1}{|V|} \int_V\Big[\underset{\C^n\times\R}{\iint} \big\vert
R_y\big[\T[f](z,t)\big]\big\vert^p\,dz\,dt\Big]\,dy
\\  &  =& \frac{1}{|V|} \int_V\Big[\underset{\C^n\times\R}{\iint} \big\vert
\T[f](z,t-y_1-\cdots-y_n)\big\vert^p\,dz\,dt\Big]\,dy
\\ &  =& \frac{1}{|V|} \int_V\Big[\underset{\C^n\times\R}{\iint} \Big\vert
\underset{\C^n\times\R^n}{\iint} \widetilde T(z,w,s-y)\,R_s[f](w,t)\,dw\,ds \Big\vert^p\,dz\,dt\Big]\,dy
\\ & \leq& \frac{1}{|V|}\int_{\R}\Big[\underset{\C^n\times\R^n}{\iint}\Big\vert
\underset{\C^n\times\R^n}{\iint} \widetilde T(z,w,s-y)\,R_s[f](w,t)\,\chi(s)\,dw\,ds\Big\vert^p dz\,dy\Big] dt
\\ &  =& \frac{1}{|V|}\int_{\R}\Big[\underset{\C^n\times\R^n}{\iint}\Big\vert
\underset{\C^n\times\R^n}{\iint} \widetilde T(z,w,y-s)\,R_{-s}[f](w,t)\,\chi(-s)\,dw\,ds\Big\vert^p
dz\,dy\Big] dt
\\ &  =& \frac{1}{|V|}\int_{\R}\norm \widetilde \T[\widetilde F_t]\norm^p_{L^p(\C^n\times\R^n)}\,dt
\\ &  \leq& A_p^p\, \frac{1}{|V|}\, \int_{\R} \norm \widetilde F_t \norm^p_{L^p(\C^n\times \R^n)}\,dt
\end{eqnarray*}
\end{small}

\noindent
where
$$
\widetilde F_t(w,s) = R_{-s}[f](w,t)\,\chi(-s)= f(w,t+s_1+\cdots +
s_n)\,\chi(-s).
$$
But then
\begin{equation*}
\begin{split}
\int_{\R} \norm \widetilde F_t\norm^p_{L^p}\,dt &= \underset{\C^n\times \R^{n+1}}{\iiint} |f(w,t+s_1+\cdots +s_n)\,\chi(-s)|^p \,dw\,ds\,dt
\\& = \big\vert V+E\big\vert\,\,\norm f \norm^p_{L^p(\C^n\times\R)}.
\end{split}
\end{equation*}
It follows that
$$
\norm \T[f] \norm^p_{L^p(\C^n\times \R)} \leq A_p^p\,\,
(1+\varepsilon)\,\norm f \norm^p_{L^p(\C^n\times \R)},
$$
which completes the proof. \hfq

\Subsec{Transference of products of projections}\label{SS:PFT} Recall that on $\widetilde M$, the 
relative\break\vskip-11pt\noindent 
fundamental solutions $\widetilde \N_J$ satisfy 
$$
\widetilde \N_J\,\bx_J = \bx_J\,\widetilde \N_J  = \prod_{j=1}^n \big(I-\S_J^{J(j)}\big) = I + \sum_{\emptyset\neq A\subset \{1,\ldots, n\}} (-1)^{|A|}\, \widetilde\S_{J,A}
$$
where $\widetilde \S_{J,A}$ is the operator with distribution kernel
$$
\widetilde S_{J,A}(p,q) = \bigotimes_{j \in A} S_j^{J(j)}(p_j,q_j).
$$
We want to understand what happens when this operator is transferred to $M$ by the mapping $\pi$. If $|A| = r$ we obtain an operator $\S_{J,A}$ on $M$ whose distribution kernel is given by
\begin{equation}\label{ESzegoProjections}
S_{J,A}(z_1, \ldots, z_n, w_1, \ldots, w_n, t) = \int_{\Sigma_r(t)} \prod_{j\in A} S_j^{J(j)}(z_j,w_j,r_j)\,d\tilde r,
\end{equation}
where $\Sigma_r(t) = \big\{(r_1,\ldots,r_q)\,\big\vert\,r_1+\cdots + r_q = t\big\}$, and $d\tilde r$ denotes the $(r-1)$-dimensional Lebesgue measure on $\Sigma_r(t)$. In other words, $S_{J,A}$ is the convolution in the $t$-variable of the $r$ functions $\{S_j^{J(j)}\}_{j\in A}$. We can study such distributions by taking the partial Fourier transform in the $t$-variable, defined by
$$
\hat f(z,w,\tau) = \mathcal F[f](z,w,\tau) = \int_{-\infty}^{+\infty} e^{-2\pi i t\tau} f(z,w,t)\,dt.
$$

We have the following basic fact about the partial Fourier transforms of the individual distributions $S_j^{\pm}$.

\begin{lemma} Let $\varphi:\C\to \R$ be  subharmonic{\rm ,}
 and let $\displaystyle Z = \frac{\partial}{\partial z} + i \frac{\partial \varphi}{\partial
z}\,\frac{\partial}{\partial t}$ and $\displaystyle \bar Z = \frac{\partial}{\partial \bar z} - i \frac{\partial
\varphi}{\partial \bar z}\,\frac{\partial}{\partial t}$. Let $f \in L^2(\C\times \R)$. If $Z[f] = 0${\rm ,} the support
of the partial Fourier transform $\hat f$ is contained in $\left\{(z,\tau)\in \C\times\R\,\Big\vert\,\tau \leq
0\right\}$. If $\bar Z[f] = 0${\rm ,}
 the support of the partial Fourier transform $\hat f$ is contained in
$\left\{(z,\tau)\in \C\times\R\,\Big\vert\,\tau \geq 0\right\}$. 
If $\S^{(+)}$ is the orthogonal projection onto
the null space of $Z${\rm ,}
 then the support of the partial Fourier transform of the distribution kernel $S^{(+)}$ is
supported where $\tau \leq 0$.  If $\S^{(-)}$ is the orthogonal projection onto the null space
 of $\bar Z${\rm ,} then
the support of the partial Fourier transform of the distribution 
  kernel $S^{(-)}$ is supported where $\tau \geq 0$. 
\end{lemma}

\Proof  Since the proofs are identical, we only deal with the case of the operator $Z$. We have
$$
\widehat Z[g](z,\tau) = e^{\tau \varphi(z)}\,\frac{\partial }{\partial z}\big[e^{-\tau \varphi}\,g\big](z,\tau).
$$
Since $\mathcal F$ is an isometry, it follows that if $Z[f] = 0$, then for almost every $\tau$ the function $z \to e^{-\tau \varphi(z)}\mathcal F[f](z,\tau)$ is an anti-holomorphic function $h_\tau$. Again since $\mathcal F$ is an isometry, it follows that  
for almost every such $\tau$,  
$$
 \int_{\C} |h_\tau(z)||^2 e^{\tau \varphi(z)}\,dm(z) <\infty.
$$
However, if $\tau > 0$, this implies that $h_\tau \equiv 0$. In fact, if $\tau > 0$,  then
$$
\frac{\partial ^2}{\partial z\,\partial \bar z}\left[|h_\tau|^2\,e^{\tau \varphi}\right] = e^{\tau \varphi}\left[ \left\vert \frac{\partial f}{\partial z} + \tau \,f\,\frac{\partial \varphi}{\partial z}\right\vert^2 + \tau\,|h_\tau|^2 \frac{\partial ^2 \varphi}{\partial z\,\partial \bar z}\right] \geq 0.
$$
Hence the function $z \to |h_\tau(z)|^2 e^{\tau \varphi(z)}$ is subharmonic, and so its value at any point is dominated by its average over a disk centered at the point of radius $r$. Letting $r \to \infty$ shows that $|h_\tau(z)|^2\,e^{\tau \varphi(z)} = 0$. Thus we have shown that if $f$ is in the null space of the operator $Z$ in $L^2(\C\times \R)$ , then the support of $\mathcal F[f]$ is contained in $\left\{(z,\tau)\,\big\vert\,\tau \leq 0\right\}$.

Next, since $\S[f](z,t) = \iint S(z,w,t-s)\,f(w,s)\,dw\,ds$, it follows that $\hat\S[f](z,\tau) = \int \hat S(z,w,\tau)\,\hat f(w,\tau)\,dw$. Since $\hat S[f](z,\tau) \equiv 0$ for all $f \in L^2(\C\times \R)$ and all $\tau > 0$, it follows that $\hat S(z,w,\tau) \equiv 0$ for $\tau > 0$. This completes the proof. \Endproof\vskip4pt

If we apply this result to the operators $\{\S_j^{J(j)}\}$ and use equation (\ref{ESzegoProjections}), we obtain

\begin{lemma}\label{LemmaS_JA} The operators $\S_{J,A} = 0$ on $M\times M$ unless
 $A \subset J$ or $A \cap J =~\emptyset$.
\end{lemma}

 \vglue-12pt
\Subsec{Relative fundamental solutions for $\bx_J$ on $M$}\label{SSDistributions} Recall that for each $J
\in \I_q$, we have an operator
$$
\bx_J = \sum_{j=1}^n \bx_j^{J(j)}
$$
which acts either on $\widetilde M$ or on $M$. We have constructed relative fundamental solutions for the operator on $\widetilde M$.
We then construct relative fundamental solutions for $\bx_J$ on $M$  by transferring operators $\widetilde \N_J$ and $\widetilde
\K_J$ on $\widetilde M$ to operators $\N_J$ and $\K_J$ on $M$.  This gives the first of the main results of this paper:

\demo{\scshape Theorem 2.4.1} \textit{For each of the $2^n$ possible operators $\{\bx_J\}${\rm ,} there is a
distribution
$K_J$ on $M\times M$ so that if $\mathcal{K}_J$ denotes the linear operator 
$$
\K_J[\varphi](p)= \int_M \varphi(q)\,K_J(p,q)\,dq,
$$
then
$$
\K_J\,\bx_J = \bx_J\,\K_J =
\begin{cases}
I-S_0 &\text{if \, $\bx_J$ acts on functions};
\\ I & \text{if \, $\bx_J$ acts on a $(0,r)$-form with $1\leq r \leq n-1$};
\\ I-S_{n} &\text{if \, $\bx_J$ acts on $(0,n)$-forms}.
\end{cases}
$$}

\Proof  Formally the kernel of the operator $\K_J$ is given by 
$$
K_J(z,w,t)= \int_{\Sigma(t)} \widetilde K_J(z,w,r)\,d\tilde r.
$$
However, we have observed that $\widetilde K_J$ can be approximated by functions which are
 smooth and have compact support in $(z,w,t)$. For such approximations, the integral converges absolutely.
Since on $\widetilde M$ 
$$
\widetilde \K_J\,\bx_J = \bx_J\,\widetilde \K_J = I - \prod_{j=1}^n S_j^{J(j)},
$$
it follows from Lemma \ref{LemmaS_JA} that on $M$ we have the desired equation for $\K_J$. This completes the proof.
\Endproof\vskip4pt  

By a similar argument we have a companion result for the transfer of the operator $\widetilde \N$.

\demo{\scshape Theorem 2.4.1a} \textit{For each of the $2^n$ possible operators $\{\bx_J\}${\rm ,}
 there is a
distribution
$N_J$ on $M\times M$ so that if $\mathcal{N}_J$ denotes the linear operator 
$$
\N_J[\varphi](p)= \int_M \varphi(q)\,N_J(p,q)\,dq,
$$
then
$$
\N_J\,\bx_J = \bx_J \,\N_J = I + \sum_{A\subset J}
(-1)^{|A|}S_{J,A} + \sum_{A\cap J = \emptyset} (-1)^{|A|}S_{J,A}.
$$
}

\Subsec{$L^p$-regularity and replacements for maximal hypoellipticity}\label{SSLPregularity} We have
the following regularity results for the operators $\N_J$ and $\K_J$.

\begin{theorem} Let $Q = Q(Z)$ be any quadratic expression
in the vector fields $\{Z_1, \bar Z_1, \ldots, Z_n, \bar Z_n\}$ on
$M$. Then the operator $Q(Z)\N_J$ is a bounded operator on
$L^p(M)$ for $1 < p < \infty$.
\end{theorem}

\Proof  The statement about $Q(Z)\N$ follows from Theorem
\ref{T:5.1.1.1}, the remark following it, and Theorem
\ref{T:7.1.1.5}.
\Endproof\vskip4pt

As we have pointed out, the operator $\K$ does not satisfy all maximal hypoelliptic estimates, and consequently $Q(Z)\K$ for a
general quadratic expression $Q$ in the vector fields $\{Z_1, \bar Z_1, \ldots, Z_n, \bar Z_n\}$ will not in general be bounded on
$L^2(M)$. 
The second main result of this paper gives the appropriate substitute
result.  It is Theorem 2.4.2 stated in Section 2.

\Proof  The proof of the theorem follows from 
\ref {T:QuadraticDerivatives} and Theorem \ref{T:7.1.1.5}, since operators on $\widetilde M$ involving $(T_k-T_l)$ map to the zero operator on $M$ under
 transference.  \hfq

\section{Pseudo-metrics on $M$}\label{S:Geometries}

So far, we have obtained various $L^p$ regularity results for the
relative fundamental solutions to the Kohn-Laplacian on $M$. We
are also interested in describing the nature of the singularities
of the corresponding kernels. As pointed out in the introduction,
we cannot expect that these kernels behave like standard
fractional integration operators on a space of homogeneous type
where there is one distinguished metric. The object of this
section is to study two different metrics or pseudo-metrics on the
hypersurface $M$ given in equation (\ref{E:1.0.2}) that are
relevant to the analysis of the $\dbar_b$-complex. We will then be
able to describe the singularities of our kernels in terms of
these metrics.

\Subsec{The sum of squares metric}\label{SS:Sum of
squares} The vector fields $\{X_1, \ldots, X_{2n}\}$ have the property that
they and all their commutators span the tangent space at each
point. Hence (\cite{NaStWa85}) they define a natural nonisotropic
metric on $M$ which we write $d_\Sigma$. This metric has the
property that balls of radius $\delta$ are essentially ellipsoids
of radius $\delta$ in the directions of the vector fields $\{X_1,
\ldots, x_{2n}\}$, but are much smaller in the missing $T$
direction.

Explicitly, if $p=(z,t)$ and $q=(w,s)$ are two points of $M$, then
\begin{eqnarray*}
d_\Sigma(p,q) &\approx& \sum_{j=1}^n |z_j-w_j| \\&& + \min_j
\Big\{ \mu_j\Big(w_j,\Big\vert t-s +
2\,\IM\Big[\sum_{j=1}^n\sum_{k=1}^{m_j} \frac{1}{k!} \frac{X^k
P_j}{X z_j^k}(w_j)(z_j-w_j)\Big]\Big\vert\Big)\Big\}.
\end{eqnarray*}
The corresponding ball centered at $(w,s)$ of radius $\delta$ is
given by
\begin{eqnarray*}
B_\Sigma\big((w,s),\delta\big) &\approx &\Bigg\{(z,t) \in
M\,\Big\vert\, |z_j-w_j| < \delta \qquad \text{and} \\&&\Big\vert
t-s + 2\, \IM\Big[\sum_{j=1}^n\sum_{k=1}^{m_j} \frac{1}{k!}
\frac{\partial^k P_j}{\partial z_j^k}(w_j)(z_j-w_j)\Big]\Big\vert
< \sum_{j=1}^n \Lambda_j(w_j,\delta)\Bigg\}.
\end{eqnarray*}
The volume of this ball is
\begin{equation}\label{E:Volume1}
\Big\vert B_\Sigma\big((w,s),\delta\big)\Big\vert \approx
\delta^{2n}\,\Big[\sum_{j=1}^n \Lambda_j(w_j,\delta)\Big].
\end{equation}

The appearance of the expression $$ \Big\vert t-s +
2\,\IM\Big[\sum_{j=1}^n\sum_{k=1}^{m_j} \frac{1}{k!}
\frac{\partial^k P_j}{\partial z_j^k}(w_j)(z_j-w_j)\Big]\Big\vert
$$ occasionally makes it difficult to work with this distance.
However, it is always possible to make a biholomorphic change of
variables so that the point $(w,s)$ is moved to the origin, and
all pure $z$ and $\bar z$ derivatives of the defining polynomials
$\{P_j\}$ vanish there.

Suppose that $w = (w_1,\ldots,w_n,w_{n+1}) \in M$. Let
$\Phi^w:\C^{n+1} \to \C^{n+1}$ be the biholomorphic mapping given
by $\Phi^w(z_1,\ldots,z_n,z_{n+1}) = (\zeta_1,
\ldots,\zeta_n,\zeta_{n+1})$ where
\begin{equation}
\begin{split}
\zeta_j &= z_j-w_j ,\qquad 1\leq j \leq n;\\ \zeta_{n+1} &=
z_{n+1}- w_{n+1} -2i\, \sum_{j=1}^n\sum_{k=1}^{m_j}
\frac{1}{k!}\,\frac{\partial^k P_j}{\partial
z_j^k}(w_j)\,(z_j-w_j)^k.
\end{split}
\end{equation}
Also, set
\begin{equation}
P_j^w(\zeta_j) = \sum_{\substack{k \geq 1\\ \ell \geq 1}}
\frac{1}{k!\,\ell!} \,\frac{\partial^{k+\ell}P_j}{\partial
z_j^k\,\partial \bar z_k^\ell}(w_j)\,\zeta_j^k\,\bar\zeta_j^\ell.
\end{equation}
Note that each $P_j^w$ is again a subharmonic, nonharmonic
polynomial of degree $m_j$, and that
\begin{equation}
\frac{\partial^k P_j^w}{\partial \zeta_j^k}(0) = \frac{\partial^k
P_j^w}{\partial \bar\zeta_j^k}(0) = 0, \quad 1 \leq k \leq m_j.
\end{equation}
The mapping $\Phi^w$ maps the point $w$ to the origin of
$\C^{n+1}$ and maps the hypersurface $M$ to the hypersurface
\begin{equation}
M^w = \Big\{(z_1,\ldots,z_n,z_{n+1}) \in \C^{n+1}\,\Big\vert
\IM[z_{n+1}] = \sum_{j=1}^n P_j^w(z_j)\Big\}.
\end{equation}

Write $w_{n+1} = s + i\sum_{j=1}^n P_j(w_j)$. Since we can
identify $M$ and $M^w$ with $\C^n\times \R$, the mapping $\Phi^w$
induces a change of variables on $\C^n \times\R$ given by
$\Phi^a(z_1,\ldots,z_n,t) = (\zeta_1,\ldots,\zeta_n,s)$ where
\begin{equation}
\begin{split}
\zeta_j &= z_j - w_j,\qquad 1\leq j \leq n;\\ s &= t-s +
2\,\IM\Big[\sum_{j=1}^n\sum_{k=1}^{m_j} \frac{1}{k!}
\frac{\partial^k P_j}{\partial z_j^k}(w_j)(z_j-w_j)\Big].
\end{split}
\end{equation}

In this case, we say that $M$ is normalized at the origin, and we
have the much simpler expressions
\begin{equation*}
B_\Sigma(0,\delta) \approx \Big\{(z,t) \in M\,\Big\vert\,
\text{$|z| < \delta$ and $|t| < \sum_{j=1}^n
\Lambda_j(0,\delta)$}\Big\}.
\end{equation*}
The volume of this ball is given by
\begin{equation*}
\big\vert B_\Sigma(0,\delta)\big\vert \approx
\delta^{2n}\Big[\sum_{j=1}^n \Lambda_j(0,\delta)\Big].
\end{equation*}
The corresponding distance from $(z,t)$ to the origin is
\begin{equation*}
d_\Sigma(z,t) = \sum_{j=1}^m |z_j| +
\min_j\Big\{\mu_j(0,|t|)\Big\}.
\end{equation*}

\Subsec{The Szeg\ho\ metric}\label{SS:Szego}  There is a second pseudo-metric $d_S$ on $M$ which also
plays an important role in our analysis. In general it is not equivalent to
the sum of squares metric. The ball of radius $\delta$ is
essentially an ellipsoid of length $\delta$ in the $T$ direction,
and of length $\mu_j(p,\delta)$ in the $z_j$ direction. Thus
unlike the sum of square balls, these balls are not isotropic in
the complex directions $z_1, \ldots,z_n$. It follows from the scaling arguments in \cite{McNeal89} and \cite{NaRoStWa88} that the
Szeg\"o kernel on $M$ behaves like a singular integral operator
relative to the metric $d_S$, and so we call this the Szeg\"o
metric.

\medskip

\begin{definition}
Let $p=(z_1,\ldots,z_n,t)$ and $q=(w_1,\ldots,w_n,s)$ be two
points in $M$. Set
\begin{eqnarray*}\label{E:2.2.1}
d_S(p,q)& = &\sum_{j=1}^n \Lambda_j(w_j,|z_j-w_j|) \\
&&+ \Big\vert
t-s+2\,{\rm Im}\Big[\sum_{j=1}^n\, \sum_{k=1}^{m_j}
\frac{1}{k!}\,\frac{\partial ^k P_j}{\partial
z^k}(w_j)\,(z_j-w_j)^k\Big]\Big\vert.
\end{eqnarray*}
In particular, if $q=(0,0)$ and if the domain is normalized at the
origin, then
\begin{equation*}
d_S(p,0) = \sum_{j=1}^n \Lambda(0,z_j) + |t|.
\end{equation*}
\end{definition}

It is not hard to check that the function $d_S$ has the properties
of a pseudo-metric given in the following proposition.

\medskip

\begin{proposition} \label{P:2.13}The function $d_S$ has
the following properties\/{\rm :}\/
\begin{enumerate}
\item For all $p,q \in M$ we have $d_S(p,q) \geq 0${\rm ,} and $d_S(p,q)
=0$ if and only if $p=q$.
 \item There exists a constant $C$ so
that for all $p,q\in M${\rm ,}
\begin{equation*}
d_S(p,q)\leq C\,d_S(q,p).
\end{equation*}
\item There exists a constant $C$ so that for all $p,q,r\in M${\rm ,}
\begin{equation*}
d_S(p,r) \leq C\big[d_S(p,q) + d_S(q,r)\big].
\end{equation*}
\end{enumerate}
\end{proposition}
\vglue8pt

The balls corresponding to this pseudo-metric are given by
\begin{equation*}
B_S(p,\delta) = \left\{q \in M \,\Big\vert\,d_S(p,q) <
\delta\right\}.
\end{equation*}
The measure of these Szeg\"o balls is then given by
\begin{equation}\label{E:Volume2}
\Big\vert B_S(p,\delta)\Big\vert \approx \delta\,\prod_{j=1}^n
\mu_j(p,\delta).
\end{equation}

\Subsec{Comparison of $d_\Sigma$ and $d_S$}\label{SS:Comparison}  There is also a more intrinsic way
of defining these balls and distances. Recall that $\bar Z_{j} = X_{j}+ i X_{n+j}$ where the $X_{j}$ are
real vector fields. The ball $B_{\Sigma}$ centered at $q$ of radius $\delta$ is essentially the set of points
$p$ to which one can flow from $q$ along piecewise smooth curves tangent to one of the vectors $\{X_1,
\ldots, X_{2n}\}$ for a total time less than $\delta$. The ball $B_{S}$ centered at $q$ of radius $\delta$ is
essentially the set of points $p$ to which one can flow from $q$ along piecewise smooth curves which are
tangent to one of the vectors $\{X_{j},X_{n+j}\}$ for a total time less than $\mu_{j}(q;\delta )$, for $1\leq j
\leq n$.

From this description, or by direct calculation of the two
distances, we obtain the following inclusion of balls and
relationship between distances.

\begin{lemma}\label{T:Ball-incl} Let $0 < \delta$. Then
\begin{align*}
B_{\Sigma}(q;\delta ) &\subset
B_{S}\big(q;\,\max_{j}\{\Lambda_{j}(q;\delta )\}\big),\\
B_{S}(q;\delta ) &\subset
B_{\Sigma}\big(q;\,\max_{j}\{\mu_{j}(q,\delta)\}\big).
\end{align*}
Also
\begin{align}
\min_{j}\Big\{\Lambda_{j}(q;\,d_{\Sigma}(p,q))\Big\} &\leq
d_{S}(p,q) \leq
\max_{j}\Big\{\Lambda_{j}(q;\,d_{\Sigma}(p,q))\Big\},\notag\\
\min_{j}\Big\{\mu_{j}(q;\,d_{S}(p,q))\Big\} &\leq d_{\Sigma}(p,q)
\leq \max_{j}\Big\{\mu_{j}(q;\,d_{S}(p,q))\Big\}.\notag
\end{align}
\end{lemma}

\Proof  Suppose that $p \in B_{\Sigma}(q;\delta)$. We
 can flow from $q$ to $p$ along a
piecewise smooth curve tangent to one of the vectors $\{X_1,
\ldots, X_n\}$ for time less than or equal to $\delta$. Hence we
can flow from $q$ to $p$ along a curve which is tangent to one of
$X_{j}, X_{n+j}$ for a total time at most
$\mu_{j}\big(q;\,\Lambda_{j}(q;\delta )\big)\leq
\mu_{j}\big(q;\,\max_{k}\{\Lambda_{k}(q;\delta )\}\big)$. It
follows that $p \in B_{S}\big(q;\max_{j}\{\Lambda_{j}(q;\delta
)\}\big)$. This proves the first inclusion of balls. The second is
proved in the same way.

The inequalities follow from these inclusions. For example,
suppose $d_{\Sigma}(p,q)\break = \delta$. Then $p\in
B_{\Sigma}(q;\delta)\subset
B_{S}\big(q,\max_{j}\{\Lambda_{j}(q;\delta )\}\big)$. Hence
$$d_{S}(p,q) \leq \max_{j}\{\Lambda_{j}(q;\delta )\} =
\max_{j}\{\Lambda_{j}\big(q;\,d_{\Sigma}[p,q]\big)\}.
$$ This is the
second part of the first inequality. Similarly, if $d_{S}(p,q) =
\delta$, then $p\in B_{S}(q;\delta)\subset
B_{\Sigma}\big(q;\max_{j}\{\mu_{j}(q;\delta )\}\big)$, and so
$d_{\Sigma}(p,q) \leq\max_{j}\{\mu_{j}(q;\delta )\} =
\max_{j}\{\mu_{j}(q;\, d_{S}(p,q))\}$, which is the second part of
the second inequality. The first half of each inequality follows
in the same way.\Endproof\vskip4pt

There is a relationship between the volumes of the balls
$B_\Sigma(p,\delta)$ and $B_S(p,\delta)$. More generally, there is
a relationship between the sizes of  fractional integral
operators with respect to these two metrics. Recall that if $d$ is
a metric, then a singular integral kernel $S(p,q)$ relative to $d$
has size
\begin{equation*}
|S(p,q)| \lesssim \big\vert B\big(p,d(p,q)\big)\big\vert^{-1}
\end{equation*}
and a fractional integral kernel $K(p,q)$ smoothing of order
$\alpha$ has size
\begin{equation*}
|K(p,q)| \lesssim d(p,q)^\alpha\, \big\vert
B\big(p,d(p,q)\big)\big\vert^{-1}.
\end{equation*}

In the case of $d_\Sigma$ and $d_S$ we have

\begin{corollary} Suppose that $\alpha \geq 0$. Then
\begin{equation*}
\frac{\big(d_{\Sigma}(p,q)\big)^{\alpha}}{\big|
B_{\Sigma}\big(q;\,d_{\Sigma}(p,q)\big)\big|} \,\,\lesssim\,\,
\frac{\big(\max_{j}\{\mu_{j}\big(q;\,d_{S}(p,q)\big)\}\big)^{\alpha}}{
\big| B_{S}\big(q;\,d_{S}(p,q)\big)\big|}.\label{E:ball-ineq}
\end{equation*}
\end{corollary}
\vskip8pt

\Proof  We shall use the abbreviations $d_\Sigma(p,q) = d_\Sigma$ and $d_S(p,q) = d_S$. Using the volumes of the balls $B_{S}$ and $B_{\Sigma}$ given in equations (\ref{E:Volume1}) and (\ref{E:Volume2}), the stated inequality is equivalent to 
\begin{align*}
d_{\Sigma}^{\alpha-2n}\,\Big[\sum_{j=1}^n \Lambda_j(p,d_\Sigma)\Big]^{-1} \lesssim
d_{S}^{-1}\,\Big[\prod_{j=1}^n \mu_j(p,d_S)^2\Big]^{-1}\, \Big[ \max_j \big\{ \mu_j(p,d_S) \}\Big]^\alpha.
\end{align*}
However, according to  Lemma \ref{T:Ball-incl}, we have $d_\Sigma
\leq \max_j\{\mu_j(p,d_S)\}$, $d_\Sigma^{-2n} \leq
\left[\prod_{j=1}^n \mu_j(p,d_S)^2\right]^{-1}$ and $d_S \leq
\sum_{j=1}^n \Lambda_j(p,d_\Sigma)$. This completes the proof.
\hfq

\section{Differential inequalities for the relative fundamental solution
$K$}\label{S:Size Estimates}

We now show that the distribution kernel $K_J$ for the relative fundamental
solution $\K_J$ is singular only on the diagonal of $M \times M$,
and we obtain estimates on the size of the kernel and its
derivatives away from the diagonal.

\Subsec{Statement of the main result} Let $d_S$ denote the Szeg\"o metric on $M$ and let $d_\Sigma$
denote the sum of squares metric on $M$. One of the main results stated in Section
2 is Theorem 2.4.4.  It gives the following estimate: 
\begin{eqnarray*}
\Big\vert \Big[\prod_{j=1}^n \partial_j^{\alpha_j}\Big]\,K_{J}(p,q)\Big\vert
&\lesssim&\frac{\left[\sum_{j=1}^n
\mu_j\big(p,d_S(p,q)\big)\right]}{\Big\vert B_S\big(p,d_S(p,q)\big)\Big\vert}^2\,\log\left[2+ \frac{\sum_{j=1}^n \mu_j\big(p,d_S(p,q)\big)}{d_\Sigma(p,q)}\right]
\\&& \times\prod_{j=1}^n \left[\mu_j\big(p,d_S(p,q)\big)^{-1} + d_\Sigma(p,q)^{-1} \right]^{|\alpha_j|}.
\end{eqnarray*}

Before beginning the proof of this, we make several remarks.

\begin{enumerate}
\item At least formally, the distribution kernel $K(z,w,t)$ is given by \end{enumerate}\vskip-14pt
\begin{equation}\label{E:XY}
K(z,w,t) = \int_{\Sigma(t)}\widetilde K(z,w,r)\,d\tilde r.
\end{equation}
\vskip-24pt
\phantom{up}
\begin{enumerate} \item[]
The kernel $\widetilde K(z,w,r)$ has singularities whenever $z_j = w_j = r_j = 0$, and our formal integral (\ref{E:XY}) runs over these nonintegrable singularities. We deal with this difficulty as follows. We have observed that $\widetilde K(z,w,r)$ is the limit as $\varepsilon \to 0$ of kernels $\widetilde K_\varepsilon(z,w,r)$ which are, for $\varepsilon > 0$, smooth, bounded functions on $\C^n\times \C^n \times \R^n$. We write $K(z,w,t)= \lim_{\varepsilon \to 0} K_\varepsilon(z,w,t)$ where 
$$
K_\varepsilon(z,w,t) = \int_{\Sigma(t)}\widetilde K_\varepsilon(z,w,r)\,d\tilde r 
$$
and where this integral now converges. Near points where the kernel $\widetilde K$ becomes singular, we integrate the corresponding $\widetilde K_\varepsilon$ by parts to obtain good estimates which are independent of $\varepsilon$. This will justify the formal calculations. In the discussion below, we will suppress the dependence on $\varepsilon$ with the understanding that all estimates are uniform in $\varepsilon$.

\item[(2)] For simplicity of exposition, we shall deal only with the case $n = 2$. The estimates in this case are
sufficiently complicated, and for larger values of $n$, the arguments require similar computations together
with appropriate induction hypotheses.

\item[(3)] We shall further simplify the notation by making our computations of 
\pagebreak
$K(z,w,t)$ at the point $w =
0$ and we shall assume that our domain is normalized at the origin. As we have seen, there is no loss of
generality in doing this, and various expressions involving distance functions become easier to write.
\end{enumerate}
\vglue-12pt

\Subsec{Szeg\ho\ kernels as derivatives} In this section we establish the estimates used in integration by
parts near points where the integrand in (\ref {E:XY}) has nonintegrable singularities. These singularities are
caused by the presence of distribution kernel  $S_j(z,w,t)$  of the Szeg\"o projection $\S_j$ either onto the
null space of $\bar Z_j$ or $Z_j$ in $L^2(M_j)$:\begin{equation*}
\S_j[f](z,t) = \underset{\C\times\R}{\iint} S_j(z,w,t-s) \, f(w,s) \, dw \,ds.
\end{equation*}
Recall that  $S_j(z,w,t)$ is a distribution on $\C \times \C\times \R$ which is singular only on the set where $z=w$ and $t=0$. If $\partial^\alpha_j$ denotes a derivative of order $|\alpha|$ in the vector fields $Z_j$ and $\bar Z_j$ acting either in the $(z,t)$ or $(w,t)$ variables, then we have the estimate
\begin{equation}\label{EAA}
\Big\vert \partial^\alpha_j\,\partial^k_t S_j(z,w,t-s)\Big\vert \leq C_{\alpha,k}\, d_j\big((z,t),(w,s)\big)^{-2-|\alpha|}
\Lambda_j\big(z, d_j\big((z,t),(w,s)\big)\big)^{-1-k}.
\end{equation}

In transferring from $\widetilde M$ to $M$, we shall be forced to integrate formally  over the nonintegrable
singularities that occur when $(z,t) = (w,s)$. What will save us is the fact that $\partial^\alpha_j\,S_j(z,w,t)$
is essentially a high derivative in $t$ of a bounded function, and so we are able to integrate by parts. To see
this, we recall that the Szeg\"o kernel is given as an integral of the corresponding Bergman kernel. If
$\partial^\alpha_j$ denotes a derivative of order $\alpha$ in the vector fields $Z_j$ and $\bar Z_j$ acting
either in the $(z,t)$ or $(w,s)$ variables, we have
\begin{equation*}
\partial^\alpha_j S_j(z,w,t-s) = \int_0^\infty \partial^\alpha_j B_j(z,w,t-s+ir)\,dr.
\end{equation*}
Moreover, the Bergman kernel $B_j$ satisfies the size estimates
\begin{equation}\label{AD}
\begin{split}
&\Big\vert \partial^\alpha_j \,\partial^k_t\,\partial^\ell_r
B_j(z,w,t-s+ir)\big\vert \\ \leq &C_{\alpha,k}\,
\left[d_j\big((z,t),(w,s)\big)+\mu_j(z,r)\right]^{-2-|\alpha|}
\left[\Lambda_j\big(z,
d_j\big((z,t),(w,s)\big)\big)+r\right]^{-2-k-r}.
\end{split}
\end{equation}
(The relationship between the Szeg\"o and Bergman projections and the estimates (\ref{EAA}) and (\ref{AD}) for the distribution kernels can be found in
\cite{NaRoStWa88}.)

 \medbreak
We shall use the following decomposition when we need to integrate by parts.

\begin{lemma}\label{L:9.3.2} Fix $\delta > 0$. Then for each $\alpha \geq 0$ and each integer
 $m > 2+ \frac{1}{2} \,\alpha${\rm ,} there is a constant $C = C(\alpha,m)${\rm ,} and there are functions
$F_j^{\m}(z,w,t)$ and $G_j^{\m}(z,w,t)$ so that
\begin{equation*}
\partial_j^\alpha S_j(z,w,t) = F_j^\m(z,w,t) + \partial^m_t G_j^\m(z,w,t)
\end{equation*}
where
\begin{align*}
|F_j^\m(z,w,t)| & \leq C\, \frac{1}{\left[d_j((z,t),(w,0))+\mu_j(z,\delta)\right]^{\alpha+2}\, \left[ \Lambda_j (z, d_j(z,t) , (w,0) ) +\delta\right]}\\ |G_j^\m(z,w,t)| &
\leq
C\,\frac{\delta^m}{\left[d_j((z,t),(w,0))+\mu_j(z,\delta)\right]^{\alpha+2} \, \left[\Lambda_j (z , d_j (z , t) , (w , 0)) +\delta\right]}
.
\end{align*}
\end{lemma}

\Proof 
Write
$$
\partial_j^\alpha S_j(z,w,t) = \int_0^\delta \partial_j^\alpha B_j(z,w,t+is)\,ds +
\int_\delta^\infty \partial_j^\alpha B_j(z,w,t+is)\,ds.
$$

Now for any smooth function $\varphi$ on the interval $[0,\delta]$, Taylor's theorem gives for any positive integer $m$
$$
\int_0^\delta \varphi(s)\,ds = \sum_{k=0}^{m-1} \frac{(-1)^k}{(k+1)!}\,\delta^{k+1}\,\varphi^{(k)}(\delta) +
\frac{(-1)^m}{m!} \int_0^\delta s^m\,\varphi^{(m)}(s)\,ds.
$$
We use this with $\varphi(s) = \partial^\alpha_j\,B_j(z,w,t_is)$, and get
\begin{eqnarray*}
\partial_j^\alpha S_j(z,w,t) &= &\sum_{k=0}^{m-1}
\frac{(-1)^k\,\delta^{k+1}}{(k+1)!}\,\partial^k_s\,\partial_j^\alpha B_1(z,t+i\delta)
+ \int_\delta^\infty \partial_j^\alpha B_j(z,w,t+is)\,ds\\&& + \frac{(-1)^m}{m!} \int_0^\delta
s^m\,\partial^m_s\,\partial_j^\alpha B_j(z,w,t+is)\,ds.
\end{eqnarray*}
Set 
$$
F_j^\m(z,t) = \sum_{j=0}^{m-1} \frac{(-1)^j\,\delta^{j+1}}{(j+1)!}\,\partial^j_s\,X^\alpha B_j(z,w,t+i\delta) + \int_\delta^\infty \partial_j^\alpha B_j(z,w,t+is)\,ds.
$$
Because of the rate of decrease of $|\partial_j^\alpha B_j(z,t+is)|$ as $s \to \infty$, we have
\begin{multline*}
 \Big\vert\int_\delta^\infty \partial_j^\alpha B_j(z,w,t+is)\,ds \Big\vert
\\ 
\lesssim \frac{1}{\left[d_j((z,t),(w,0))+\mu_j(z,\delta)\right]^{\alpha+2}\, \left[\Lambda_j(z,d_j(z,t),(w,0)) +\delta\right]}.
\end{multline*}
Also, the estimates in (\ref{AD}) show that 
$$
\big\vert
\partial^k_s\,\partial_j^\alpha B_j(z,w,t+i\delta) \big\vert \lesssim 
\frac{1}{\left[d_j((z,t),(w,0))+\mu_j(z,\delta)\right]^{\alpha+2}\, \left[\Lambda_j(z,d_j(z,t),(w,0)) +\delta\right]}.
$$
 Thus we have established the correct estimate for $\displaystyle \big\vert F_j^\m(z,t)\big\vert$.

On the other hand, since $B(z,w,\zeta)$ is holomorphic in $\zeta$, we have
$$
\frac{(-1)^m}{m!} \int_0^\delta s^m\,\partial^m_s\,\partial_j^\alpha  B_j(z,w,t+is)\,ds = \partial^m_t
\Big[\frac{(-i)^m}{m!} \int_0^d s^m\,\partial_j^\alpha B_j(z,w,t+is)\,ds\Big].
$$
Set 
$$
G_j^\m(z,t) = \frac{(-i)^m}{m!} \int_0^\delta s^m\,X^\alpha B_j(z,w,t+is)\,ds.
$$
We have
\begin{align*}
\Big\vert \int_0^\delta s^m\,&\partial_j^\alpha B_j(z,w, t+is)\,ds \Big\vert
\\
&\leq \int_0^\delta \frac{s^{m}}{\left[d_j((z,t), (w,0))+\mu_j(z,s)\right]^{\alpha+2}\left[\Lambda_j (z,d_j((z,t),(w,0)) + s\right]^2}ds
\\ &\leq
\int_0^{\mu_j(z,\delta)} \frac{\Lambda_j(z,s)^{m}\,\Lambda_j'(z,s)\,ds}
{\left[d_j((z,t),(w,0))+s\right]^{\alpha+2}\,
\left[\Lambda_j(z,d_j((z,t),(w,0)) + \Lambda_j(z,s)\right]^2}.
\end{align*}
Now $\Lambda_j(z,s)^{m-2}\,\Lambda_j'(z,s) = s^{2m-3}\varphi(s)$ where $\varphi$ is increasing. Thus if $2m-3 > \alpha + 1$, we can estimate this last integral by the length of the interval times the value of 
the integrand at the right-hand endpoint. Thus we get
\begin{multline*}
\Big\vert \int_0^\delta s^m X^\alpha B_j(z,w, t+is)\,ds \Big\vert
\\ \lesssim \frac{\delta^m}{\left[d_j((z,t),(w,0))+ \mu_j(z,\delta)\right]^{\alpha+2}\,
\left[\Lambda_j(z,d_j((z,t),(w,0))) + \delta\right]}.
 \end{multline*}
This gives us the required estimate for $G_j^\m(z,t)$, and completes the proof of Lemma \ref{L:9.3.2}. 
\Endproof\vskip4pt

We remark that if $\alpha = 0$ and $m = 2$, then the term
corresponding to $G_2(z,w,t)$ is not bounded, but involves a
logarithm. In fact, suppose $A, B > 0$. Then
\begin{align*}
\int_0^B \frac{ds}{\mu_j(z, A+s)^2} &=
\int_A^{A+B}\frac{ds}{\mu_j(z,s)^2}\\ &=
\int_{\mu_j(z,A)}^{\mu_j(z,A+B)}\frac{ \Lambda_j'(z,t)}{t^2}\,dt\\
&\lesssim \frac{(A+B)}{\mu_j(z,A+B)^2}
\,\log\left[\frac{\mu_j(z,A+B)}{\mu_j(z,A)}\right].
\end{align*}

\Subsec{Proof of Theorem {\rm \ref{TH4}}} When $n = 2$, according to Theorem \ref{T:MainEstimate},
we need to consider four integrals:
\begin{eqnarray*}
K(z_1,z_2,0,0,t) &=& \int_{-\infty}^{+\infty}
K_0(z_1,r,z_2,t-r)\,dr \\& &+ \int_{-\infty}^{+\infty}
S_1(z_1,r)\,K_1(z_1,r,z_2,t-r)\,dr \\&& +
\int_{-\infty}^{+\infty}
S_2(z_2,t-r)\,K_2(z_1,r,z_2,t-r)\,dr\\&&  +
\int_{-\infty}^{+\infty}S_1(z_1,r)\,S_2(z_2,t-r)\,
 K_\infty(z_1,r,z_2,t-r)\,dr.
 \end{eqnarray*}
The first of these is the easiest to analyze since the integral
does not run across any singularities unless $z_1 = z_2 = t = 0$.
The last of the four is the hardest since there are two possible
places where we will need to integrate by parts: when $z_1 = r =
0$ and when $z_2 = t-r = 0$. We shall do the computations only in
the first and last cases, since the other two integrals are then
very similar.

 \Subsubsec{The estimate for $K_0$}\label{SS:K_0}  We shall use the abbreviations
$$d_\Sigma((z_1,r,z_2,t-r),(0,0,0,0)\big) = d_\Sigma(z_1,r,z_2,t-r)
= d_\Sigma.$$
  Recall that on $\widetilde M$,
\begin{equation*}
d_\Sigma((z_1,r,z_2,t-r), (0,0,0,0)\big) \approx |z_1| + \mu_1(0,r)
+ |z_2| + \mu_2(0,t-r).
\end{equation*}
Also, recall that on $M$
\begin{equation*}
d_\Sigma(z_1,z_2,t) \approx |z_1|+|z_2|+
\min\big\{[\mu_1(0,t),\mu_2(0,t)\big\}.
\end{equation*}
Thus to estimate $K_0$ we have
\begin{multline*}
X^\alpha \Big[\int_{-\infty}^{+\infty} K_0(z_1, r,z_2,t-r)\,dr
\Big]
\\ \lesssim
\int_{-\infty}^{+\infty}
d_\Sigma(z_1,r,z_2,t-r)^{-2-|\alpha|}\big[\Lambda_1(0,d_\Sigma)+
\Lambda_2(0,d_\Sigma)\big]^{-2}\,dr.
\end{multline*}
We split this integral into three parts: the first is where $|r|
\leq \frac{1}{2}|t|$ in which case $|t-r| \approx |t|$; the second
is where $|t-r| \leq \frac{1}{2}|t|$ in which case $|r| \approx
|t|$, and finally the complement where $|t-r| \approx |r|$.

For the first integral we get the estimate
\begin{eqnarray*}
&&(|z_1|+|z_2| +
\mu_2(t))^{-2-|\alpha|}\,\big[\Lambda_1(0,|z_1|+|z_2| +
\mu_2(t))+\Lambda_2(0,|z_1|+|z_2| + \mu_2(t))\big]^{-1}\\  &&\hskip1.5in
\lesssim \frac{d_\Sigma(z_1,z_2,t)^{2-|\alpha|}}{\Big\vert
B_\Sigma\big(0,d_\Sigma(z_1,z_2,t)\big)\Big\vert}
\end{eqnarray*}
which is an estimate of the correct sort. For the second integral,
we interchange the roles of $r$ and $t-r$, and obtain the same
estimate. Finally for the third integral, we can make the estimate
$$
\int_{|r| \geq \frac{1}{2}|t|}
d_\Sigma(z_1,r,z_2,r)^{-2-|\alpha|}\,\Big[\Lambda_1(0,d_\Sigma(z_1,r,z_2,r))
+ \Lambda_2(0,d_\Sigma(z_1,r,z_2,r))\big]^{-1}.
$$
Because of the decay of the integral as $|r| \to \infty$, we can
estimate this by $|t|$ times the value of the integrand at $r =
|t|$. This again gives the correct sort of estimate, and completes
the analysis of the integral involving $K_0$.

\Subsubsec{The estimate for $K_\infty$}\label{SS:Kinfty} 
 We now turn to the hardest estimate involving $K_\infty$. Let $X$
denote a derivative using either $Z_1$ or $\bar Z_1$, and let $Y$
denote a derivative using either $Z_2$ or $\bar Z_2$. We want to
estimate integrals of the form
\begin{equation}\label{E:10.1.1}
\int_{-\infty}^{+\infty}
X^{\alpha_1}S_1(z_1,r)\,Y^{\beta_1}S_2(z_2,t-r)\,X^{\alpha_2}\,Y^{\beta_2}
K_{\infty}(z_1,r,z_2,t-r)\,dr
\end{equation}
by
\begin{equation}\label{E:13.2.2z}
\begin{split}
\frac{\left(\mu_1(d_S) +
\mu_2(d_S)\right)^2}{d_S\,\mu_1(d_S)^2\,\mu_2(d_S)^2}&\,
\log\left[2 + \frac{d_S}{\min\Big\{\Lambda_1(d_\Sigma),
\,\Lambda_2(d_\Sigma)\Big\}}\right]\\
&\cdot\left[d_\Sigma^{-1}+\mu_1(d_S)^{-1}\right]^{\alpha_1}
\,\left[d_\Sigma^{-1} + \mu_2(d_s)^{-1}\right]^{\beta_1}
\,d_\Sigma^{-(\alpha_2+\beta_2)}\,.
\end{split}
\end{equation}

The integrand is dominated by
\begin{eqnarray*}
&&\left[|z_1|+
\mu_1(0,|r|)\right]^{-2-\alpha_1}\left[\Lambda_1(0,|z_1|)
+|r|\right]^{-1} \left[|z_2|
+\mu_2(0,|t-r|)\right]^{-2-\beta_1}\\&&\quad\cdot\left[\Lambda_2(0,|z_2|)
+ |t-r|\right]^{-1}\, \left[|z_1|+ \mu_1(0,|r|)+|z_2|
+\mu_2(0,|t-r|)\right]^{2-(\alpha_2+\beta_2)}
\end{eqnarray*}
which in turn is dominated by a sum
\begin{eqnarray*}
&&\left[|z_1|+
\mu_1(0,|r|)\right]^{-\alpha_1}\left[\Lambda_1(0,|z_1|)
+|r|\right]^{-1} \left[|z_2|
+\mu_2(0,|t-r|)\right]^{-2-\beta_1}\\&&\quad\cdot\left[\Lambda_2(0,|z_2|)
+ |t-r|\right]^{-1}\, \left[|z_1|+ \mu_1(0,|r|)+|z_2|
+\mu_2(0,|t-r|)\right]^{-(\alpha_2+\beta_2)}\\
&&\qquad+ \left[|z_1|+
\mu_1(0,|r|)\right]^{-2-\alpha_1}\left[\Lambda_1(0,|z_1|)
+|r|\right]^{-1} \left[|z_2|
+\mu_2(0,|t-r|)\right]^{-\beta_1}\\&&\quad\cdot\left[\Lambda_2(0,|z_2|)
+ |t-r|\right]^{-1}\, \left[|z_1|+ \mu_1(0,|r|)+|z_2|
+\mu_2(0,|t-r|)\right]^{-(\alpha_2+\beta_2)}.
\end{eqnarray*}
In estimating the integral (\ref{E:10.1.1}), the integrand can
possibly become infinite when $ r = 0$ (if $|z_1| = 0)$, and when
$r = t$ (if $|z_2| = 0)$.

We shall assume
\begin{equation}
\Lambda_1(0,|z_1|) \leq \Lambda_2(0,|z_2|).
\end{equation}
Also, we have
\begin{equation}
d_S(z_1,z_2,t) \approx \Lambda_1(0,|z_1|) + \Lambda_2(0,|z_2|)
+|t|,
\end{equation}
and
\begin{equation}
d_\Sigma(z_1,z_2,t) \approx
|z_1|+|z_2|+\min\big\{\mu_1(0,|t|),\,\mu_2(0,|t|)\big\}.
\end{equation}

\demo{Case $1$}  Suppose $\Lambda_1(0,|z_1|) \leq
\Lambda_2(0,|z_2|) \leq 2|t|$. In this case $d_S(z_1,z_2,t)\break
\approx |t|$.
  First consider the integral over the region where $|r|
\geq \frac{1}{2}|t|$ and $|t-r| \geq \frac{1}{2}|t|$. In this
region, $|t-r| \approx |r|$, and we can dominate this part of the
integral (\ref{E:10.1.1}) by
\begin{eqnarray*}
&&\hskip-24pt  \underset{|r| \geq \frac{1}{2}|t|}{\int} \left[|z_1| +
\mu_1(r)\right]^{-2-\alpha_1}\,\left[\Lambda_1(z_1) +
|r|\right]^{-1}\,\left[|z_2| +
\mu_2(r)\right]^{-2-\beta_1}\\&&\qquad\cdot\left[\Lambda_2(z_2) +
|r|\right]^{-1}\,\left[|z_1|+|z_2|+\mu_1(r) +
\mu_2(r)\right]^{2-(\alpha_2+\beta_2)}\,dr
\\
&&\quad \lesssim |t|\, \left[|z_1| +
\mu_1(t)\right]^{-2-\alpha_1}\,\left[\Lambda_1(z_1) +
|t|\right]^{-1}\,\left[|z_2| +
\mu_2(t)\right]^{-2-\beta_1}\\&&\qquad \cdot \left[\Lambda_2(z_2) +
|t|\right]^{-1}\,\left[|z_1|+|z_2|+\mu_1(t) +
\mu_2(t)\right]^{2-(\alpha_2+\beta_2)}\\ &&\quad \lesssim
\mu_1(d_S)^{-2-\alpha_1}\,\mu_2(d_S)^{-2-\beta_1}
\,d_S^{-1}\,\left[\mu_1(d_S) +
\mu_2(d_S)\right]^{2-(\alpha_2+\beta_2)}
\end{eqnarray*}
since $|t| \approx d_S$. But since $\mu_1(d_S) + \mu_2(d_S) \geq
d_\Sigma$, this part of the integral (\ref{E:10.1.1}) is dominated
by a constant times \ref{E:13.2.2z}.

Next, we split the integral for $|r| \leq \frac{1}{2}|t|$ into the
region $A = \{\Lambda_1(0,\mu_2(0,|t|))\break \leq |r| \leq
\frac{1}{2}|t|\}$ and the region $B= \{|r| \leq \min\big\{
\Lambda_1(0,\mu_2(0,|t|)),\,\frac{1}{2}|t|\big\}\}$. In both of
these regions we have $|t-r|\approx |t|$.

The integral of (\ref{E:10.1.1}) over the region $A$ is dominated
by
\begin{eqnarray*}
&&\hskip-14pt \int_A \left[|z_1|+
\mu_1(0,|r|)\right]^{-\alpha_1}\left[\Lambda_1(0,|z_1|)
+|r|\right]^{-1} \left[|z_2|
+\mu_2(0,|t|)\right]^{-2-\beta_1}\\&&\qquad\cdot \left[\Lambda_2(0,|z_2|)
+ |t|\right]^{-1}\, \left[|z_1|+ \mu_1(0,|r|)+|z_2|
+\mu_2(0,|t|)\right]^{-(\alpha_2+\beta_2)}\,dr
\\&&\qquad + \int_A\left[|z_1|+
\mu_1(0,|r|)\right]^{-2-\alpha_1}\left[\Lambda_1(0,|z_1|)
+|r|\right]^{-1} \left[|z_2|
+\mu_2(0,|t|)\right]^{-\beta_1}\\&&\qquad\cdot\left[\Lambda_2(0,|z_2|) +
|t|\right]^{-1}\, \left[|z_1|+ \mu_1(0,|r|)+|z_2|
+\mu_2(0,|t|)\right]^{-(\alpha_2+\beta_2)}\,dr.
\end{eqnarray*}
Note that $d_\Sigma(z_1,z_2,t) \lesssim |z_1| + \mu_2(0,|t|)$.
Hence the first term in this sum is dominated by a constant times
\begin{equation*}
d_\Sigma^{-\alpha_1} \mu_2\big(d_S\big)^{-2-\beta_1}\,d_S^{-1}\,
d_\Sigma^{-(\alpha_2+\beta_2)}\, \log\left[2 +
\frac{d_S}{\Lambda_1(d_\Sigma)}\right].
\end{equation*}
The second integral in this sum is dominated by
\begin{equation*}
d_\Sigma^{-2-\alpha_1} \mu_2\big(d_S\big)^{-\beta_1}\,d_S^{-1}
\,d_\Sigma^{-(\alpha_2+\beta_2)}.
\end{equation*}
since the integral converges at infinity. Thus since
$\min\{\mu_1(d_S),\mu_2(d_S)\} \lesssim d_\Sigma$, both terms are
dominated by \ref{E:13.2.2z}.

In the integral of (\ref{E:10.1.1}) over the region $B$, we
integrate the term\break $X^\alpha_1\,S(z_1,r)$ by parts. This means
replacing $X^\alpha_1\,S(z_1,r)$ by a sum $F^\m_1 +
\partial^m_r\,G_1^\m$. In estimating the term with $F^\m_1$, we
can still replace $|t-r|$ by~$|t|$. Let $r^* =
\min\{\Lambda_1(\mu_2(t)),\frac{1}{2}|t|\}$. Then this term is
dominated by a constant multiple of
\begin{eqnarray*}
&&\underset{|r|\leq r^*}{\int} [|z_1| +\mu_1(r^*)]^{-2-\alpha_1}\,
\left[\Lambda_1(z_1)+|r^*|\right]^{-1}
\left[|z_2|+\mu_t(t)\right]^{-2-\beta_1}\left[\Lambda_2(z_2)
+|t|\right]^{-1}\\ &&\quad \cdot\left[|z_1|+|z_2| +
\mu_2(t)\right]^{-(\alpha_2+\beta_2)}\,
\left[|z|_1+|z_2|+\mu_1(t)+\mu_2(t)\right]^2\,dr.
\end{eqnarray*}
If $r^* = \Lambda_1(\mu_2(t))$, we get the estimate
$$
d_\Sigma^{-2-\alpha_1}\,\mu_2(d_S)^{-2-\beta_1}\,d_S^{-1}
\,d_\Sigma^{-(\alpha_2+\beta_2)}.
$$
If $r^* = \frac{1}{2}|t|$ we get the estimate
$$
d_S^{-2-\alpha_1}\,\mu_2(d_S)^{-2-\beta_1}\,d_S^{-1}
\,d_\Sigma^{-(\alpha_2+\beta_2)}.
$$
Either case is  still dominated by \ref{E:13.2.2z}.
 
Now consider the actual integration of $G^\m_1$ by parts. Each
time we integrate, we introduce a factor of $(r^*)^{+1}$ from the
integration and also introduce $\min\big\{ (\Lambda_2(z_2)
+|t|)^{-1},\,\Lambda_1(|z_1|+|z_2| + \mu_2(t))^{-1}\big\}$ from
the differentiation. In either case, the ratio is bounded, and we
are reduced to the same estimate as before.
 
Finally, we need to consider the integral of (\ref{E:10.1.1}) over
the region where $|t-r| \leq \frac{1}{2}|t|$. But now we can
interchange $r$ and $t-r$, and the arguments are the same as those
for the region $|r| \leq \frac{1}{2}|t|$ with the roles of the
subscripts $1$ and $2$ interchanged. We again get the estimate
\ref{E:13.2.2z}. This completes the analysis of Case 1.

\demo{Case $2$} Suppose $2|t| \leq
\Lambda_2(0,|z_2|)$. In this case $$d_S \approx \Lambda_2(z_2).$$
Also, $\min\{\mu_1(t),\,\mu_2(t)\} \leq \mu_2(t) \leq \mu_2(2t)
\leq |z_2|$. Thus $$d_\Sigma(z_1,z_2,t) \approx |z_1|+|z_2|.$$

First consider the integral over the region where $|r| \geq
\Lambda_2(z_2)$. Because of our hypothesis on $|t|$, it follows in
this region that $|t-r|\approx |r|$, and so we need to estimate
\begin{eqnarray*}&&
\underset{|r| \geq \Lambda_2(z_2)}{\int}  \left[|z_1|+
\mu_1(0,|r|)\right]^{-2-\alpha_1}\left[\Lambda_1(0,|z_1|)
+|r|\right]^{-1} \left[|z_2|
+\mu_2(0,|r|)\right]^{-2-\beta_1}\\&&\qquad\qquad\cdot \left[\Lambda_2(0,|z_2|) +
|r|\right]^{-1}\, \left[|z_1|+ \mu_1(0,|r|)+|z_2|
+\mu_2(0,|r|)\right]^{2-(\alpha_2+\beta_2)}\,dr.
\end{eqnarray*}
We break up the term $ \left[|z_1|+ \mu_1(0,|r|)+|z_2|
+\mu_2(0,|r|)\right]^2$ into two parts, and integrate the
resulting integrals separately. Because of the rate of decay of
the integrand, we can estimate the integrals by $\Lambda_2(z_2)$
times the value of the integrand at the left endpoint. We obtain
the estimate
$$
\mu_1(d_S)^{-2-\alpha_1} \,\mu_2(d_S)^{-\beta_2} \,d_S^{-1}
\,d_\Sigma^{-(\alpha_2+\beta_2)} + \mu_1(d_S)^{-\alpha_1}
\,\mu_2(d_S)^{-2-\beta_2} \,d_S^{-1}
\,d_\Sigma^{-(\alpha_2+\beta_2)}.
$$
This gives us the desired estimate \ref{E:13.2.2z}.

Thus it remains to integrate over the region $|r| \leq
\Lambda_2(z_2)$ when $2|t| \leq \Lambda_2(z_2)$ and hence $d_S =
\Lambda_2(z_2)$ and $d_\Sigma = |z_1|+|z_2|$. There are four
subcases to deal with.

\demo{Case {\rm 2a}}  Assume that $|z_1|\leq |z_2|$ so that
$d_\Sigma \approx |z_2|$ and assume that $\Lambda_2(z_2)\leq
\Lambda_1(z_2)$. In this case we integrate by parts, integrating
$X^\alpha_1 S_1(z,r)$, and differentiating the other parts of the
integrand over the whole interval $|r| \leq \Lambda_2(z_2)$. We
must estimate
\begin{equation}\label{E:13.1.6y}
\begin{split}
\int_0^{\Lambda_2(z_2)} \Big[F^\m_1(z_1,r,d) &+
\partial^m_r[G^\m_1(z_1,r,d)]\Big] \\&\cdot Y^{\beta_1}S_2(z_2,t-r)
X^{\alpha_2}Y^{\beta_2}K_\infty(z_1,r,z_2,t-r)\,dr.
\end{split}
\end{equation}
We choose $d = \Lambda_2(z_2)= d_S$, and write
\begin{equation*}
H(z_1,r,z_2,t-r) = Y^{\beta_1}S_2(z_2,t-r)
X^{\alpha_2}Y^{\beta_2}K_\infty(z_1,r,z_2,t-r).
\end{equation*}
We can make the following estimates when $|r| \leq
\Lambda_2(z_2)$:
\begin{eqnarray*}
\big\vert F^\m_1(z_1,r,d)\big\vert & \leq& [|z_1| +
\mu_1(d_S)]^{-2-\alpha_1}\,[\Lambda_1(z_1)+ d_S]^{-1},\\ \big\vert
G^\m_1(z_1,r,t)\big\vert & \leq &d_S^m\, [|z_1| +
\mu_1(d_S)]^{-2-\alpha_1}\,[\Lambda_1(z_1)+ d_S]^{-1},\\ \big\vert
H(z_1,r,z_2,t-r) \big\vert &\lesssim&
\Big[|z_2|^{-\beta_1}\,\Lambda_2(z_2)^{-1}
\,[|z_1|+|z_2|]^{-(\alpha_2+\beta_2)}\\ && +
|z_2|^{-2-\beta_1}\,\Lambda_2(z_2)^{-1}
[|z_1|+|z_2|]^{-(\alpha_2+\beta_2)}\,[|z_1|+ \mu_1(d_S)]^2\Big],\\
 \big\vert \partial^m_r
H(z_1,r,z_2,t-r)\big\vert & \lesssim& \Big[
|z_2|^{-\beta_1}\,\Lambda_2(z_2)^{-1}
\,[|z_1|+|z_2|]^{-(\alpha_2+\beta_2)}\\ && +
|z_2|^{-2-\beta_1}\,\Lambda_2(z_2)^{-1}
[|z_1|+|z_2|]^{-(\alpha_2+\beta_2)}\,[|z_1|+ \mu_1(d_S)]^2\Big]
\\&&\times
\max\{\Lambda_1(z_2)^{-1},\Lambda_2(z_2)^{-1}\}^m.
\end{eqnarray*}

Now the term in the integral (\ref{E:13.1.6y}) that does not
require integration by parts (the part involving $F^\m_1$) can be
estimated by
$$
\mu_1(d_S)^{-2-\alpha_1}\,\mu_2(d_S)^{-\beta_1}\,d_S^{-1}\,
d_\Sigma^{-(\alpha_2+\beta_2)} + \mu_1(d_S)^{-\alpha_1}\,
\mu_2(d_S)^{-2-\beta_1}\,d_S^{-1}\, d_\Sigma^{-(\alpha_2+\beta_2)}
$$
  which is dominated by the estimate \ref{E:13.2.2z}.
After we integrate by parts, the term involving $G^\m_1$ can be
estimated in the same way since $\Lambda_2(z_2) \leq
\Lambda_1(z_2)$ and hence
$$
\max\big\{\Lambda_1(z_2)^{-1},\,\Lambda_2(z_2)^{-1}\big\} =
\Lambda_2(z_2)^{-1} = d_S^{-1}.
$$

\demo{Case {\rm 2b}}   Assume that $|z_1|\leq |z_2|$ so
that $d_\Sigma \approx |z_2|$, and assume that $\Lambda_1(z_2)
\leq \Lambda_2(z_2) = d_S$. In this case, we integrate by parts
only over the interval $|r| \leq \Lambda_1(z_2)$, and just make
elementary estimates on the part where $\Lambda_1(z_2) \leq |r|
\leq \Lambda_2(z_2)$.

For the first integral, we must estimate
\begin{multline}\label{E:13.1.6x}
\int_0^{\Lambda_1(z_2)} \Big[F^\m_1(z_1,r,d)  +
\partial^m_r[G^\m_1(z_1,r,d)]\Big] \\ \cdot  Y^{\beta_1}S_2(z_2,t-r)
X^{\alpha_2}Y^{\beta_2}K_\infty(z_1,r,z_2,t-r)\,dr.
\end{multline}
We choose $d = \Lambda_1(z_2)$, and again write
$$
H(z_1,r,z_2,t-r) = Y^{\beta_1}S_2(z_2,t-r)
X^{\alpha_2}Y^{\beta_2}K_\infty(z_1,r,z_2,t-r).
$$
We can make the following estimates when $|r| \leq
\Lambda_1(z_2)$:
\begin{eqnarray*}
\big\vert F^\m_1(z_1,r,d)\big\vert & \leq& [|z_1| +
|z_2|]^{-2-\alpha_1}\,[\Lambda_1(z_1)+ \Lambda_1(z_2)]^{-1}\\
&\approx &d_\Sigma^{-2-\alpha_1}\,\Lambda_1(d_\Sigma)^{-1},\\
\big\vert G^\m_1(z_1,r,t)\big\vert & \leq& (\Lambda_1(z_2))^m\,
[|z_1| + |z_2|]^{-2-\alpha_1}\,[\Lambda_1(z_1)+
\Lambda_1(z_2)]^{-1}\\ &\approx&
(\Lambda_1(z_2))^m\,d_\Sigma^{-2-\alpha_1}\,\Lambda_1(d_\Sigma)^{-1},\\
 \big\vert
H(z_1,r,z_2,t-r) \big\vert &\lesssim& \Lambda_2(z_2)^{-1}
\,[|z_1|+|z_2|]^{-(\beta_1+\alpha_2+\beta_2)}\\ &\approx&
d_S^{-1}\,d_\Sigma^{-(\beta_1+\alpha_2+\beta_2)},\\
 \big\vert \partial^m_r
H(z_1,r,z_2,t-r)\big\vert & \lesssim& {\bf }\Lambda_2(z_2)^{-1}
\,[|z_1|+|z_2|]^{-(\beta_1+\alpha_2+\beta_2)} 
\\
&&\cdot \max\{\Lambda_1(z_2)^{-1},\Lambda_2(z_2)^{-1}\}^m\\ &\approx&{\bf }
d_S^{-1}\,d_\Sigma^{-(\beta_1+\alpha_2+\beta_2)}\,\Lambda_1(z_2)^{-m}.
\end{eqnarray*}
The part of the integral in (\ref{E:13.1.6x}) is estimated by
$$
d_S^{-1}\,d_\Sigma^{-(2+\alpha_1+\beta_1+\alpha_2+\beta_2)}
$$
which is dominated by the estimate \ref{E:13.2.2z} since  
$\min\big\{\mu_1(d_S),\,\mu_2(d_S)\big\}\lesssim d_\Sigma$.

For the second integral, we must estimate
\begin{eqnarray*}
&&\int_{\Lambda_1(z_2)}^{d_S} d_\Sigma^{-\alpha_1} [\Lambda_1(z_1)
+ r]^{-1}\,[\Lambda_2(z_2)]^{-1}\,
d_\Sigma^{-(2+\beta_1+\alpha_2+\beta_2)}\,dr\\ &&\qquad+ 
\int_{\Lambda_1(z_2)}^{d_S}  [|z_1|+\mu_1(r)]^{-\alpha_1}
[\Lambda_1(z_1) + r]^{-1}\,[\Lambda_2(z_2)]^{-1}\,
d_\Sigma^{-(2+\beta_1+\alpha_2+\beta_2)}\,dr
.
\end{eqnarray*}
Both of  these integrals are dominated by
$$
d_S^{-1}\,d_\Sigma^{-(2+\alpha_1+\beta_1+\alpha_2+\beta_2)}
\,\log\left[2+\frac{d_S}{\Lambda_1(d_\Sigma)}\right],
$$
although the logarithm term only appears in the first integral.
This again gives the estimate \ref{E:13.2.2z}.

\demo{Cases {\rm 2c} and {\rm 2d}} Here we assume that $|z_2|
\leq |z_1|$. This time we integrate the term $Y^\beta_1
S_2(z_2,t-r)$, and everything goes as in Cases 2a and 2b with the
roles of the subscripts $1$ and $2$ interchanged.
 
This completes the analysis of the integral involving $K_\infty$,
and consequently completes the proof of the estimates in Theorem
\ref{TH4} when $n = 2$.

\section{H\"older regularity for $\K$}\label{S:Further estimates}

We now turn to the study of the smoothing properties of the relative fundamental solutions $\K_J$ on the scale of (isotropic) H\"older spaces, and the proof of Theorem \ref{TH3}.

\Subsec{$L^1$ modulus of continuity of $K$} We first observe that the presence of the logarithm term in the
estimate for $K$ does not affect the $L^1$ norm of $K$ over small
balls.

\begin{proposition}\label{P:11.1.1} There is a constant $C$ such
that for $\delta > 0${\rm ,}
$$
\int_{d_S(p,q)< \delta} \big\vert K(p,q)\big\vert dq \leq
C\Big[\sum_{j=1}^n \mu_j(p,\delta)\Big]^2.
$$
\end{proposition}

\Proof  Without loss of generality, we can assume that $p$
is the origin, and that the domains $M_j$ are normalized at the
origin. We abbreviate $d_S\big((z,t),(0,0)\big)$ by $d_S(z,t)$,
and the same for $d_\Sigma$. We write
\begin{multline*}
 \int_{d_S\big(z,t) < \delta} \big\vert
K\big((z,t),(0,0)\big)\big\vert\,dz\,dt
\\  =
\sum_{j=0}^\infty \,\sum_{k=0}^\infty
\,\underset{\substack{d_S(z,t)
\approx\, 2^{-j}\,\delta\\
\inf_{j}\{\Lambda_j(0,d_\Sigma(z,t))\}\approx\, 2^{-k}d_S(z,t)}}{\int}
K\big((z,t),(0,0)\big)\,dz\,dt.
\end{multline*}
The size of the integrand of the $(j,k)^{\rm th}$ integral is
dominated by a constant times
$$
k\,\frac{\sum_{j=1}^n \mu_j(0, 2^{-j}\,\delta)^2}{\big\vert
B_S(0,2^{-j}\,\delta)\big\vert} 
$$
where the factor $k$ comes from the logarithm term. On the other
hand, the region of integration is given by
\begin{align*}
\sum_{j=1}^n \Lambda_j(0,|z_j|) + |t| &\approx 2^{-j}\,\delta,\\
\inf_j\Big\{\Lambda_j\Big(0,\sum|z_j| +
\min\{\mu_j(0,t)\}\Big)\Big\}
 &\approx
2^{-k-j}\delta.
\end{align*}
The volume of this is dominated by
$$
2^{-\varepsilon\,k}\,\big\vert B_S(0,2^{-j}\,\delta)\big\vert.
$$
Hence the integral is dominated by
\vglue12pt
\hfill $
\displaystyle{\sum_{j=0}^\infty \left[\sum_{\ell=1}^n\mu_\ell(0,2^{-j}\delta)^2
\,\sum_{k=0}^\infty k\,2^{-\varepsilon\,k}\right]\approx
\left[\sum_{\ell = 1}^n\mu_\ell(0,\delta) \right]^2.}
$
\Endproof\vskip12pt

Next, we have the following estimates on the $L^1$ modulus of
continuity of the kernel $K$.

\begin{proposition}\label{P:7.1.2r} Let $h \in \C^2\times \R$ be a vector with
Euclidean length $|h|$. Then
$$
\int_{B_S(p,1)}\big\vert K(p+h,q) - K(p,q)\big\vert dq \lesssim
\sum_{j=1}^n \mu^\sh_j(p,|h|)^2
$$
where
$$
\mu_j^\sh(p,\delta)^2 = \delta\,\int_{\delta}^1
\mu_j(p,t)^2\,\frac{dt}{t^2}\geq \mu_j(p,\delta)^2.
$$
\end{proposition}

\Proof  We split the integral into two parts. The first is
where $d_S(p,q) \leq 10 |h|$, and here we simply estimate the
$L^1$ norm of $K$ as in Proposition \ref{P:11.1.1}. In the second
integral, where $d_S(p,q) \geq 10 |h|$, we use the estimate
\begin{equation*}
\begin{split}
\big\vert K(p+h,q) - K(p,q)\big\vert &\leq h \big\vert\partial_t
K(p,q)\big\vert\\ & \leq
h\,\frac{\sum_{j=1}^n\mu_j(p,d_S(p,q))^2}{d_S(p,q)\,V(p,q)}.
\end{split}
\end{equation*}
 A similar estimate then gives the
desired result, and completes the proof.
\Endproof\vskip4pt  

If we use a second difference, we can improve on Proposition
\ref{P:7.1.2r}.

\begin{proposition}\label{P:11.2.2} Let $h \in \C^2\times \R$.
Then
$$
\int_{B_S(p,1)}\Big\vert K(p+h,q) +K(p-h,q) -2 K(p,q)\Big\vert
\,dq \lesssim \sum_{j=1}^n\mu_j(p,|h|)^2.
$$
\end{proposition}
\vglue-12pt

\Subsec{Applications to smoothness}   As an immediate consequence of
the above we have the following result:

\begin{theorem} Suppose that $f$ is a bounded function supported
in\break $B_S(p_0,10)$. Then for $p \in B_S(p_0,1)${\rm ,}
\begin{equation*}
\begin{split}
\big\vert\K[f](p+h) - \K[f](p)\big\vert &\leq
C\,\left[\sum_{j=1}^n\mu_j^\sh(p,|h|)^2\right],\\
\big\vert\K[f](p+h)+\K[f](p-h) -2 \K[f](p)\big\vert &\leq
C\,\left[\sum_{j=1}^n\mu_j(p,|h|)^2 \right].
\end{split}
\end{equation*}
In particular{\rm ,} if $m = \max\{m_1,m_2\}>2$ 
is the maximum type{\rm ,}
then
\begin{equation*}
\big\vert\K[f](p+h) - \K[f](p)\big\vert \leq C\,|h|^{\frac{2}{m}}.
\end{equation*}
In general\/{\rm , (}\/when $m \geq 2${\rm ),}  
\begin{equation*}
\big\vert\K[f](p+h)+\K[f](p-h) -2 \K[f](p)\big\vert  \leq C \,
| h |^{\frac{2}{m}}.
\end{equation*} 
\end{theorem}

Indeed, $\mu_j ( p , \delta )^2 \lesssim \delta^{\frac{2}{m_j}}$ by
Definition 3.1.2 of $\mu_j$.  Next, since 
$( \mu_j^{\#} )^2 = \delta
\displaystyle{\int^{1}_{\delta}} \, \mu_j ( t )^2 \, \frac{dt}{t}$,
we see that $( \mu_j^{\#} ( \delta ) )^2 \lesssim \delta^{\frac{2}{m_j}}$, 
when $m_j > 2$.  The desired conclusions are therefore established.

\section{Examples}

In this section we provide examples that show where our regularity results for $\bx_b$ are optimal. We study the same decoupled boundary as in Section \ref{SSSpecialCase}, and let
\begin{equation*}
M = \left\{(z_1,z_2,z_3) \in \C^3\,\Big\vert\,\IM[z_3] = |z_1|^n + |z_2|^m\right\}
\end{equation*}
where $m$ and $n$ are positive even integers. Assume that $m \geq n$. We consider the operator $\bx_b$ as it acts on $(0,1)$-forms of the form $f\,d\bar z_2$. Then
\begin{equation*}
\bx_b[f\,d\bar z_2] = (\bx_1+\bx_2)[f]\,d\bar z_2
\end{equation*} 
where $\bx_1 = -Z_1\,\bar Z_1$ and $\bx_2 = - \bar Z_2\,Z_2$. Since we are not in degree zero or two, the operator $\bx_b$ has no null space in $L^2(M)$, and we have constructed an operator $\K$ so that $\K\,\bx_b = \bx_b\,\K = I$.

If we identify $M$ with $\C^2\times \R$ then as in Section \ref{SSSpecialCase} we have
$$
\bar{Z}_1 = \frac{\partial}{\partial \bar{z}_1} - i\,\frac{n}{2} |z_1 |^{n-2} \, z_1 \frac{\partial}{\partial t}, \qquad \bar{Z}_2 = \frac{\partial}{\partial \bar{z}_2} - i\,\frac{m}{2} | z_2|^{m-2} z_2 \frac{\partial}{\partial t}
.
$$
For $\gamma \in \R$ set
$$
F_\gamma(z_1,z_2,t) = (t + i|z_1|^n+i|z_2|^m)^\gamma.
$$
It is easy to check that
\begin{align*}
\bx_1[F_\gamma] &= 0,\\
\bx_2[F_\gamma] &= -im(m-2)|z_2|^{m-4}z_2^2\,F_{\gamma-1},\\
\overline\bx_1[F_\gamma] &= -in(n-2) |z_1|^{n-4}z_1^2\,F_{\gamma-1},\\
\overline\bx_2[F_\gamma] & = 0,
\end{align*}
and hence that
\begin{equation*}
|\bx_b[F_\gamma]| \approx |z_2|^{m-2}\,|F_{\gamma-1}|.
\end{equation*}

We now show that the H\"older regularity established in Theorem \ref{TH3} is optimal. Observe that $\bx_b[F_\gamma]$ is bounded near the origin on $M$ if and only if $m-2 + m(\gamma-1) \geq 0$, that is, if and only if $\gamma \geq 2/m$. However, if $F_\gamma$ satisfies an isotropic H\"older condition of order $\alpha$, then so does its restriction to the line $z_1 = z_2 = 0$. But $F_\gamma(0,0,t) = t^\gamma$ which does not satisfies a H\"older condition of order greater than $\gamma$. Thus on $M$, if a function $g$ is bounded, the equation $\bx_b[u] = g$ has solutions which do not satisfy H\"older conditions of any order greater than $2/m = 2/\max\{m,n\}$.

Next, Theorem \ref{TH2} (b) provides a replacement for estimates giving maximal hypoellipticity. 
Suppose that $\bx_b[u] \in L^p(M)$. It need not follow that $\overline \bx_1[u] \in L^p(M)$, but we show
that if $|z_1|^{n-2}\,|B(z_1,z_2,t)|\lesssim |z_2|^{m-2}$, it does follow that $B\,\bx_1[u] \in L^p(M)$.   Let
$\chi \in \mathcal C^\infty_0(M)$ be with $\chi \equiv 1$ near the origin. Then  
\begin{align*}
B\,\overline \bx_1 [\chi F_\gamma] &\approx B\,|F_{\gamma_1}|\,|z_1|^{n-2} \\&= \Big(B\,\frac{|z_1|^{n-2}}{|z_2|^{m-2}}\Big)\,\Big[|z_2|^{n-2}\,|F_{\gamma-1}|\Big]\approx  \Big(B\,\frac{|z_1|^{n-2}}{|z_2|^{m-2}}\Big)\,\big|\bx_b[\chi F_\gamma]\big|
.\end{align*}
Since $F_\gamma \in L^p_{\rm loc}(M)$ if and only if
 $\displaystyle \gamma\,p > -\big(1 + \frac{2}{m} + \frac{2}{n}\big)$, we see that our condition is
essentially optimal.

\end{document}